\documentclass[10pt]{article}

\usepackage{geometry} 
\usepackage{graphicx}

\usepackage{latexsym}
\usepackage{amsmath,amssymb,mathtools,amsfonts}
\usepackage{amsthm}
\usepackage{float}
\usepackage{url}
\usepackage{enumerate}
\usepackage{subcaption}
\usepackage{caption}

\usepackage{setspace}

\theoremstyle{definition}

\newtheorem{thm}{Theorem}
\newcommand{\argmin}{\mathop{\rm arg~min}\limits}

\makeatletter
\newcommand{\figcaption}[1]{\def\@captype{figure}\caption{#1}}
\newcommand{\tblcaption}[1]{\def\@captype{table}\caption{#1}}
\makeatother

\title{\huge{Approximation and application of minimizing movements for surface PDE}}

\author{Elliott Ginder \\{\emph{School of Interdisciplinary Mathematical Sciences, Meiji University}}
\vspace{10pt}
\\ Karel Svadlenka \\{\emph{Department of Mathematics, Tokyo Metropolitan University}}
\vspace{10pt}
\\Takuma Muramatsu 
\\{\emph{Graduate School of Advanced Mathematical Sciences, Meiji University}}}

\date{} 

\begin{document}

\maketitle
\thispagestyle{empty}
\addtocounter{page}{-1}
\newpage

\newgeometry{top=2cm}
\tableofcontents
\clearpage
\restoregeometry
\setstretch{1.2}
\section{Introduction}\label{序論}
We extend the applicability of minimizing movements (MM) to approximating solutions of surface partial differential equations (SPDE) and apply the technique to approximate mean curvature flow (MCF)\cite{mcf} and hyperbolic MCF (HMCF)\cite{hmcf} on surfaces.
The MBO algorithm \cite{mcf} and the HMBO algorithm \cite{hmcf}, both based on the level set method, are well-known approximation methods for such curvature flows.
We recall that the MBO algorithm is an approximation method for mean curvature flow and is based on solving the heat equation.
On the other hand, the HMBO algorithm is an approximation method for the hyperbolic mean curvature flow and involves solving the wave equation.
It has been shown that the MBO algorithm and HMBO algorithm can approximate curvature flow under area preservation constraints, as well as in the multiphase setting \cite{hmcf, multi_volume}.
Minimizing movements \cite{newbook} are used to realize the area conservation condition, and the signed distance vector field \cite{hmcf} is used for calculations involving multiphase regions.

On the other hand, for curvature flow on curved surfaces, approximation methods for the mean curvature flow were presented in \cite{mbo_surf,mc_surf}.
The authors show that the Closest point method (CPM) \cite{cpm} can be used to approximate mean curvature flow on surfaces.
This is done by extending the values of functions defined on a surface to the ambient space of the surface. In turn, the CPM enables the approximation of surface gradients and other differential quantities by making use of the surround Euclidean space. However, no previous studies have treated approximation of solutions to  curvature flow under area preservation involving interfaces in the multiphase and  curved surface setting.

This study develops approximation methods for mean curvature flow and hyperbolic mean curvature flow under multiphase area preservation conditions for interfaces moving on curved surfaces.

Our approach is to extend MM to the case of surface PDE and to use their framework to apply the MBO and HMBO algorithms. Similar to \cite{cpm}, our generalizations make use the CPM which we combine with the surface-constrained signed distance vector field \cite{aokishibata}.

The outline of this paper is as follows. In section \ref{目的}, we describe the research background and our objectives. Our objectives are based on conventional approximation methods such as the MBO algorithm and HMBO algorithm. To achieve our goals, approximation methods for constrained partial differential equations on surfaces are required. We therefore demonstrate that partial differential equations on surfaces can be computed using the Closest Point Method and introduce surface-type minimizing movements to handle partial differential equations with constraints. In section \ref{曲面上偏微分方程式の数値計算}, we discuss computational techniques related to partial differential equations on surfaces. In particular, we create an approximation method by combining the Closest Point Method with minimizing movements and perform numerical error analyses for heat and wave equations defined on surfaces. In section \ref{曲面上界面運動の数値計算}, we discuss our approximation method for mean curvature flow on surfaces and hyperbolic mean curvature flow on surfaces. This includes an explanation of the signed distance vector field, which is required to handle multiphase domains, and the area preservation condition that is achieved through the use of minimizing movements. In section \ref{数値計算結果と考察}, we discuss mean curvature flow and hyperbolic mean curvature flow on surfaces and describe the method for enforcing area preservation conditions in multiphase environments. We summarize the contents of this paper and discuss future challenges in section \ref{まとめ}.

\section{Background}\label{目的}
In this section, we will briefly explain our goals, and the mathematical frameworks used in our research. In section \ref{sub目的}, we will touch upon our research objectives. Then, in section \ref{cpm目的}, we introduce the Closest Point Method (CPM), and in section \ref{MM} we will introduce the minimizing movements (MM). Again remark that, by combining these methods, we obtain an approximation method for partial differential equations with constraints on curved surfaces.
\subsection{Objectives}\label{sub目的}
As stated in Section \ref{序論}, the goal of this research is to create an approximation method for interfacial motion on curved surfaces. Here, we will explain two representative examples of interfacial motions on curved surfaces: mean curvature flow \cite{mcf} and hyperbolic mean curvature flow \cite{hmcf}. The approximation method for mean curvature flow that we emply is known as the MBO algorithm, which alternates between solving the heat equation and constructing level set functions \cite{mcf}. The numerical solution method for hyperbolic mean curvature flow is known as the HMBO algorithm, which alternates between solving the wave equation and constructing level set functions \cite{hmcf}. Examples of more complex interfacial motions involving area preservation and extendtions to multiphase regions are illustrated. Here, multiphase regions refer to regions where the domain is divided into three or more regions by the interface. We use the fact that the area preserving condition can be realized by imposing a constraint on the partial differential equations used in the MBO and HMBO algorithms \cite{multi_volume, s2}. We again remark that the case of multiphase regions is possible to treat by using the signed distance vector field \cite{hmcf}.

Based on the above considerations, the purpose of this research is as follows:
\begin{itemize}
\item establish an approximation method for partial differential equations with constraints on curved surfaces
\item extend the MBO and HMBO algorithma to curved surfaces and realize numerical compute mean curvature flow and hyperbolic mean curvature flow on surfaces
\item generalize our famework to treat the above problems in the multiphase setting.
\end{itemize}
\subsection{Closest Point Method}\label{cpm目的}
When numerically solving partial differential equations on surfaces, approximation of surface gradients (SG) on the surface are necessary. 
For a smooth surface $S$ embedded in $n$-dimensional Euclidean space, the SG of a function $u$ on the surface $S$ is given by
\begin{align}\label{SG}
\nabla_S u=\nabla u-\boldsymbol{n}(\boldsymbol{n}\cdot\nabla u),
\end{align}
where $\boldsymbol{n}$ is the unit normal vector of the surface $S$, and $\nabla$ is the usual gradient in the Euclidean space.

\begin{figure}[H]
    \begin{center}
        \fbox{\includegraphics[bb=0cm 0cm 8cm 7cm,scale=0.7]{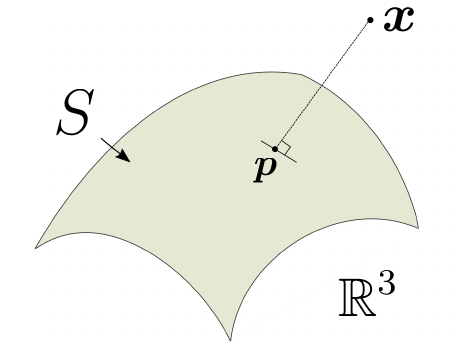}}
        \caption{Example of the closest point.\\( $\boldsymbol{p}=C_S(\boldsymbol{x}),\quad\boldsymbol{x}\in\mathbb{R}^3$ )}\label{ex:cpm}
    \end{center}
\end{figure}

In the CPM, an approximation of the SG on the surface is obtained by smoothly extending the values of the function defined on the surface in the direction of the surface normal vector \cite{cpm}. 
This is enabled by constructing a function that provides the closest point on the surface $S$ to any external point of the ambient space. 
For any point $\boldsymbol{x}$ in $n$-dimensional space, the function $C_S$ that gives the closest point on the surface $S$ is defined as follows:
\begin{align}\label{CPF}
    C_{S}(\boldsymbol{x})=\argmin_{\boldsymbol{y}\in S}|| \boldsymbol{x}-\boldsymbol{y}||
\end{align}
For example, if $S$ is the unit circle ($n=2$), then $C_S$ is given by:
$$
    C_S(x,y)=\left(\frac{x}{\sqrt{x^2+y^2}},\frac{y}{\sqrt{x^2+y^2}}\right),\quad (x,y)\in\mathbb{R}^2
$$
If $S$ is the unit sphere ($n=3$), then $C_S$ is given by:
$$
    C_S(x,y,z)=\left(\frac{x}{\sqrt{x^2+y^2+z^2}},\frac{y}{\sqrt{x^2+y^2+z^2}},\frac{z}{\sqrt{x^2+y^2+z^2}}\right),\quad (x,y,z)\in\mathbb{R}^3
$$
Note that, in the case of the unit circle or sphere, the closest point to the origin is not uniquely determined. 
In this case, since the distance between the origin and the surface is constant, CPM assigns an arbitrary point on the surface as the closest point to the origin. 
Figure \ref{ex:cpm} shows the relationship between a point $\boldsymbol{x}$ in three-dimensional space and its closest point $\boldsymbol{p}$ on the surface $S$. 
The closest point $\boldsymbol{p}$ is given by the closest point function $C_S$ as $\boldsymbol{p}=C_S(\boldsymbol{x})$.
Note that, if a surface $S$ is represented by a parameterization, it is easy to find the closest point. In such a case, we just need to solve the optimization problem \eqref{CPF}. On the other hand, if the surface $S$ is not parameterized but is represented by a point cloud or a triangulated surface, then a bit more ingenuity is required. For example, when the surface $S$ is discretized by a point cloud, the closest point on the surface $S$ may involve constructing implicit surfaces defined by distance functions, or require one to accept a certain level of loss of uniquess. 
The following theorems forms the basis for approximations using the CPM \cite{cpm}.

\begin{thm}\label{na0}
    Let $S \subset \mathbb{R}^3$ be a smooth surface. 
    Let $u: \mathbb{R}^3 \rightarrow \mathbb{R}$ be an arbitrary smooth function that has a constant value in the normal direction of the surface $S$ near the surface. 
    Then, on the surface $S$,
    \begin{align}
        \nabla u= \nabla_S u
    \end{align}
    holds \cite{cpm}.
    Here, $\nabla_S$ designates the surface gradient on the surface $S$.
\end{thm}
\begin{thm}\label{na1}
    Let $v$ be an arbitrary smooth vector field in $\mathbb{R}^3$ that is tangent to the surface $S$ or to a surface that is at a fixed distance from the surface $S$.
     Then, on the surface $S$,
    \begin{align}
        \nabla\cdot v=\nabla_S\cdot v
    \end{align}
    holds.
\end{thm}
If $u$ is a function defined on a surface $S$, then the function $u(C_S(\boldsymbol{x}))$ given by the closest point function $C_S$ is a function that has a constant value in the direction of the normal vector of the surface.
Therefore, by Theorem \ref{na0}, we have:
$$
    \nabla u(C_S(\boldsymbol{x})) = \nabla_S u(\boldsymbol{x}),\quad \boldsymbol{x}\in S
$$
Furthermore, since $\nabla u(C_S(\boldsymbol{x}))$ is always tangent to a surface that is at a short and fixed distance from the surface $S$, Theorem \ref{na1} implies that:
$$
    \nabla\cdot\nabla u(C_S(\boldsymbol{x})) = \nabla_S\cdot\nabla_S u(\boldsymbol{x}),\quad \boldsymbol{x}\in S
$$
The operator $\nabla_S\cdot\nabla_S$ on the right-hand side of the above equation is intrinsic to the surface $S$ and is denoted by $\Delta_S$ (the Laplace-Beltrami operator).  
From the above theorem, it can be seen that by using the extension given by $C_S$, approximation methods used in Euclidean space can be applied to approximate differential operators such as $\nabla_S$ and $\Delta_S$.

In \cite{cpm}, examples of solving partial differential equations on surfaces using CPM and numerical methods are presented. 

As mentioned earlier, in numerical calculations of curvature flow using MBO and HMBO algorithms, an approximation method for surface partial differential equations with constraints is required to achieve the area-preservation condition.
One effective approximation method for constrained partial differential equations is the method of minimizing movements.
In the next section \ref{MM}, we will the extention of minimizing movements to the case of surface PDE.

\subsection{Minimizing movements}\label{MM}
Here we will explain the method of minimizing movements (MM), also known as Discrete Morse Flow \cite{newbook}. 
Minimizing movements are a method for approximating the gradient flow of an energy functional
$$
    \mathcal{E}(u)=\int_\Omega L(\nabla u(\boldsymbol{x}),u(\boldsymbol{x}),\boldsymbol{x}) d\boldsymbol{x}
$$
by iteratively minimizing functionals of the form
$$
    \mathcal{F}_n(u)=\int_{\Omega}\frac{|u-u_{n-1}|^2}{2h}d\boldsymbol{x}+\mathcal{E}(u)
$$
within a suitable function space and where $h>0$ is a suitable time step. 
Here, $\Omega$ is a region with given boundary conditions, and $u_n$ is an approximation of $u$ at time $t=nh$. 
The Euler-Lagrange equation of each functional $\mathcal{F}_n$ represents an approximation of the gradient flow of the energy functional. 
By changing $\mathcal{F}_n(u)$, various approximations of solutions to partial differential equation can be obtained as the Euler-Lagrange equation of $\mathcal{F}_n$. 
For example, for a given $\alpha > 0$, if we set,
\begin{align}\label{mmheat_flat}
    \mathcal{F}_n(u)=\int_{\Omega} \frac{|u-u_{n-1}|^2}{2h}+\alpha\frac{|\nabla u|^2}{2}d\boldsymbol{x}
\end{align}
we obtain an approximation of a heat equation, and if we set
\begin{align}\label{mmwave_flat}
    \mathcal{F}_n( u)
    = \int_{\Omega}\frac{| u-2 u_{n-1}+ u_{n-2}|^2}{2h^2}
    +\alpha\frac{|\nabla  u|^2}{2}\hspace{5pt}d\boldsymbol{x}
\end{align}
we obtain an approximation of the wave equation.
Minimizing movements are based on energy minimizations, so it is possible to naturally handle constrained partial differential equations by adding penalty terms to the functional. 

In the next section we will create an approximation methods for constrained partial differential equations on surfaces by combining minimizing movements with the CPM.

\section{Numerical calculation of PDEs on surfaces}\label{曲面上偏微分方程式の数値計算}
Here we will provide an overview of an approximation method that combines CPM and MM, and explain its algorithm.
We also perform a numerical convergence analysis for the surface heat and wave equations using our method.
\subsection{Approximation of PDEs on surfaces by CPM.}\label{CPMによる曲面上偏微分方程式の近似}
We will explain the method of discretization in time when approximating solutions to the heat and wave equations on a curved surfaces using CPM.
To this end, let $S$ be a closed smooth surface without boundary in $\mathbb{R}^3$.
We consider the following surface heat equation \eqref{spde}:
\begin{align}
    \begin{cases}\label{spde}
        u^S_{t}(t,\boldsymbol{x})=\alpha\Delta_S u^S(t,\boldsymbol{x}),\hspace{25pt} & \boldsymbol{x}\in S, \quad t>0 \\
	u^S(0,\boldsymbol{x})=f(\boldsymbol{x}),\hspace{25pt}                        & \boldsymbol{x}\in S
    \end{cases}
\end{align}
and surface wave equation \eqref{spde2}: 
\begin{align}
    \begin{cases}\label{spde2}
        u^S_{tt}(t,\boldsymbol{x})=\alpha\Delta_S u^S(t,\boldsymbol{x}), \hspace{25pt} & \boldsymbol{x}\in S,\quad t>0 \\
        u^S_{t}(0,\boldsymbol{x})=V_0(\boldsymbol{x}), \hspace{25pt}                   & \boldsymbol{x}\in S           \\
        u^S(0,\boldsymbol{x})=f(\boldsymbol{x}),\hspace{25pt}                          & \boldsymbol{x}\in S
    \end{cases}
\end{align}
Here, $\alpha>0$ is a constant, $f(\boldsymbol{x})$ is the initial condition, $V_0$ is the initial velocity, and $\Delta_S$ is the Laplace-Beltrami operator on the surface $S$.
We remark that boundary conditions are not included in equations \eqref{spde} and \eqref{spde2} because the surface $S$ is without boundary.
For a given time step size $h>0$, we approximate the time derivative in \eqref{spde} using a forward difference, and in \eqref{spde2} we use a centered difference approximation with respect to time.
By defining $u^S_n$ to be an approximation of $u^S$ at time $nh$ where $n=0,1,\cdots$, we obtain:
\begin{align}\label{heat_cpm}
    \begin{cases}
        u^S_{n+1}(\boldsymbol{x})=u^S_{n}(\boldsymbol{x})+h\alpha \Delta_S u_n^S(\boldsymbol{x}), \hspace{15pt} \\
        u_{0}^S(\boldsymbol{x})=f(\boldsymbol{x}), \hspace{15pt}
    \end{cases}\quad \boldsymbol{x}\in S
\end{align}
which is a approximation scheme for equation \eqref{spde}.
Similarly, the result for equation \eqref{spde2} is given by:
\begin{align}\label{wave_cpm}
    \begin{cases}
        u^S_{n+1}(\boldsymbol{x})=2u^S_{n}(\boldsymbol{x})-u^S_{n-1}(\boldsymbol{x})+h^2\alpha \Delta_S u_n^S(\boldsymbol{x}), \hspace{15pt} \\
        u^S_{-1}(\boldsymbol{x})=u^S_0(\boldsymbol{x})-hV_0(\boldsymbol{x}), \hspace{15pt}                                                   \\
        u_{0}^S(\boldsymbol{x})=f(\boldsymbol{x}), \hspace{15pt}
    \end{cases}\quad \boldsymbol{x}\in S
\end{align}

Since $\Delta_S$ is included in the right-hand side of equations \eqref{heat_cpm} and \eqref{wave_cpm}, they are difficult to compute in the general setting.
However, in the CPM, the following equations \eqref{heat_cpm_ext} and \eqref{wave_cpm_ext} are used to compute the solutions in the space $\Omega$ surrounding the surface $S$.
Here, $u_n$ is a function value defined on $\Omega$ at time $nh$.
\begin{align}\label{heat_cpm_ext}
    \begin{cases}
        u_{n+1}(\boldsymbol{x})=u_{n}(C_S(\boldsymbol{x}))+h\alpha \Delta u_n(C_S(\boldsymbol{x})), \\
        u_0(\boldsymbol{x})=f(C_S(\boldsymbol{x})),
    \end{cases}\boldsymbol{x}\in \Omega
\end{align}
\begin{align}\label{wave_cpm_ext}
    \begin{cases}
        u_{n+1}(\boldsymbol{x})=2u_{n}(C_S(\boldsymbol{x}))-u_{n-1}(C_S(\boldsymbol{x}))+h^2\alpha \Delta u_n(C_S(\boldsymbol{x})), \\
        u_{-1}(\boldsymbol{x})=u_0(C_S(\boldsymbol{x}))-hV_0(C_S(\boldsymbol{x})),                                                  \\
        u_0(\boldsymbol{x})=f(C_S(\boldsymbol{x})),
    \end{cases}\boldsymbol{x}\in \Omega
\end{align}
Here, $C_S$ is defined by equation \eqref{CPF}, and $\Delta = \nabla \cdot \nabla$.
Since equations \eqref{heat_cpm_ext} and \eqref{wave_cpm_ext} do not contain $\Delta_S$, it is possible to apply standard numerical approximation techniques in the surrounding Euclidean space calculate surface gradient quantities.
Also, since the surface is given by a point cloud, interpolation can be used to define the numerical solution restricted to the surface $S$ or at any other point in the domain $\Omega$.

Although explicit methods were used to discretize the time derivatives in equations \eqref{spde} and \eqref{spde2}, implicit methods can also be used \cite{Implicitcpm}.
The combination of CPM and MM for the calculation of equations \eqref{spde} and \eqref{spde2} are described in detail in Section \ref{CPMとMMの組み合わせ}.

\subsection{Combination of CPM and MM}\label{CPMとMMの組み合わせ}
As mentioned in Section \ref{目的}, when performing calculations for curvature flow with an area preservation constraint via the MBO or HMBO algorithms, an approximation method for the constrained partial differential equation is necessary.
Here we explain the approximation method for the constrained partial differential equations on surfaces by combining CPM and MM.

We will explain our method for applying CPM to minimizing movement for the surface heat equation \eqref{spde}, and the surface wave equation \eqref{spde2}.
As described in Section \ref{CPMによる曲面上偏微分方程式の近似}, applying CPM yields the approximations for the surface heat equation \eqref{spde} and the surface wave equation \eqref{spde2}, given by equations \eqref{heat_cpm_ext} and \eqref{wave_cpm_ext}, respectively.
As a numerical method, utilizing the method of minimizing movements requires one to approximate functional values. In particular, for $n=0,1,\cdots$, using a time step size $h>0$ and a constant $\alpha>0$, the following functional values are required and can be approximated, for example, by means of the finite element method:
\begin{align}\label{mmheat}
    \mathcal{F}_{n+1}(u)
    = \int_{\Omega}\frac{| u(\boldsymbol{x})- u_{n}(C_{S}(\boldsymbol{x}))|^2}{2h}
    +\alpha\frac{|\nabla  u(\boldsymbol{x})|^2}{2}d\boldsymbol{x}
\end{align}
\begin{align}\label{mmwave}
    \mathcal{F}_{n+1}(  u)
    = \int_{\Omega}\frac{|  u(\boldsymbol{x})-2  u_{n}(C_{S}(\boldsymbol{x}))+  u_{n-1}(C_{S}(\boldsymbol{x}))|^2}{2h^2}
    +\alpha\frac{|\nabla   u(\boldsymbol{x})|^2}{2}d\boldsymbol{x}
\end{align}
Here, $\Omega$ is a sufficiently large region that covers the surface $S$, and $u_{n}$ minimizes functional $\mathcal{F}_{n}$.

In the following, we will show that the Euler-Lagrange equations for equations \eqref{mmheat} and \eqref{mmwave} lead to the implicitly discretized equations using the CPM method for partial differential equations on surfaces.

Let $\phi$ be an arbitrary function from $C_0^\infty(\Omega)$ and $\epsilon$ be a real number.
We compute the first variation of equation \eqref{mmheat} as follows:
\begin{align}\label{mmheat1}
    \left.\frac{d}{d\epsilon}\mathcal{F}_{n+1}(u+\epsilon\phi)\right|_{\epsilon=0}=0.
\end{align}
The first variation is:
\begin{align}
    \frac{d}{d\epsilon}\mathcal{F}_{n+1}(u+\epsilon\phi)
     & =\frac{d}{d\epsilon}\int_{\Omega}\frac{| (u+\epsilon\phi)- u_{n}(C_{S})|^2}{2h}
    +\alpha\frac{|\nabla  (u+\epsilon\phi)|^2}{2}d\boldsymbol{x}\notag                 \\
     & =\int_{\Omega}\frac{(u+\epsilon\phi)- u_{n}(C_{S})}{h}\phi
    +\alpha{\nabla  (u+\epsilon\phi)}\cdot\nabla\phi d\boldsymbol{x}\label{mmheat2}.
\end{align}
Substituting $\epsilon=0$ into equation \eqref{mmheat2}, we obtain:
\begin{align}
    \begin{split}
        \left.\frac{d}{d\epsilon}\mathcal{F}_{n+1}(u)\right|_{\epsilon=0}
        & =\int_{\Omega}\frac{u- u_{n}(C_{S})}{h}\phi
        +\alpha{\nabla u}\cdot\nabla\phi d\boldsymbol{x}\label{mmheat3} \\
        & =\int_{\Omega}\left( \frac{u- u_{n}(C_{S})}{h}
        -\alpha\Delta{u}\right)\phi d\boldsymbol{x}+\alpha\int_{\partial\Omega}\frac{\partial u}{\partial \boldsymbol{\nu}}\phi dS,
    \end{split}
\end{align}
where $\partial\Omega$ is the boundary of $\Omega$, and $\partial u/\partial \boldsymbol{\nu}$ is the outer normal derivative of $u$ on $\partial\Omega$.
Since $\phi$ is an arbitrary $C_0^\infty(\Omega)$ function, the boundary integral in  \eqref{mmheat3} is zero, and we have:
\begin{align*}
    \left.\frac{d}{d\epsilon}\mathcal{F}_{n+1}(u)\right|_{\epsilon=0}
     & =\int_{\Omega}\left( \frac{u- u_{n}(C_{S})}{h}
    -\alpha\Delta{u}\right)\phi d\boldsymbol{x}
\end{align*}
A weak form of the Euler-Lagrange equation \eqref{mmheat1} is therefore:
\begin{align*}
    \int_{\Omega}\left( \frac{u- u_{n}(C_{S})}{h}
    -\alpha\Delta{u}\right)\phi d\boldsymbol{x}=0
\end{align*}
Since $\phi$ is arbitary, the fundamental lemma of the calculus of variations applies to obtain:
\begin{align*}
    \frac{u- u_{n}(C_{S})}{h}\notag
    -\alpha\Delta{u}=0 
\end{align*}
which, written as an approximation scheme states:
\begin{align}
    u=u_{n}(C_{S})+h\alpha\Delta{u}\label{ELheat}
\end{align}

Equation \eqref{ELheat} is an implicit form of the time-discrete surface heat equation \eqref{heat_cpm_ext} obtained by using the CPM.

The functional \eqref{mmwave} can be treated in the same fashion. We obtain:
\begin{align}
    \frac{d}{d\epsilon}\mathcal{F}_{n+1}(u+\epsilon\phi)
     & =\frac{d}{d\epsilon}\int_{\Omega}\frac{| (u+\epsilon\phi)- 2u_{n}(C_{S})+u_{n-1}(C_{S})|^2}{2h^2}
    +\alpha\frac{|\nabla  (u+\epsilon\phi)|^2}{2}d\boldsymbol{x}\notag                                   \\
     & =\int_{\Omega}\frac{(u+\epsilon\phi)- 2u_{n}(C_{S})+u_{n-1}(C_{S})}{h^2}\phi
    +\alpha{\nabla  (u+\epsilon\phi)}\cdot\nabla\phi d\boldsymbol{x}\label{mmwave2}.
\end{align}
Setting $\epsilon=0$ in \eqref{mmwave2}, we have:
\begin{align}
    \begin{split}
        \left.\frac{d}{d\epsilon}\mathcal{F}_{n+1}(u)\right|_{\epsilon=0}
        & =\int_{\Omega}\frac{u- 2u_{n}(C_{S})+u_{n-1}(C_{S})}{h^2}\phi
        +\alpha{\nabla u}\cdot\nabla\phi d\boldsymbol{x}\label{mmwave3} \\
        & =\int_{\Omega}\left( \frac{u- 2u_{n}(C_{S})+u_{n-1}(C_{S})}{h^2}
        -\alpha\Delta{u}\right)\phi d\boldsymbol{x}+\alpha\int_{\partial\Omega}\frac{\partial u}{\partial \boldsymbol{\nu}}\phi dS.
    \end{split}
\end{align}
As before, $\partial u/\partial \boldsymbol{\nu}$ is the derivative of $u$ in the direction of the outer normal vector $\boldsymbol{\nu}$ on $\partial\Omega$.
Since $\phi$ is an arbitrary in $C_0^\infty(\Omega)$, we obtain the following equation:
\begin{align*}
    \left.\frac{d}{d\epsilon}\mathcal{F}_{n+1}(u)\right|_{\epsilon=0}
     & =\int_{\Omega}\left( \frac{u- 2u_{n}(C_{S})+u_{n-1}(C_{S})}{h^2}
    -\alpha\Delta{u}\right)\phi d\boldsymbol{x}.
\end{align*}
It follows that
\begin{align*}
    \frac{u- 2u_{n}(C_{S})+u_{n-1}(C_{S})}{h^2}\notag
    -\alpha\Delta{u}=0,
\end{align*}
which can be written \eqref{ELwave}:
\begin{align}
    u=2u_{n}(C_{S})-u_{n-1}(C_{S})+h^2\alpha\Delta{u}\label{ELwave}.
\end{align}
Equation \eqref{ELwave} is an implicit approximation of the time-discretized surface wave equation \eqref{spde2} using the CPM (compare to Equation \eqref{wave_cpm_ext}).

Having shown that the minimizing schemes above produce approximate solutions to the surface PDE \eqref{spde} and \eqref{spde2} we now turn to discussing related numerical considerations..
Next, we introduce the computational algorithms for implementing the CPM and MM.
\subsection{Computational methods for the heat and wave equations on surfaces}\label{アルゴリズム}
Here we will explain the computational notions used in our numerical methods.
For simplicity, we will first explain in setting of the surface heat equation \eqref{spde} on a smooth closed surface $S$ without boundary in three-dimensional Euclidean space. For the sake of clarity, we will also explain the detail in the setting of the surface wave equation \eqref{spde2}.

Let $\alpha>0$ denote the diffusion coefficient and $\Delta_S$ denote the Laplace-Beltrami operator on $S$. Given a time step of $h>0$, the algorithm for the surface-type minimizing movement that we developed is as follows.\\\\

\underline{{\bf{Surface-type minimizing movements for the surface heat equation \eqref{spde}}\label{spde1_al}}}
\begin{enumerate}
    \setlength{\leftskip}{1.0cm}
    \item \label{get_dxdydz}Create the computational domain $\Omega^D$ by preparing a sufficiently large Cartesian grid covering the surface $S$.
          Let $x_\text{min}$, $x_\text{max}$, $y_\text{min}$, $y_\text{max}$, $z_\text{min}$, and $z_\text{max}$ be the coordinates of the grid boundaries, as shown in Figure \ref{al1}.
          Let the grid spacing in the three spatial directions be given by $\Delta x$, $\Delta y$, and $\Delta z$, respectively.
          Then $\Omega^D$ is defined as follows:
          \begin{align*}
              \Omega^D=\left.\left\{\boldsymbol{x}_{i,j,k} =\left(
              \begin{matrix}x_i\\y_j\\z_k\end{matrix}
              \right)\right|0\le i \le N_x, 0\le j \le N_y, 0\le k \le N_z\right\}
          \end{align*}
          where $i$, $j$, and $k$ are natural numbers, and $N_x, N_y,$ and $N_z$ denote the number of grid points along the axes of the coordinate system.
          The grid points in the computational domain are expressed as follows:
          $$
              x_i=x_\text{min}+i\Delta x, \quad y_j=y_\text{min}+j\Delta y, \quad
              z_k=z_\text{min}+k\Delta z,
          $$
          $$
              N_x=\frac{x_\text{max}-x_\text{min}}{\Delta x}, \quad N_y=\frac{y_\text{max}-y_\text{min}}{\Delta y},\quad N_z=\frac{z_\text{max}-z_\text{min}}{\Delta z}
          $$
          For simplicity, we assume a uniform grid $\Delta x=\Delta y=\Delta z$.
          \begin{figure}[H]
              \centering
              \fbox{\includegraphics[bb=0cm 0cm 10cm 8cm,scale=0.7]{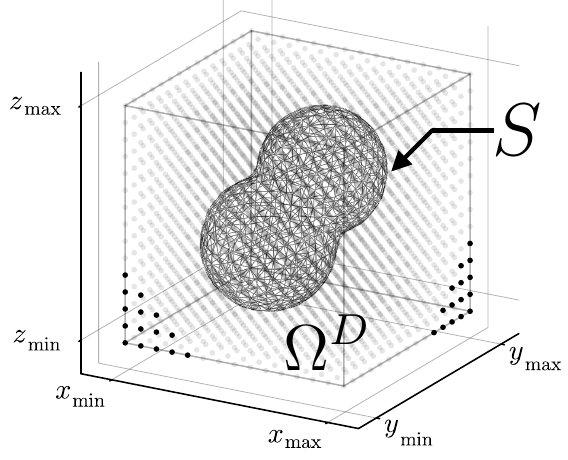}}
              \caption{Surface $S$ and computational domain $\Omega^D$}
              \label{al1}
          \end{figure}

    \item Using the closest point function $C_S$, compute and record the closest point on the surface $S$ for each point in $\Omega^D$.
    \item To reduce computational cost, calculations are performed only in a vicinity near the surface $S$.
          This process is called {\emph{banding}}.
          In particular, we extract a set of points from $\Omega^D$ whose Euclidean distance to the surface $S$ is less than or equal to a constant value $\lambda >0$ and denote the region by $\Omega^D_\lambda$.
          This is expressed as follows, where $||\cdot||$ represents the Euclidean norm.
          \begin{align}\label{getband}
              \Omega_\lambda^D=\left\{\boldsymbol{x}\in \Omega^D\hspace{5pt} |\hspace{5pt}||\boldsymbol{x}-C_S(\boldsymbol{x})||\le\lambda\right\}
          \end{align}
          \\
          {\bf{Remark:}}
          $\Omega_\lambda^D$ is a {\emph{point cloud}}; it consists of discrete points. In the continuous case, a sufficiently large region $\Omega\subset\mathbb{R}^3$ covering the surface $S$ is taken, and the region $\Omega_\lambda$ around the surface is defined as follows:
          \begin{align}\label{getband_cont}
              \Omega_\lambda=\left\{\boldsymbol{x} \in \Omega\hspace{5pt} |\hspace{5pt}||\boldsymbol{x}-C_S(\boldsymbol{x})||\le\lambda\right\}
          \end{align}
          {\bf{Remark:}}
          The value of $\lambda$ needs to be chosen appropriately, depending on the interpolation method used in Step \ref{heat_interp} below.
          If a polynomial interpolation is used, $\lambda$ depends on the degree of the interpolation.
          Here, we explain a method for determining $\lambda$ when performing a linear interpolation in a two-dimensional space (higher dimensions can be treated analogously).
          We assume that the grid points in the computational domain $\Omega^D$ have equal spacing in both the horizontal and vertical directions (Figure \ref{bandings}(a)).
          To obtain the interpolated value at the point denoted by ``$\star$" in Figure \ref{bandings}(b), four points denoted by ``$\bullet$" are required.
          In this case, the maximum distance between the interpolation point and the grid points is $\sqrt{2(\Delta x/2)^2}$.
          The maximum distance occurs in Figure \ref{bandings}(c), and its value is $\sqrt{2(\Delta x)^2}$.
          Therefore, $\lambda$ must be larger than $\sqrt{2(\Delta x)^2}$.
          Thus, one choice is to set $\lambda=\sqrt{(\Delta x)^2+2(\Delta x)^2}$ when a linear interpolation is used in a two-dimensional space.
          This discussion can be generalized to the case of a $d$-dimensional $p$th degree polynomial interpolation, then we obtain \cite{cpm}: $$\lambda=\sqrt{(d-1)\left(\frac{p+1}{2}\right)^2+\left(1+\frac{p+1}{2}\right)^2}\Delta x$$
          Since we are considering surfaces in three-dimensional space, we select $\lambda$ using the interpolation degree $p$ as follows:
          \begin{align}\label{lambda}
              \lambda=\sqrt{2\left(\frac{p+1}{2}\right)^2+\left(1+\frac{p+1}{2}\right)^2}\Delta x
          \end{align}
          \begin{figure}[H]
              \fbox{
                  \begin{tabular}{ccc}
                      \begin{minipage}[t]{0.3\hsize}
                          \centering
                          \includegraphics[bb=0cm 0cm 10cm 4cm]{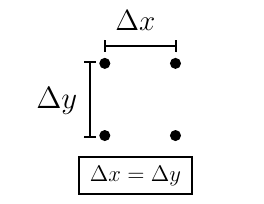}
                          \subcaption{}\label{banding1}
                      \end{minipage} &
                      \begin{minipage}[t]{0.3\hsize}
                          \centering
                          \includegraphics[bb=0cm 0cm 10cm 4cm]{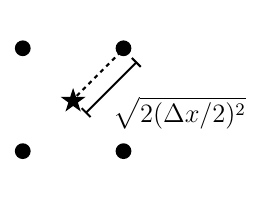}
                          \subcaption{}\label{banding2}
                      \end{minipage} &
                      \begin{minipage}[t]{0.3\hsize}
                          \centering
                          \includegraphics[bb=0cm 0cm 10cm 4cm]{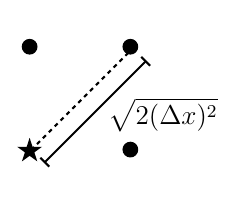}
                          \subcaption{}\label{banding3}
                      \end{minipage}
                  \end{tabular}
              }
              \caption{How to determine $\lambda$ in $\mathbb{R}^2$}
              \label{bandings}
          \end{figure}
    \item When approximating the gradient of a function in $\Omega_\lambda^D$, information about the boundary points is necessary. 
          We define the characteristic function $\phi_{i,j,k}$ as follows:
          \begin{align*}
              \phi_{i,j,k}=\begin{cases}
                  0,\quad \boldsymbol{x}_{i,j,k}\in{\Omega_\lambda^D} \\
                  1, \quad\text{otherwise},
              \end{cases}
          \end{align*}
          from which we define the boundary points $\partial\Omega_\lambda^D$ of $\Omega_\lambda^D$ as follows:
          \begin{align*}
              \partial\Omega_\lambda^D=\left\{\boldsymbol{x}_{i,j,k}\in\Omega^D \hspace{5pt}|\hspace{10pt}\phi_{i,j,k}|\nabla_{D}\phi_{i,j,k}|\neq0\right\}
          \end{align*}
          where $\nabla_{D}\phi_{i,j,k}=(\phi_{i+1,j,k}-\phi_{i-1,j,k},
              \phi_{i,j+1,k}-\phi_{i,j-1,k},
              \phi_{i,j,k+1}-\phi_{i,j,k-1})/(2\Delta x)$.
          We then join $\partial\Omega_\lambda^D$ with $\Omega_\lambda^D$ and define it as $\hat\Omega_\lambda^D$, that is,
          $$\hat\Omega_\lambda^D=\Omega_\lambda^D\cup\partial\Omega_\lambda^D.$$
          Figure \ref{al2} shows the relationship between $S$, $\Omega_\lambda^D$, and $\partial\Omega_\lambda^D$.
          Figure \ref{al2} is a schematic diagram of the section of Figure \ref{al1}.
    \item \label{initialize_u_heat}
          Using the closest point function $C_S$, extend the initial condition given on the surface $S$ to $\Omega^D$ as follows, where the initial condition at point $\boldsymbol{x}_{i,j,k}$ is denoted by $u^0_{i,j,k}$.
          \begin{align*}
              u^0_{i,j,k}=
              \begin{cases}
                  u^S_0(C_S(\boldsymbol{x}_{i,j,k})),\quad \boldsymbol{x}_{i,j,k} \in \hat \Omega_\lambda^D \\
                  0, \quad \boldsymbol{x}_{i,j,k}\in \Omega^D \setminus \hat \Omega_\lambda^D
              \end{cases}
          \end{align*}

          {\bf{Remark:}}
          In order to simplify the calculations, the initial values of the grid points outside of $\hat{\Omega}_\lambda^D$ are set to 0.
          \begin{figure}[H]
              \begin{center}
                  \fbox{\includegraphics[bb=0cm 0cm 11.5cm 9cm,scale=0.7]{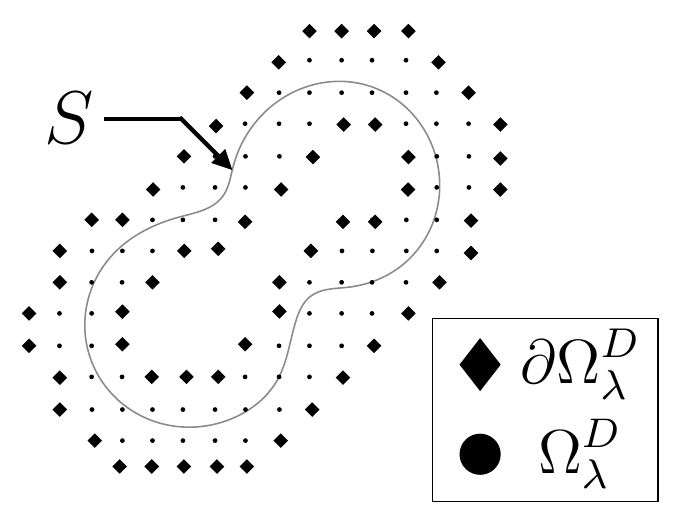}}
                  \caption{The relationship between $S$, $\Omega_\lambda^D$, and $\partial \Omega_\lambda^D$}
                  \label{al2}
              \end{center}
          \end{figure}
    \item \label{get_min_heat}Obtain an approximate solution of the heat equation on $\Omega_\lambda^D$ using MM.
          To approximate $u$ on $\Omega_\lambda^D$ in equation \eqref{mmheat_flat}, let $u_{i,j,k}=u(\boldsymbol{x}_{i,j,k})$.
          Approximate the functional values of \eqref{mmheat_flat} by means of an expression such as:
          \footnotesize
          \begin{align}\label{functional}
              {\mathcal{F}}_n(\boldsymbol{u})\approx\Delta x^3\sum_{\boldsymbol{x}_{i,j,k}\in \Omega_\lambda^D}\left\{\frac{|u_{i,j,k}-u_{i,j,k}^{n-1}|^2}{2h}+\alpha \frac{(\nabla_{D,x}u_{i,j,k})^2+(\nabla_{D,y}u_{i,j,k})^2+(\nabla_{D,z}u_{i,j,k})^2}{2}\right\}
          \end{align}
          \normalsize
          Denote the minimizer of this functional by $\boldsymbol{u}_{n}$, where $\boldsymbol{u}=(u_{i,j,k})$.
          Here, $\Delta x^3$ is the volume of the element, and $\nabla_{D,x}u_{i,j,k}, \nabla_{D,y}u_{i,j,k}, \nabla_{D,z}u_{i,j,k}$ are calculated by difference approximations as follows:
          \begin{align*}
              \nabla_{D,x}u_{i,j,k}=\frac{u_{i+1,j,k}-u_{i-1,j,k}}{2\Delta x}
          \end{align*}
          \begin{align*}
              \nabla_{D,y}u_{i,j,k}=\frac{u_{i,j+1,k}-u_{i,j-1,k}}{2\Delta x}
          \end{align*}
          \begin{align*}
              \nabla_{D,z}u_{i,j,k}=\frac{u_{i,j,k+1}-u_{i,j,k-1}}{2\Delta x}
          \end{align*}
          Note that, as mentioned earlier, we assume $\Delta x=\Delta y=\Delta z$.\\
          {\bf{Remark:}}
          Various methods can be used to obtain the minimizer of \eqref{functional}. Among them, from the viewpoint of computational cost, the L-BFGS method is often used \cite{Optim}.

    \item \label{heat_interp}Create an interpolating function  $I_n(\boldsymbol{x})$ defined on $\Omega_\lambda$ for the minimizer obtained in step \ref{get_min_heat}. 
          Using $I_n(\boldsymbol{x})$, define $u_n^S(\boldsymbol{x})=I_n(\boldsymbol{x})$ for $\boldsymbol{x}\in S$.
          It should again be noted that $I_n(\boldsymbol{x})$ is defined on $\Omega_\lambda$. There are multiple methodologies for its construction.
          One example is to use trilinear interpolation, which is a linear interpolation in 3D \cite{gomes2009implicit}.
          The computations in this study have used polynomial interpolations.

    \item \label{extrapolate}Using the closest point function $C_S$, extend $u_n^S$ onto $\hat\Omega_\lambda^D$ as follows:
          $$u^n_{i,j,k}=u_{n}^S(C_S(\boldsymbol{x}_{i,j,k})), \quad \boldsymbol{x}_{i,j,k} \in \hat \Omega_\lambda^D$$
    \item Repeat steps $\ref{get_min_heat}$ to $\ref{extrapolate}$ for $n=1,2,\cdots$ until the desired final time is reached.
\end{enumerate}

Next, we will explain the computational algorithm for the surface wave equation \eqref{spde2}.
\\\\
\underline{\bf{Surface-type minimizing movements for the surface wave  equation \eqref{spde2}}}
\begin{enumerate}
    \setlength{\leftskip}{1.0cm}
    \item Perform the computations in Steps \ref{get_dxdydz} to \ref{initialize_u_heat} of the previous algorithm.
    \item Assign $u_{i,j,k}^{-1}$ using the initial velocity $V_0$ of equation \eqref{spde2}.
        This can be done, for example, by means of the backward difference approximation $u_{i,j,k}^{-1}=u_{i,j,k}^{0}-hV_0(\boldsymbol{x}_{i,j,k})$.
          Note that $u_{i,j,k}^{-1}$ represents the value at the grid point $\boldsymbol{x}_{i,j,k}$ at time $-h$.
    \item \label{get_min_wave} Compute an approximate solution to the wave equation on $\Omega_\lambda^D$ using MM.
          Similar to the case of the heat equation, the functional values in \eqref{mmwave_flat} can be approximated as follows, for $n=1,2,\cdots$:
          \footnotesize
          \begin{align}\label{functional2}
              {\mathcal{F}}_n(\boldsymbol{u})\approx\Delta x^3\sum_{\boldsymbol{x}_{i,j,k}\in \Omega_\lambda^D}\left\{\frac{|u_{i,j,k}-2u_{i,j,k}^{n-1}+u_{i,j,k}^{n-2}|^2}{2h^2}+\alpha \frac{(\nabla_{D,x}u_{i,j,k})^2+(\nabla_{D,y}u_{i,j,k})^2+(\nabla_{D,z}u_{i,j,k})^2}{2}\right\}
          \end{align}
          \normalsize
          The minimizer of this functional is denoted by $\boldsymbol{u}_{n}$, where $\boldsymbol{u}=(u_{i,j,k})$.
          Here, $\Delta x^3$ is the volume of the element, and $\nabla_{D,x}u_{i,j,k}, \nabla_{D,y}u_{i,j,k}, \nabla_{D,z}u_{i,j,k}$ are calculated in the same way as in Step \ref{get_min_heat} of the previous algorithm.
    \item Define $u_n^S$ using the minimizer obtained in Step \ref{get_min_wave} by employing the same procedure as in Step \ref{heat_interp} of the previous algorithm.
    \item \label{extrapolate2}Using the closest point function $C_S$, extend $u_n^S$ onto $\hat\Omega_\lambda^D$.
    \item Repeat steps \ref{get_min_wave} to \ref{extrapolate2} for $n=1,2,\cdots$ until the desired final time is reached.
\end{enumerate}

In the next section, we perform a numerical error analysis using the above algorithms for the heat and wave equations on the surface of the unit sphere.
We will begin by treating the case of the surface heat equation.

\subsection{Numerical error analysis of MM for the heat equation on a surface}
Here, we will perform an numerical error analysis of the algorithm for solving the surface heat equation, described in the previous section.
Unsing MM, we numerically solve the surface heat equation \eqref{spde} on the unit sphere $S$, and examine the error between the numerical solution and the exact solution.
We define the unit sphere $S$ in the 3D space as follows:
\begin{align}\label{sphere_S}
    S=\left\{
    \left.\left(
    \begin{array}{c}
            \sin\theta\cos\phi \\[3pt]
            \sin\theta\sin\phi \\[3pt]
            \cos\theta
        \end{array}
    \right)\right|
    0\le \theta \le \pi,0\le \phi \le 2\pi
    \right\}
\end{align}

We perform two numerical experiments by changing the initial condition $f$ in equation \eqref{spde} on the unit sphere $S$.
First, we explain the initial conditions used and their corresponding exact solutions.
The results of the numerical error analysis are presented in Section \ref{数値誤差解析の結果(単位球面上熱方程式)}.
\subsection{MM and and surface heat equation: initial condition 1}
Setting the diffusion coefficient to $\alpha=1$, we take the initial condition $f$ as
$$
    f(\theta)=\cos\theta
$$
The exact solution of equation \eqref{spde} is then given by
$$
    u(\theta ,\phi,t)=e^{-2t}\cos\theta,\quad t\geq  0
$$
\subsection{MM and and surface heat equation: initial condition 2}
Setting the diffusion coefficient to $\alpha=1/42$, we take the initial condition $f$ as
$$
    f(\theta,\phi)=Y^{0}_6(\theta,\phi)+\sqrt{\frac{14}{11}}Y^{5}_6(\theta,\phi)
$$
where $Y^m_{l}(\theta,\phi)$ are the eigenfunctions of the Laplacian on the unit sphere, known as spherical harmonics.
The exact solution of equation \eqref{spde} is then given by
$$
    u(\theta ,\phi,t)=e^{-t}\left\{Y^{0}_6(\theta,\phi)+\sqrt{\frac{14}{11}}Y^{5}_6(\theta,\phi)\right\},\quad t\geq 0
$$
as shown in \cite{heatex1,heatex2}.

The results of the numerical error analysis using initial conditions 1 and 2 (described above) for the surface heat equation are described in the next section.
\subsection{MM numerical error analysis results (heat equation on the unit sphere)}\label{数値誤差解析の結果(単位球面上熱方程式)}
We investigate the relationship between $\Delta x$ and the numerical error of the MM approximation to the solution of the surface heat equation. Computations follow the computational algorithm for the heat equation on surfaces \eqref{spde}, presented in Section \ref{アルゴリズム}, where the spatial discretization $\Delta x$ is varied.
The L-BFGS method is used to minimize the discretized functional \eqref{functional}. We implement the method using Optim.jl \cite{Optim}, and calculate the functional gradient using automatic differentiation (ReverseDiff.jl \cite{ReverseDiff} is used for this purpose).
The time step $h$ is set to $h = \Delta x^2/6$, and polynomial interpolation of order $p=2$ is used (see Section \ref{補間}).
For both initial conditions 1 and 2, we calculate the maximum absolute error $L_{\infty}$ on $S$ at the closest point to each point in $\Omega_\lambda^D$ at time $t_e=0.25$. We note that, since the exact solution of the surface heat equation converges to 0, it becomes difficult to evaluate the error between the numerical solution and the exact solution. For this reason, we have chosen such a value of $t_e$ (i.e., so that the $L_{\infty}$-error of the absolute value of the exact solution at time $t_e$ is sufficiently large for both initial conditions 1 and 2).
The $L_{\infty}$-error is defined as
\begin{gather*}
    L_{\infty}\mathchar`-\text{error}=\sup_{\boldsymbol{x}\in \Omega_\lambda^D}|u(C_S(\boldsymbol{x}),t_e)-\hat u(C_S(\boldsymbol{x}),t_e)|
\end{gather*}
where $\hat u$ is the numerical solution and $u$ denotes the exact solution.

\begin{figure}[H]
    \begin{minipage}{0.5\hsize}
        \begin{center}
            \fbox{\includegraphics[bb=0 0 354 248,scale=0.5]{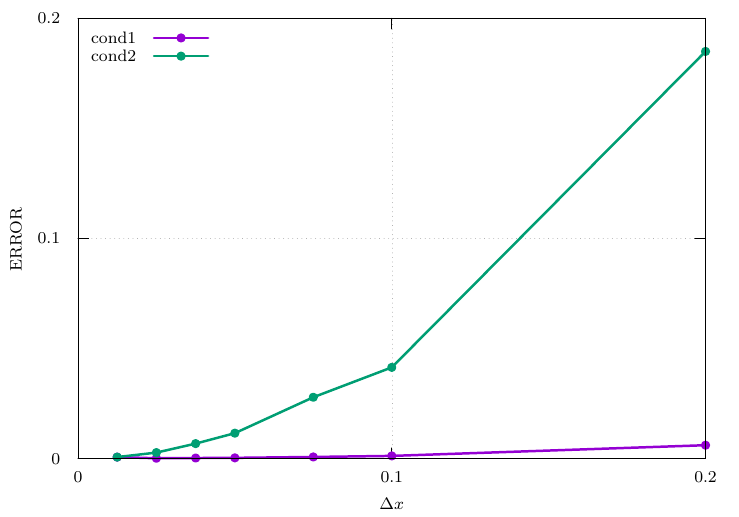}}
            \subcaption{}
            \label{heaterr}
        \end{center}
    \end{minipage}
    \begin{minipage}{0.5\hsize}
        \begin{center}
            \fbox{\includegraphics[bb=0 0 354 248,scale=0.5]{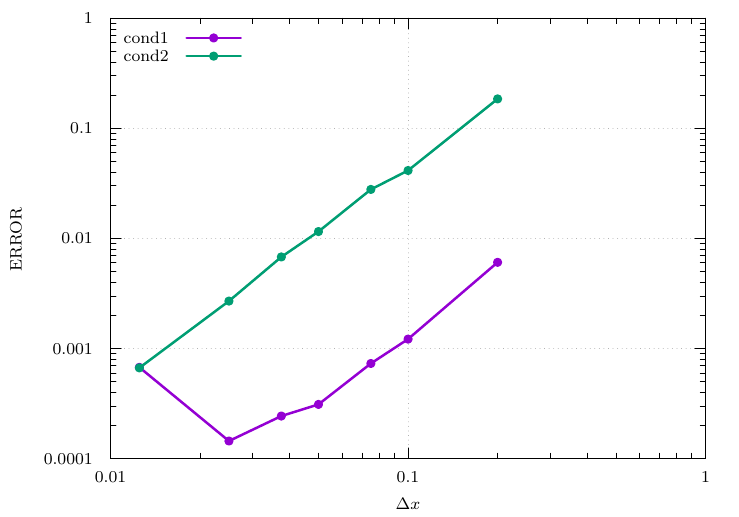}}
            \subcaption{
                }
                \label{heaterr_log}
        \end{center}
    \end{minipage}

    \caption{ (\subref{heaterr}) Numerical error for the surface heat equation, (\subref{heaterr_log}) Numerical error for the surface heat equation (log-log plot). "cond1" corresponds to the initial condition 1 and "cond2" corresponds to the initial condition 2. Except for $\Delta x = 0.0125$, it is observed that the error decreases as $\Delta x$ decreases. The numerical error is approximately proportional to $\Delta x^2$ squared.
        }
        \label{heatERR}
\end{figure}

\begin{figure}[H]
    \fbox{
        \begin{tabular}{ccc}
            \begin{minipage}[t]{0.3\hsize}
                \centering
                \includegraphics[width=\linewidth]{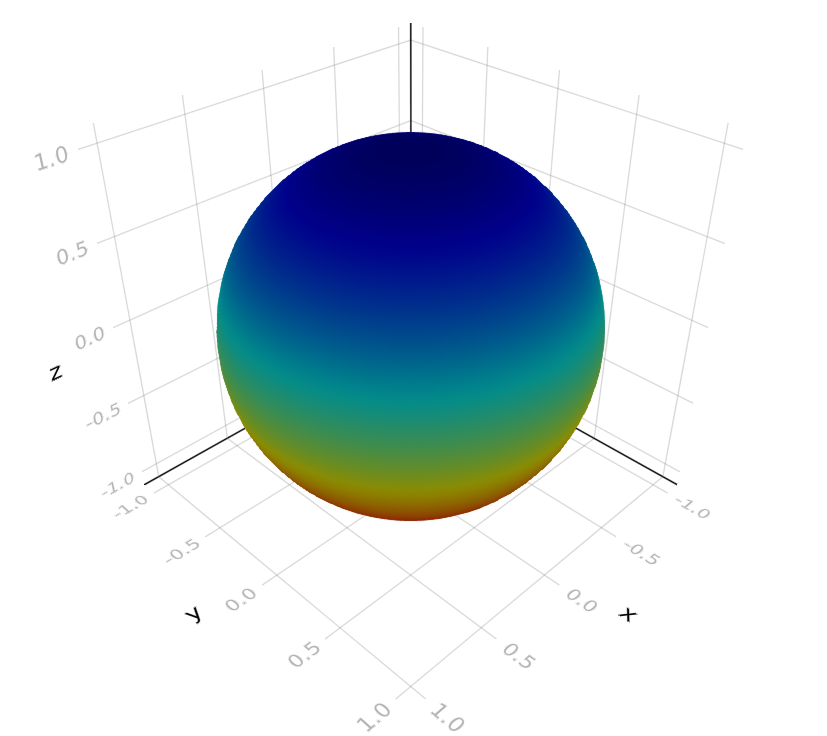}
                \subcaption{}
                \label{heatcond1_ini}
            \end{minipage} &
            \begin{minipage}[t]{0.3\hsize}
                \centering
                \includegraphics[width=\linewidth]{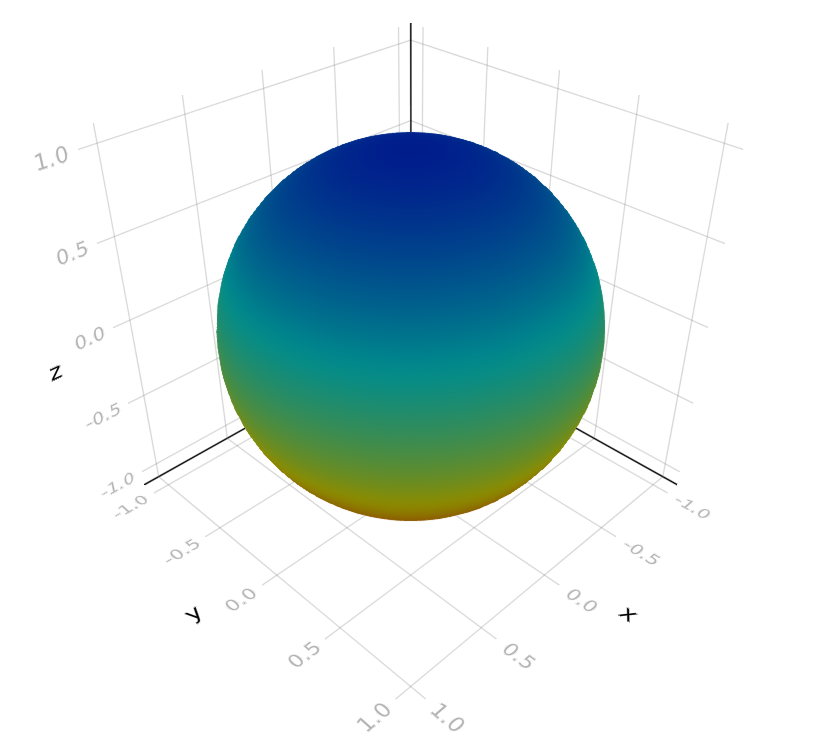}
                \subcaption{}
                \label{gomi1}
            \end{minipage} &
            \begin{minipage}[t]{0.3\hsize}
                \centering
                \includegraphics[width=\linewidth]{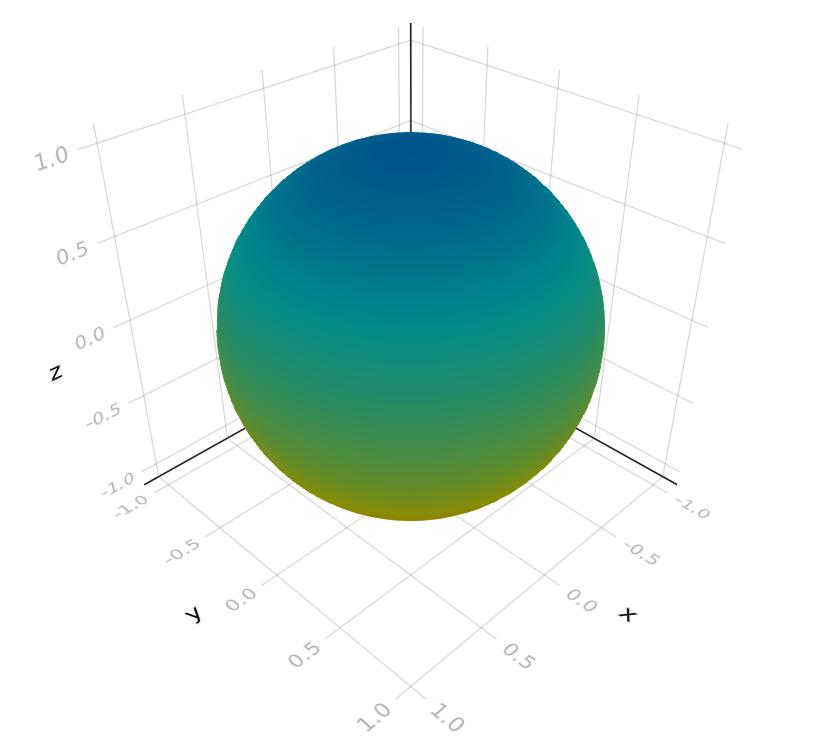}
                \subcaption{}
                \label{gomi2}
            \end{minipage}
            \\
            \begin{minipage}[t]{0.3\hsize}
                \centering
                \includegraphics[width=\linewidth]{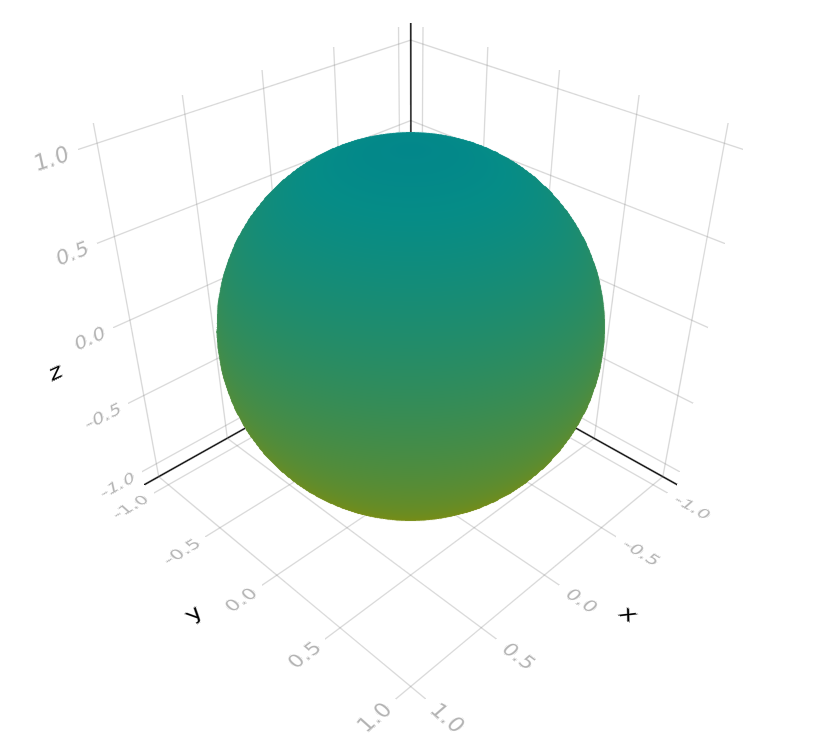}
                \subcaption{}
                \label{gomi3}
            \end{minipage} &
            \begin{minipage}[t]{0.3\hsize}
                \centering
                \includegraphics[width=\linewidth]{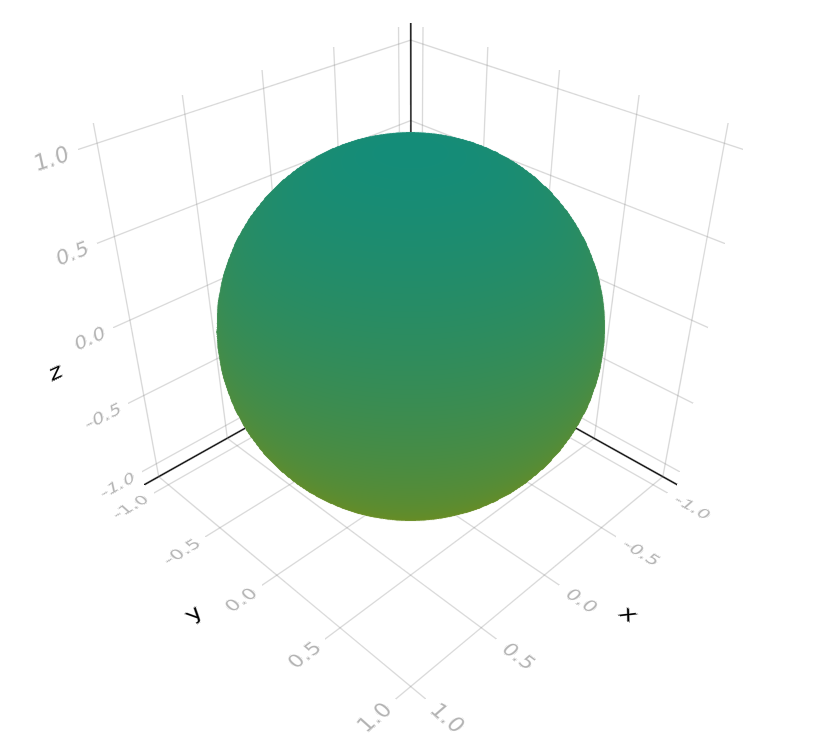}
                \subcaption{}
                \label{gomi4}
            \end{minipage} &
            \begin{minipage}[t]{0.3\hsize}
                \centering
                \includegraphics[width=\linewidth]{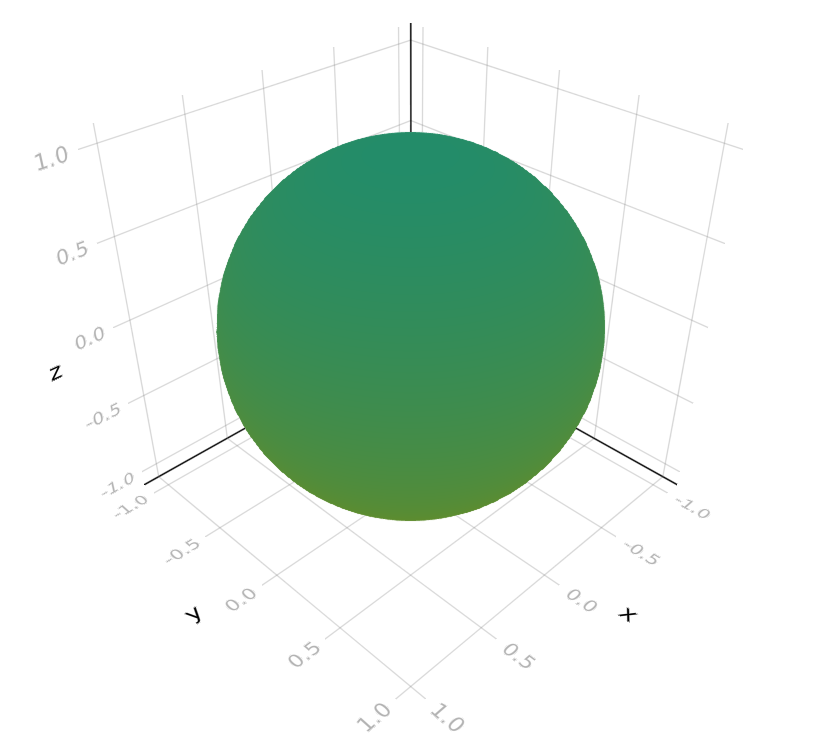}
                \subcaption{}
                \label{gomi5}
            \end{minipage}
        \end{tabular}
        }
    \caption{Initial condition and computation result (Initial condition 1): (\subref{heatcond1_ini}) shows the initial condition, and the subsequent subfigures show the time evolution, in alphabetical order.}
    \label{heatcond1}
\end{figure}

\begin{figure}[H]
    \fbox{
        \begin{tabular}{ccc}
            \begin{minipage}[t]{0.3\hsize}
                \centering
                \includegraphics[width=\linewidth]{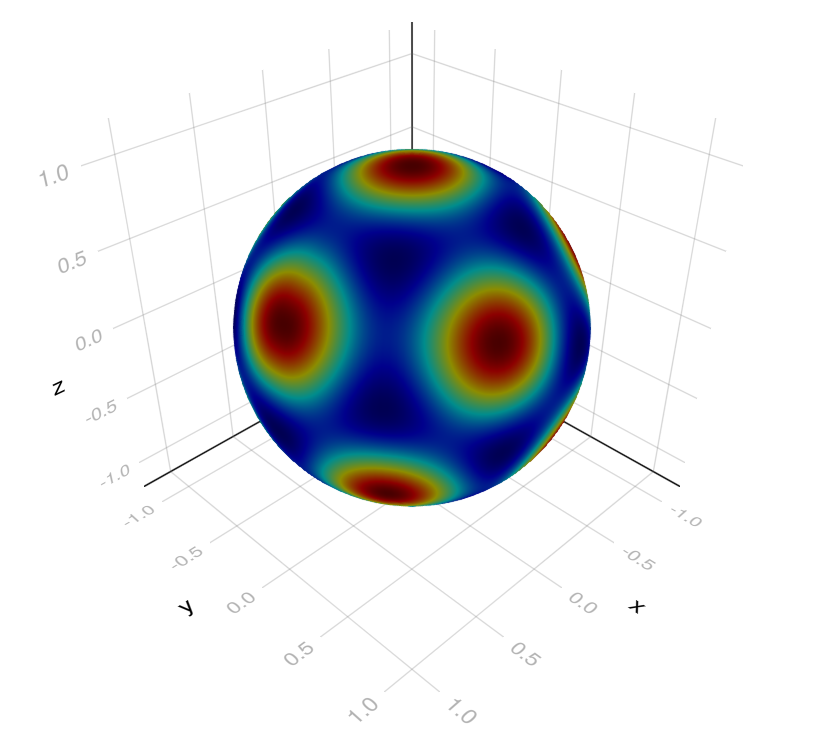}
                \subcaption{}
                \label{heatcond2_ini}
            \end{minipage} &
            \begin{minipage}[t]{0.3\hsize}
                \centering
                \includegraphics[width=\linewidth]{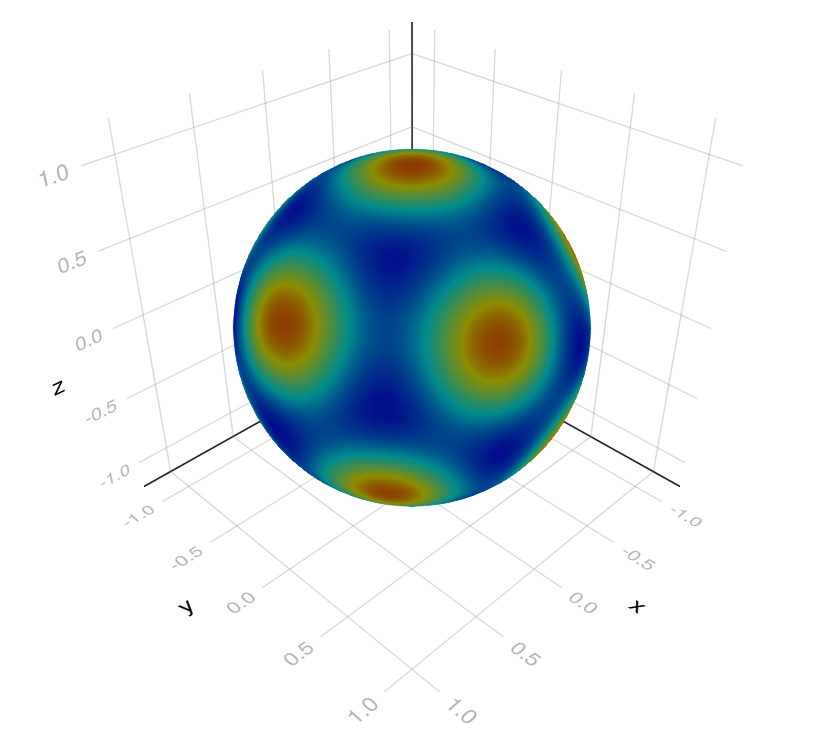}
                \subcaption{}
                \label{}
            \end{minipage} &
            \begin{minipage}[t]{0.3\hsize}
                \centering
                \includegraphics[width=\linewidth]{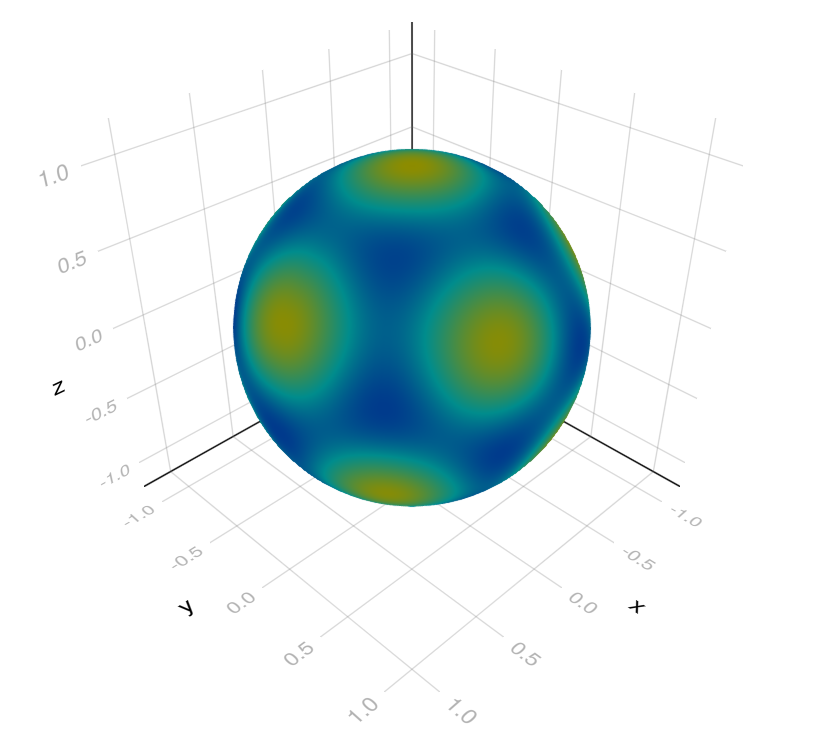}
                \subcaption{}
            \end{minipage}   \\
            \begin{minipage}[t]{0.3\hsize}
                \centering
                \includegraphics[width=\linewidth]{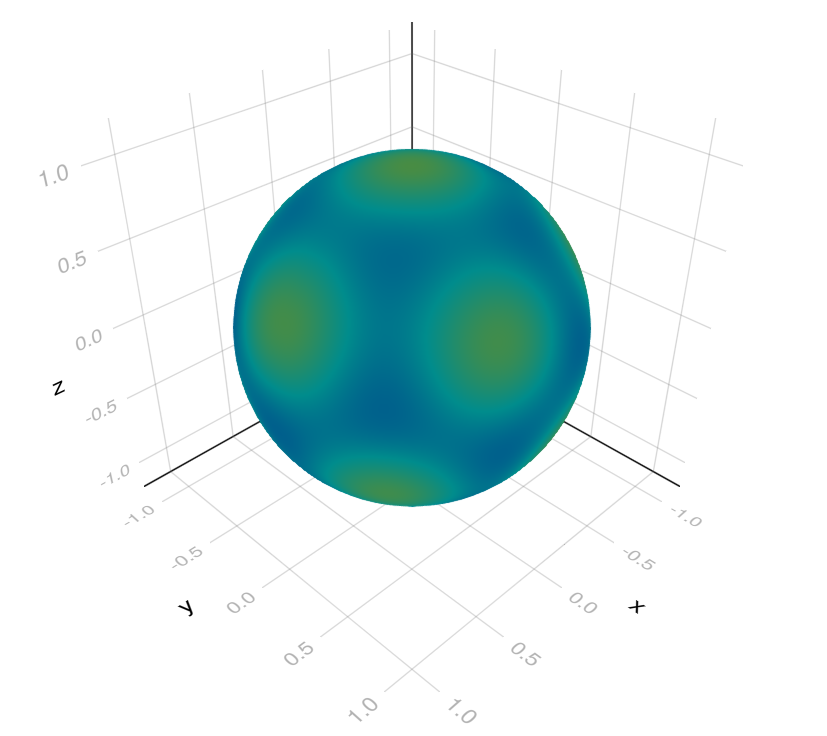}
                \subcaption{}
            \end{minipage} &
            \begin{minipage}[t]{0.3\hsize}
                \centering
                \includegraphics[width=\linewidth]{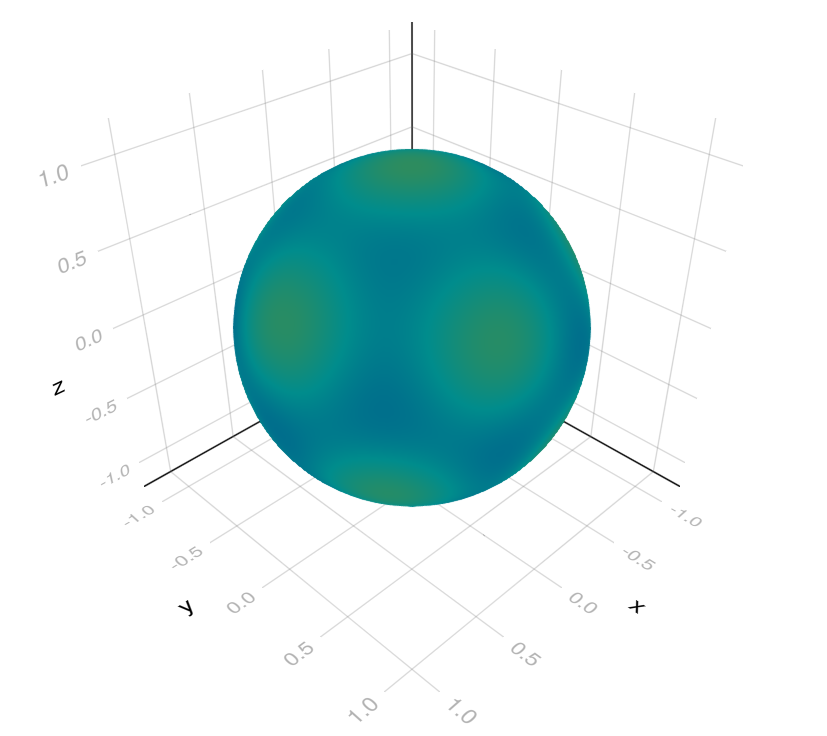}
                \subcaption{}
                \label{}
            \end{minipage} &
            \begin{minipage}[t]{0.3\hsize}
                \centering
                \includegraphics[width=\linewidth]{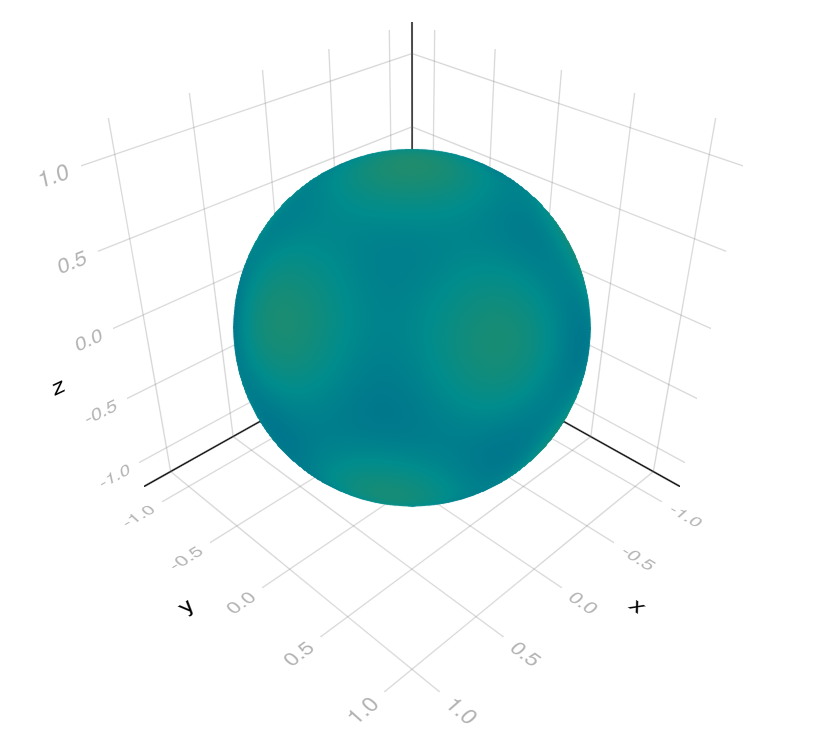}
                \subcaption{}

            \end{minipage}
        \end{tabular}
    }
    \caption{Initial condition and computation result (Initial condition 2): (\subref{heatcond2_ini}) shows the initial condition, and the subsequent subfigures show the time evolution, in alphabetical order.}
    \label{heatcond2}
\end{figure}

The results obtained for each $\Delta x$ are shown in Table \ref{tab:res1} and Table \ref{tab:res2}.
The results are plotted in Figure \ref{heatERR}(\subref{heaterr}) and Figure \ref{heatERR}(\subref{heaterr_log}) using both regular and log-log scales, respectively.
The legend in the figure denotes initial condition 1 by {\emph{cond1}}, and initial condition 2 by {\emph{cond2}}.
The time evolution is shown in Figure \ref{heatcond1} and Figure \ref{heatcond2}.
The results confirm that the $L_{\infty}$-error decreases as $\Delta x$ decreases, except for the case where $\Delta x=0.0125$.
Except for this case, the numerical error is roughly proportional to the square of $\Delta x$ when $\Delta x$ is reduced by a factor of 2.
The reason for the larger numerical error in initial condition 2 may be due to insufficient resolution relative to the initial condition.

These results confirm that the numerical solution obtained by MM converges to the exact solution of the surface heat equation as the spatial discretization converges to zero.
We also note that we observe the numerical error increases when $\Delta x$ becomes numerically too small (Figure \ref{heatERR}(\subref{heaterr_log})).
Next, we will perform a numerical error analysis for the wave equation on a curved surface.
\begin{center}
\begin{table}[h]
    \begin{minipage}{0.5\hsize}
        \caption{Results for initial condition 1}
        \centering
        \begin{tabular}{c|c}
            $\Delta x$ & $L_{\infty}$-error \\ \hline \hline 
            0.2        & 6.061e-03          \\ \hline
            0.1        & 1.218e-03          \\ \hline 
            0.075      & 7.310e-04          \\ \hline 
            0.05       & 3.103e-04          \\ \hline 
            0.0375     & 2.434e-04          \\ \hline 
            0.025      & 1.443e-04          \\ \hline  
            0.0125     & 6.747e-04          \\ \hline
        \end{tabular}
        \label{tab:res1}
    \end{minipage}
    \begin{minipage}{0.5\hsize}
        \caption{Results for initial condition 2}
        \centering
        \begin{tabular}{c|c}
            $\Delta x$ & $L_{\infty}$-error \\ \hline \hline 
            0.2        & 1.849e-01          \\ \hline
            0.1        & 4.142e-02          \\ \hline 
            0.075      & 2.784e-02          \\ \hline
            0.05       & 1.154e-02          \\ \hline 
            0.0375     & 6.783e-03          \\ \hline
            0.025      & 2.696e-03          \\ \hline  
            0.0125     & 6.666e-04          \\ \hline
        \end{tabular}
        \label{tab:res2}
    \end{minipage}
\end{table}
\end{center}


\subsection{Numerical error Analysis for the surface wave equation}
Using the MM approximation scheme previously described, we numerically solve the wave equation on a curved surface \eqref{spde2} and examine the error between the numerical and the exact solution.
We perform two computational experiments by changing the initial condition $f$ of equation \eqref{spde2} on the unit sphere $S$ expressed by equation \eqref{sphere_S}.
First, we explain the initial conditions used and their exact solutions.
The results of the numerical error analysis are presented in Section \ref{数値誤差解析の結果(球面上波動方程式)}.
\subsection{Surface wave equation: initial condition 1}
With the constant $\alpha=1$, we set the initial condition $f$ as
$$
    f(\theta)=-\cos\theta
$$
and the initial velocity $V_0$ as
$$
    V_0=0
$$
In this case, the exact solution of equation \eqref{spde2} is given by
$$
    u(\theta,\phi,t)=-\cos(\sqrt{2}t)\cos\theta,\quad t\geq  0
$$
\subsection{Surface wave equation: initial condition 2}
Let $\alpha = 1$ and the initial condition $f$ be
$$
    f(\theta ,\phi)=Y^{0}_6(\theta,\phi)+\sqrt{\frac{14}{11}}Y^{5}_6(\theta,\phi)
$$
with initial velocity $V_0 = 0$.
Then, the exact solution of equation \eqref{spde2} is given by
$$
    u(\theta,\phi,t)=\cos(\sqrt{42}t)\left\{Y^{0}_6(\theta,\phi)+\sqrt{\frac{14}{11}}Y^{5}_6(\theta,\phi)\right\},\quad t\geq  0
$$
Next, we will explain the results of the numerical error analysis using initial conditions 1 and 2 for the surface wave equation.
\subsection{Numerical error analysis results (wave equation on the unit sphere)}\label{数値誤差解析の結果(球面上波動方程式)}
We performed several numerical simulations using the algorithm presented in Section \ref{アルゴリズム} for the wave equation \eqref{spde2} and investigated the relationship between the numerical error and the spatial discretization $\Delta x$.
We used the same optimization methods and interpolation degree $p$ as those used for the surface heat equation.
The time step $h$ was set to $\Delta x/10$.
For initial condition 1, the absolute error $L_{\infty}$ is calculated at the closest point on $S$ for each point in $\Omega_\lambda^D$ at the time $t=2\pi/\sqrt{2}$.
For initial condition 2, the absolute error $L_{\infty}$ is calculated at the closest point on $S$ for each point in $\Omega_\lambda^D$ at the time $t=2\pi/\sqrt{42}$.
The maximum value of the absolute error is then obtained for each case.
We evaluated the $L_\infty$ error at the time $t_e$ such that the exact solution using initial conditions 1 and 2 have both  oscillated once over $0\leq t\leq t_e$.
The $L_\infty$ error is defined as follows:
\begin{gather*}
    L_{\infty}\mathchar`-\text{error}=\sup_{\boldsymbol{x}\in \Omega_\lambda^D}|u(C_S(\boldsymbol{x}),t_e)-\hat u(C_S(\boldsymbol{x}),t_e)|
\end{gather*}
where $\hat u$ denotes the numerical solution, and $u$ denotes the exact solution.

\begin{figure}[H]
    \begin{minipage}{0.5\hsize}
        \begin{center}
            \fbox{\includegraphics[bb=0 0 354 248,scale=0.5]{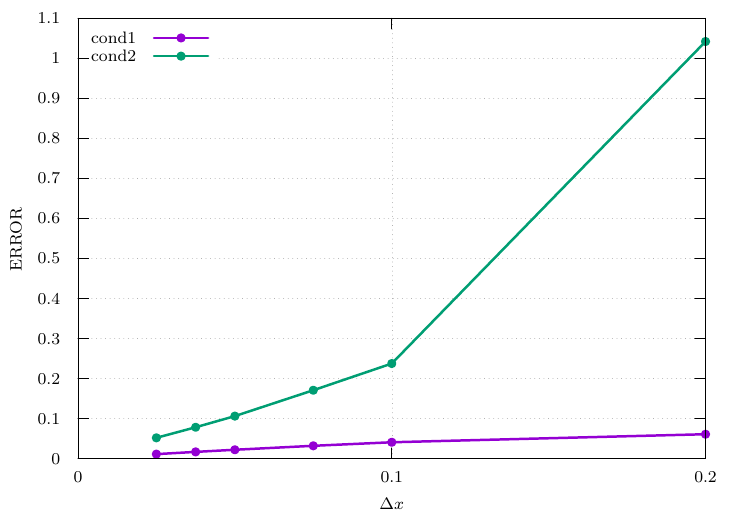}}
            \subcaption{}
            \label{waveerr}
        \end{center}
    \end{minipage}
    \begin{minipage}{0.5\hsize}
        \begin{center}
            \fbox{\includegraphics[bb=0 0 354 248,scale=0.5]{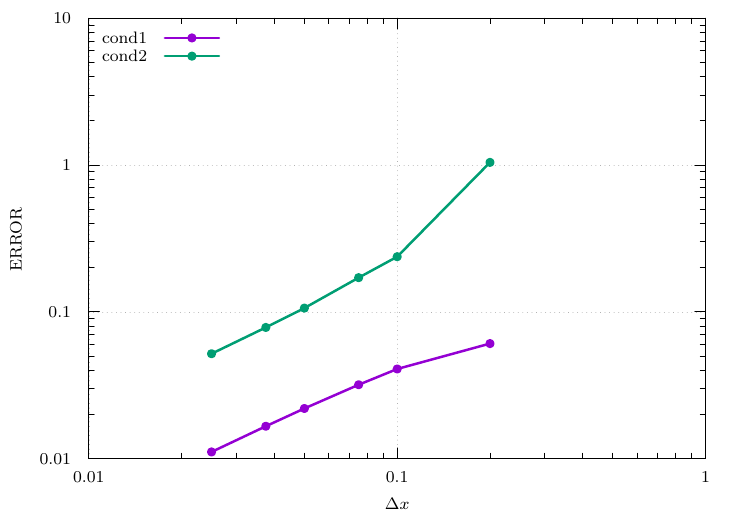}}
            \subcaption{}
            \label{waveerr_log}
        \end{center}
    \end{minipage}
    \caption{ (\subref{waveerr}) Numerical error for the wave equation on a curved surface, (\subref{waveerr_log}) Numerical error for the wave equation on a curved surface (log-log plot).
        "cond1" corresponds to initial condition 1 and "cond2" corresponds to initial condition 2. It can be observed that the error decreases as the value of $\Delta x$ becomes smaller. The numerical error is approximately proportional to $\Delta x$.}
        \label{waveERR}
\end{figure}

\begin{figure}[H]
    \fbox{
        \begin{tabular}{ccc}
            \begin{minipage}[t]{0.3\hsize}
                \centering
                \includegraphics[width=\linewidth]{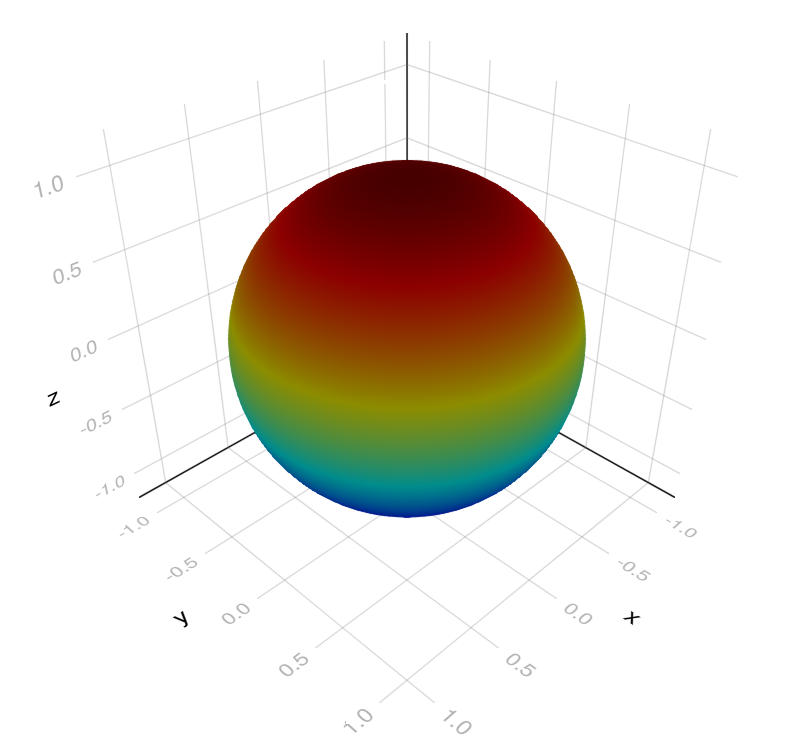}
                \subcaption{}
                \label{wavecond1_ini}
            \end{minipage} &
            \begin{minipage}[t]{0.3\hsize}
                \centering
                \includegraphics[width=\linewidth]{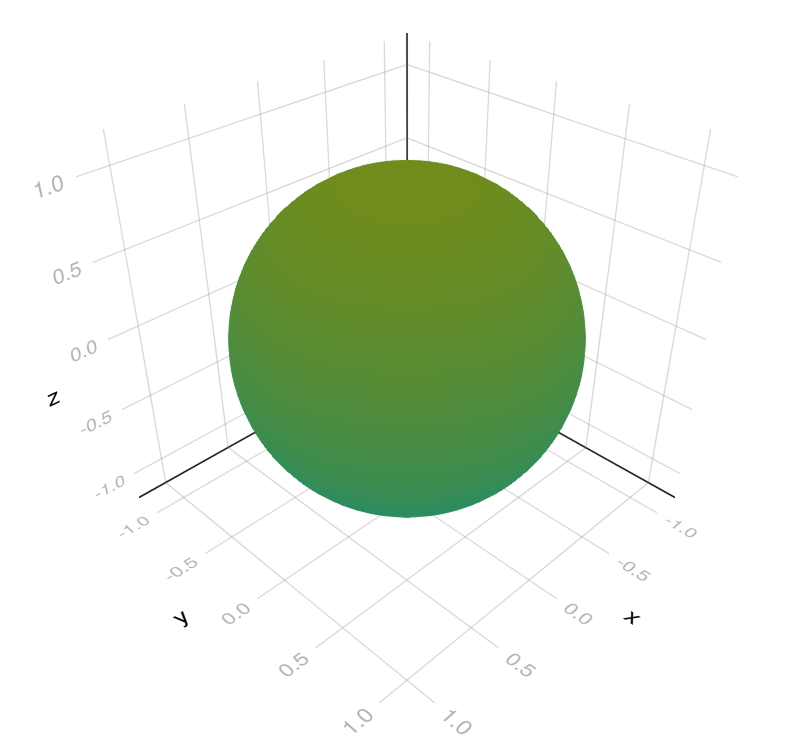}
                \subcaption{}
                \label{}
            \end{minipage} &
            \begin{minipage}[t]{0.3\hsize}
                \centering
                \includegraphics[width=\linewidth]{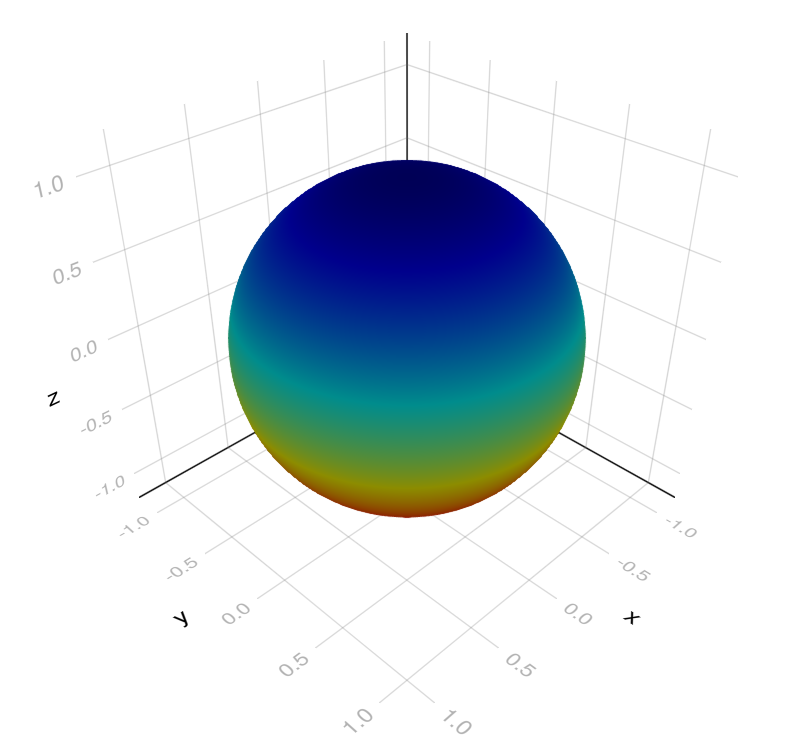}
                \subcaption{}
            \end{minipage}   \\
            \begin{minipage}[t]{0.3\hsize}
                \centering
                \includegraphics[width=\linewidth]{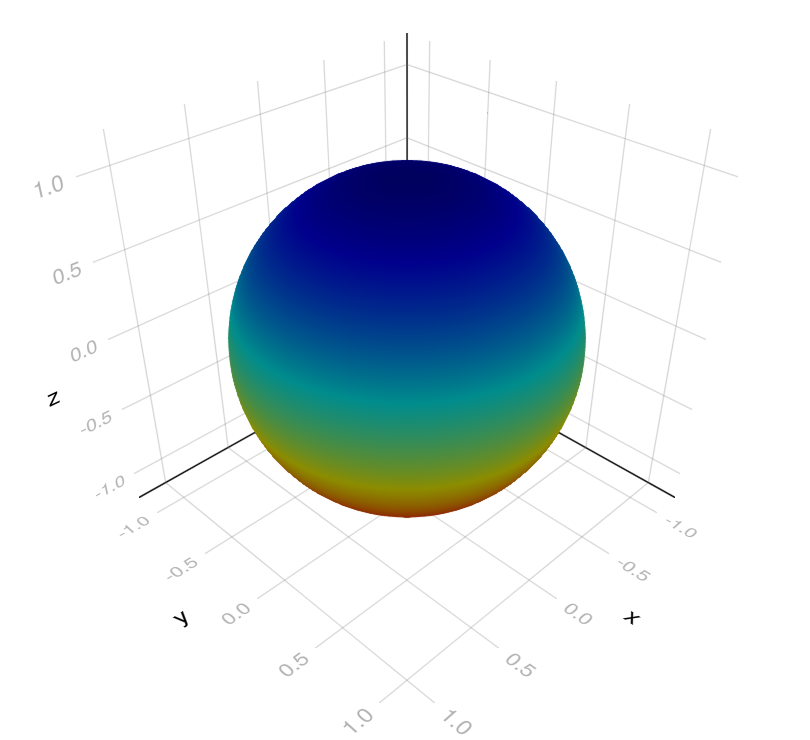}
                \subcaption{}
            \end{minipage} &
            \begin{minipage}[t]{0.3\hsize}
                \centering
                \includegraphics[width=\linewidth]{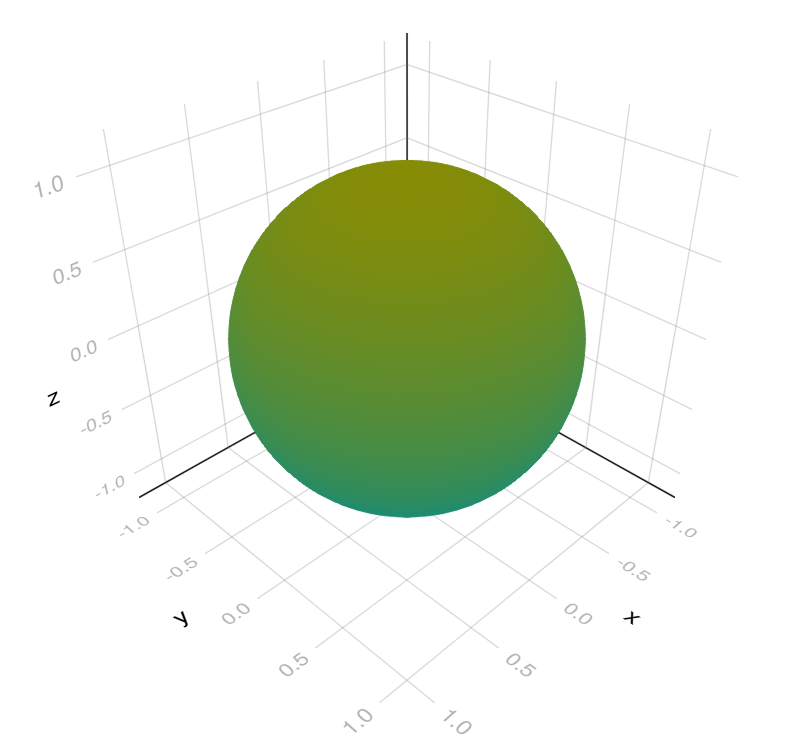}
                \subcaption{}
                \label{}
            \end{minipage} &
            \begin{minipage}[t]{0.3\hsize}
                \centering
                \includegraphics[width=\linewidth]{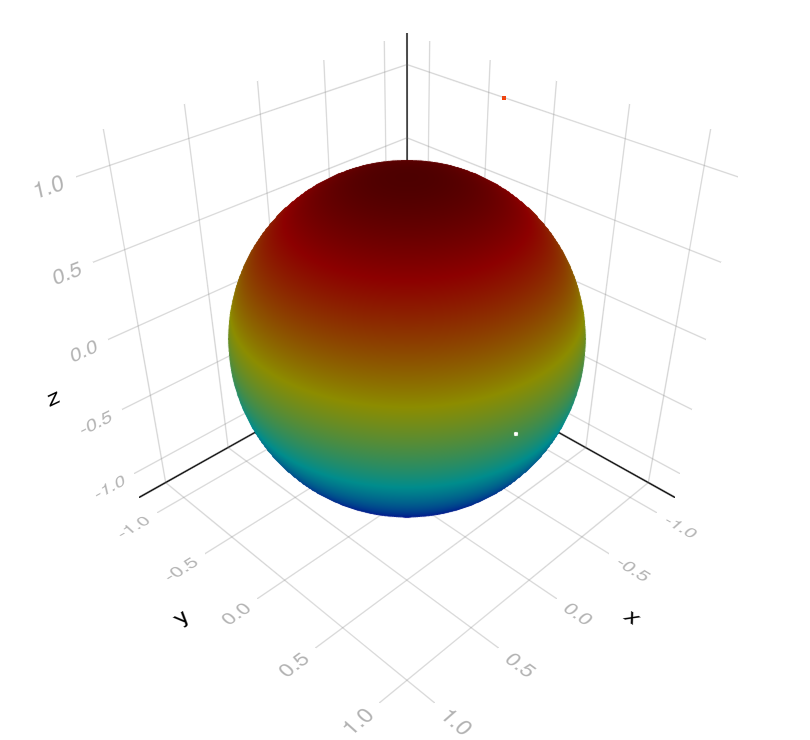}
                \subcaption{}
            \end{minipage}
        \end{tabular}
    }
    \caption{Initial condition and computation results (initial condition 1): (\subref{wavecond1_ini}) shows the initial condition, and the subsequent subfigures show the time evolution, in alphabetical order.}
    \label{wavecond1}
\end{figure}

\begin{figure}[H]
    \fbox{
        \begin{tabular}{ccc}
            \begin{minipage}[t]{0.3\hsize}
                \centering
                \includegraphics[width=\linewidth]{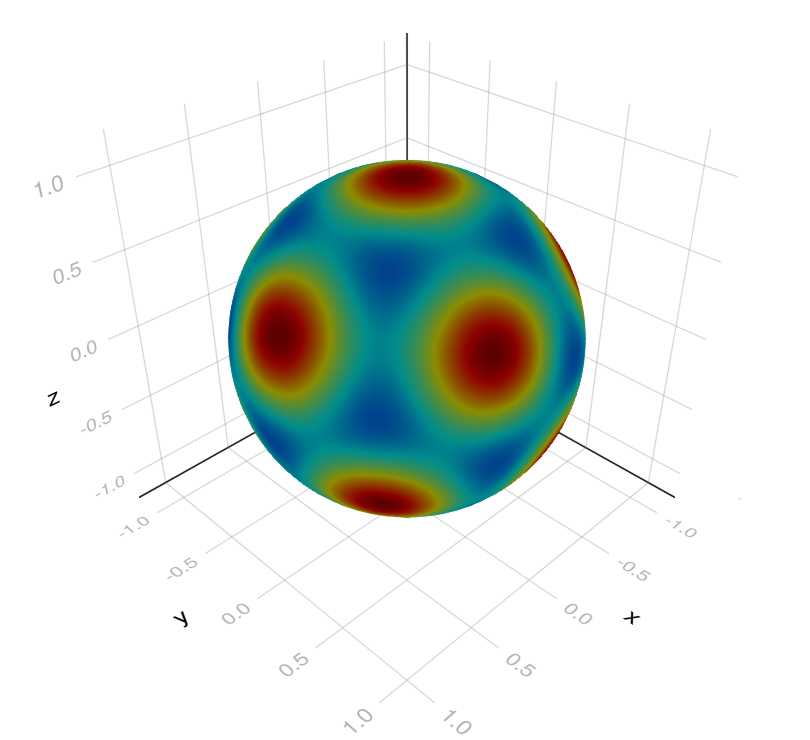}
                \subcaption{}
                \label{wavecond2_ini}
            \end{minipage} &
            \begin{minipage}[t]{0.3\hsize}
                \centering
                \includegraphics[width=\linewidth]{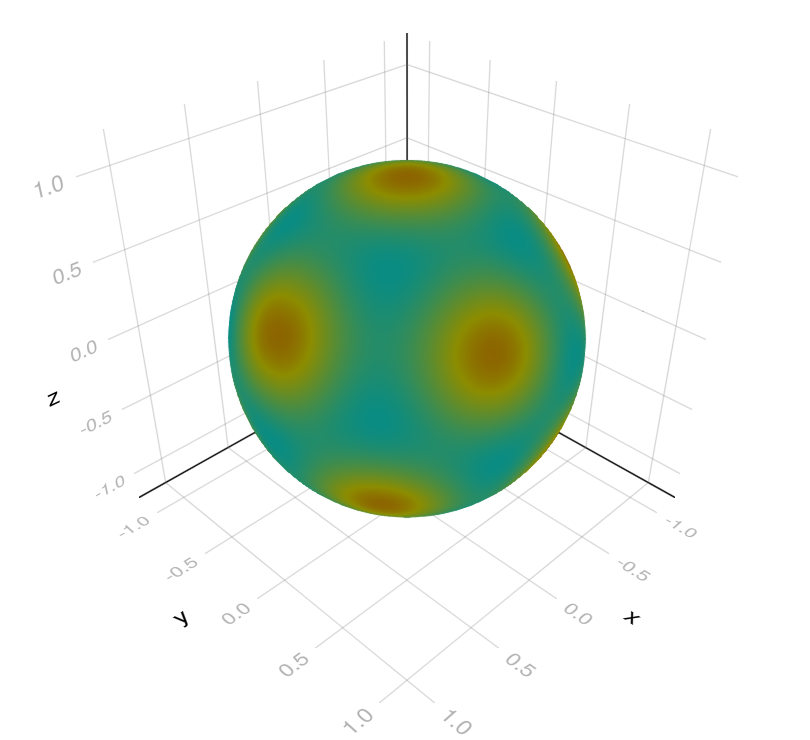}
                \subcaption{}
                \label{}
            \end{minipage} &
            \begin{minipage}[t]{0.3\hsize}
                \centering
                \includegraphics[width=\linewidth]{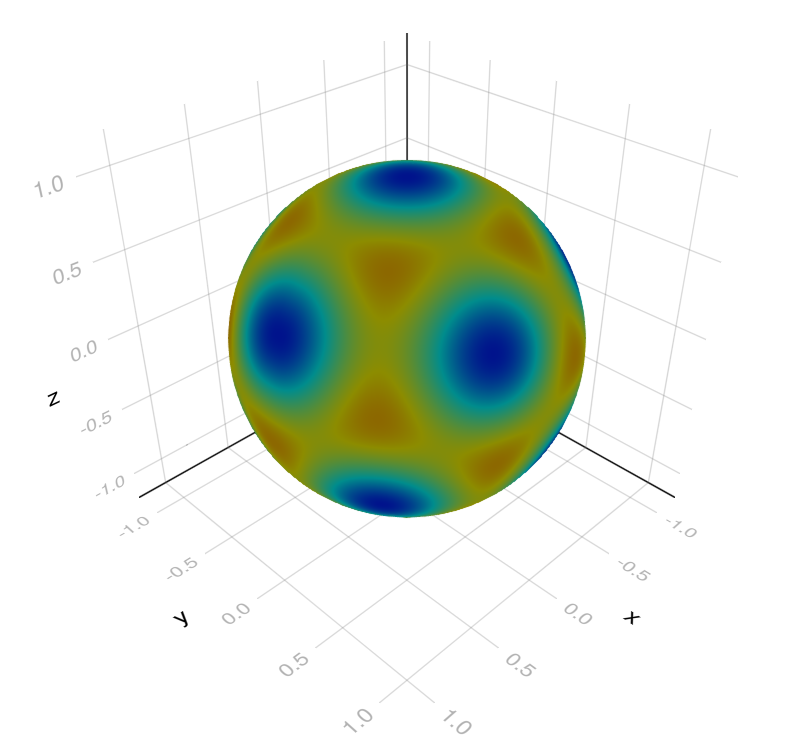}
                \subcaption{}
            \end{minipage}   \\
            \begin{minipage}[t]{0.3\hsize}
                \centering
                \includegraphics[width=\linewidth]{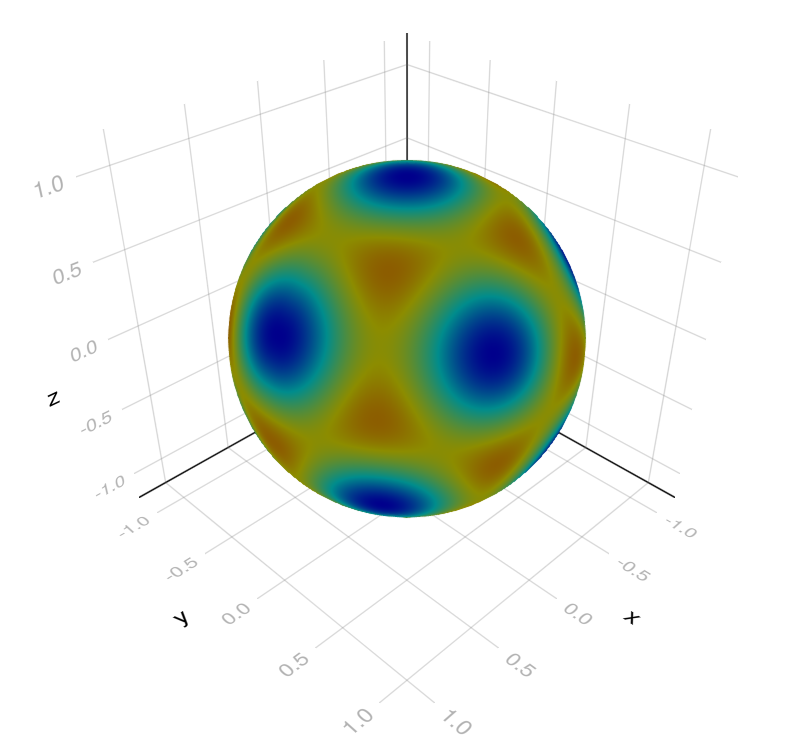}
                \subcaption{}
            \end{minipage} &
            \begin{minipage}[t]{0.3\hsize}
                \centering
                \includegraphics[width=\linewidth]{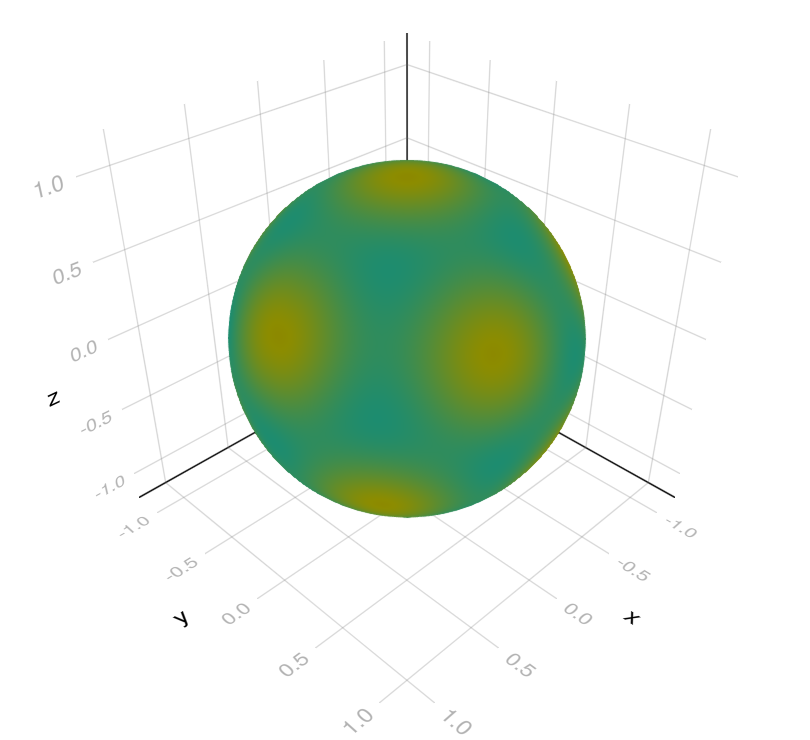}
                \subcaption{}
                \label{}
            \end{minipage} &
            \begin{minipage}[t]{0.3\hsize}
                \centering
                \includegraphics[width=\linewidth]{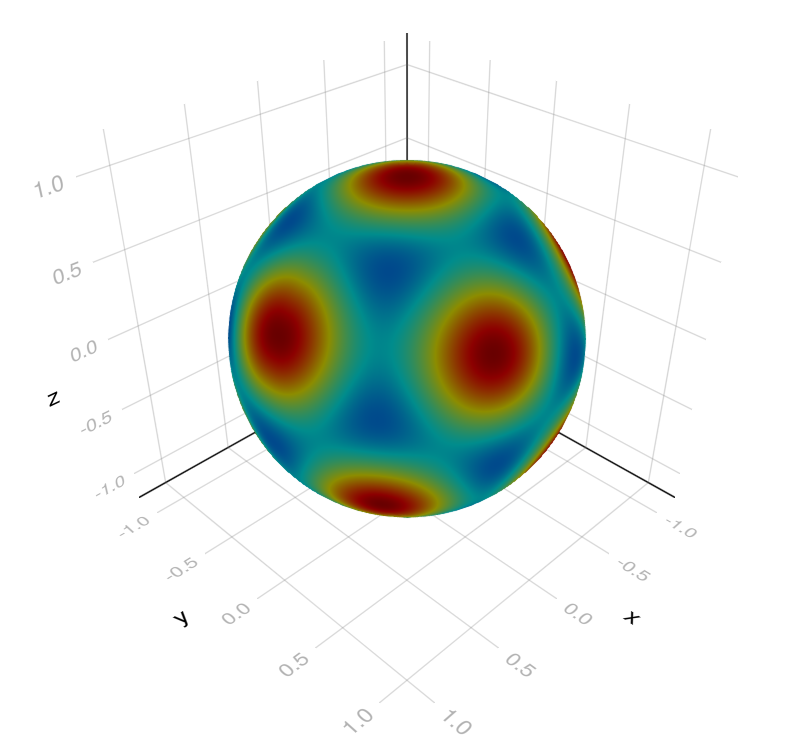}
                \subcaption{}
            \end{minipage}
        \end{tabular}
    }
    \caption{Initial condition and computation results (initial condition 1): (\subref{wavecond2_ini}) shows the initial condition, and the subsequent subfigures show the time evolution, in alphabetical order.}
    \label{wavecond2}
\end{figure}

The results obtained for $\Delta x$ and $L_{\infty}\mathchar`-\text{error}$ are presented in Table \ref{tab:wave_res1} and Table \ref{tab:wave_res2}.
Figure \ref{waveERR}(\subref{waveerr}) shows a graphical representation of the contents of Table \ref{tab:wave_res1} and Table \ref{tab:wave_res2}, while Figure \ref{waveERR}(\subref{waveerr_log}) shows a representation on log-log scale.
In the legend of the figures, {\emph{cond1}} corresponds to initial condition 1, and {\emph{cond2}} corresponds to initial condition 2.
The evolution of the solution over time is shown in Figure \ref{wavecond1} and Figure \ref{wavecond2}.
The results show that the $L_{\infty}\mathchar`-\text{error}$ decreases with $\Delta x$.
Except for when $\Delta x=0.2$, the numerical error is approximately halved when $\Delta x$ is halved.
That is, it is observed that the numerical error converges proportional to $\Delta x$.
Similar to the numerical error analysis for the heat equation, it was found that the numerical error for the initial condition 2 is greater than that for initial condition 1.
These results indicate that when the spatial grid spacing $\Delta x$ is sufficiently small, the numerical solution obtained by the developed numerical method for the wave equation converges to the exact solution.

In this section, we explained the approximation methods for partial differential equations on  curved surfaces using the CPM and MM methods, and presented the results of their numerical error analysis.
In the following sections, we will illustrate applications of the approximation methods to the simulation of interfacial motions on curved surfaces.
\begin{center}
\begin{table}[h]
    \begin{minipage}{0.5\hsize}
        \caption{Results for initial condition 1}
        \centering
        \begin{tabular}{c|c}
            $\Delta x$ & $L_{\infty}$-error \\ \hline \hline 
            0.2        & 6.073e-02          \\ \hline
            0.1        & 4.079e-02          \\ \hline 
            0.075      & 3.183e-02          \\ \hline 
            0.05       & 2.195e-02          \\ \hline 
            0.0375     & 1.660e-02          \\ \hline 
            0.025      & 1.109e-02          \\ \hline  

        \end{tabular}
        \label{tab:wave_res1}
    \end{minipage}
    \begin{minipage}{0.5\hsize}
        \caption{Results for initial condition 2}
        \centering
        \begin{tabular}{c|c}
            $\Delta x$ & $L_{\infty}$-error \\ \hline \hline 
            0.2        & 1.042e+00          \\ \hline
            0.1        & 2.373e-01          \\ \hline 
            0.075      & 1.707e-01          \\ \hline
            0.05       & 1.060e-01          \\ \hline 
            0.0375     & 7.819e-02          \\ \hline
            0.025      & 5.178e-02          \\ \hline  
        \end{tabular}
        \label{tab:wave_res2}
    \end{minipage}
\end{table}
\end{center}


\section{Numerical simulation of interfacial motions on surfaces}\label{曲面上界面運動の数値計算}
In this section, we discuss approximation methods for realizing mean curvature flow and hyperbolic mean curvature flow on surfaces.
After explaining the interfacial motions on curved surfaces, we present the numerical results of our schemes.
We also deal with the mean curvature flow and the hyperbolic mean curvature flow on curved surfaces in the multiphase setting and under area preservation conditions.
Multiphase regions are realized by means of a signed distance vector field, which allows us to incorporate area preservation constraints into the functional minimizations of the MM approach.
\subsection{Surface-constrained interfacial motions}\label{曲面上界面運動}
Here we will discuss the approximation methods for the mean curvature flow (MCF) \cite{mc_surf} and the hyperbolic mean curvature flow (HMCF) \cite{hmcf} on curved surfaces.
The surface MCF and surface HMCF are described by the following nonlinear partial differential equations, respectively.

\begin{align}\label{surface_MCF}
    \begin{cases}
        \gamma^S_t(t,s)=-\kappa^S(t,s)\nu^S(t,s), \\
        \gamma^S(0,s)=\gamma^S_0(s)\tag{Surface MCF}
    \end{cases}
\end{align}
\begin{align}\label{surface_HMCF}
    \begin{cases}
        \gamma^S_{tt}(t,s)=-\kappa^S(t,s)\nu^S(t,s), \\
        \gamma^S(0,s)=\gamma^S_0(s),                 \\
        \gamma_t^S(0,s)=v_0\nu^S(t,s)\tag{Surface HMCF}
    \end{cases}
\end{align}
Here, $S$ is a smooth surface, $\gamma^S:[0,T)\times [a,b]\rightarrow S$ is a smooth simple curve on the surface $S$, satisfying $\gamma^S(t,a)=\gamma^S(t,b)$, $\kappa^S$ is the curvature of the curve, $\gamma_0$ is the initial shape of the curve, $v_0$ is the initial velocity of the curve, and $\nu^S$ represents the outward unit normal vector of the curve on the surface $S$.
Here, $\gamma^S_t=\partial \gamma^S/\partial t,\quad\gamma^S_{tt}=\partial \gamma^S/\partial tt$.
Surface MCF and surface HMCF are equations that generalize the motion of surfaces following mean curvature flow or hyperbolic mean curvature flow in the Euclidean space to the setting of interfaces moving on curved surfaces.
Interfaces moving by surface MCF tend to decrease their length and smoothing their shape over time, while interfaces moving by surface HMCF tend to oscillate.

In the case that surface MCF and surface HMCF are subject to the area-preservation conditions, the interfaces should move while preserving the areas of the regions enclosed by the interfaces.
We handle such motions by extending the MBO algorithm and the HMBO algorithm to the surface PDE setting using CPM, MM, and a surface version of the signed distance vector field, to develope approximate solutions for surface MCF (surface MBO) and surface HMCF (surface HMBO).
We remark that our methods can also handle interfacial motions with area-preservation conditions  in the multiphase setting.
This is enabled by means of a signed distance vector field that is used to encod the shape of interfaces atop the surface.
\subsection{The signed distance vector field on surfaces}\label{曲面版符号付き距離ベクトル場}
In this section, we discuss the signed distance vector field \cite{hmcf} and its extension to the surface setting.
The signed distance vector field (SDVF) is used to encode the shape of multiphase regions by means of vector directions.
When performing numerical calculations under area preserving conditions, the signed distance can be used in the two-phase setting.
However, for interface motions involving three or more phases and area preservation conditions, it is not possible to distinguish each phase using a single signed distance function.
On the other hand, the SDVF can be used to distinguish phase locations and shapes even in the case of three or more phases.
The SDVF is constructed by assigning a special vector to each phase of the multiphase region. Each vector is weighted by its signed distance from each interface \cite{hmcf,multi_volume}. The SDVF is described below.

\begin{figure}[H]
    \begin{center}
        \fbox{\includegraphics[bb=0cm 0cm 5cm 4.5cm,scale=1.0]{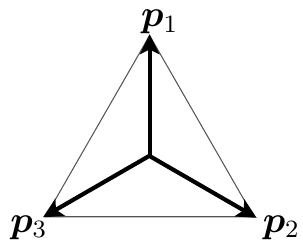}}
        \caption{Vectors pointing from the barycenter of a 2D simplex to its vertices $\boldsymbol{p}_1, \boldsymbol{p}_2,$ and $\boldsymbol{p}_3$}\label{seitantai}
    \end{center}
\end{figure}

Let $K$ be the number of phases, $\epsilon > 0$ be an interpolation parameter, $P^i$ be the region of phase $i$, $\boldsymbol{p}_i$ be the vector from the barycenter of a $(K-1)$-dimensional simplex to each vertex (refer to Figure \ref{seitantai} for $K=3$, \cite{hmcf}), $d_S^i(\boldsymbol{x})$ be the signed distance function to phase $i$, and $\chi_E$ be the characteristic function of the set $E$.
The signed distance vector field on the surface $S$ (Surface SDVF) $z_S^\epsilon$ is given by the following equation:
\begin{align}\label{sdvf}
    z_S^\epsilon (\boldsymbol{x})=\sum_{i=1}^K\left(\boldsymbol{p}_i\chi_{\{d_S^i\geq \epsilon/2\}}+\frac{1}{\epsilon}\left( \frac{\epsilon }{2}+d_S^i\right)\boldsymbol{p}_i\chi_{\{-\epsilon/2<d_S^i <\epsilon/2\}}\right), \quad \boldsymbol{x} \in S
\end{align}
where,
\begin{align}\label{chids}
    \chi_E(\boldsymbol{x})=\begin{cases}
        1\quad \boldsymbol{x}\in E, \\
        0\quad \text{otherwise},
    \end{cases}
    d_S^i(\boldsymbol{x})=\begin{cases}
        \displaystyle{\inf_{\boldsymbol{y}\in\partial P^i}}||\boldsymbol{x}-\boldsymbol{y}||_S\quad \boldsymbol{x}\in P^i, \\
        -\displaystyle{\inf_{\boldsymbol{y}\in\partial P^i}}||\boldsymbol{x}-\boldsymbol{y}||_S\quad \text{otherwise}.
    \end{cases}
\end{align}
Here, $||\boldsymbol{x}-\boldsymbol{y}||_S$ represents the geodesic distance on the surface between the two points $\boldsymbol{x}$ and $\boldsymbol{y}$, and can be expressed as the value that minimizes the following length functional:
$$
    ||\boldsymbol{x}-\boldsymbol{y}||_S=\min_{\Gamma\subset S}\text{Length($\Gamma$)}
$$
where $\Gamma$ is a curve on the surface connecting $\boldsymbol{x}$ and $\boldsymbol{y}$, and $\text{Length($\Gamma$)}$ is the length of $\Gamma$ along the surface.
\\{\bf{Remark:}}
From here on we will omit the $S$ in $z_S^\epsilon$ and $d_S^i$, and denote them simply by $z^\epsilon$ and $d^i$, respectively.

In the next Section \ref{algo-MBO}, we introduce an approximation method for Surface MCF, based on the MBO algorithm using the surface SDVF. We refer to appoximation method as the surface MBO algorithm.
In Section \ref{algo-HMBO}, we introduce an approximation method for the surface HMCF, which is based on the HMBO algorithm and the surface SDVF. Similarly, we refer to our approximation method as the surface HMBO.

\subsection{Surface MBO}\label{algo-MBO}
In this section, we discuss an approximation method for surface MCF.
We begin by describing the surface MBO in the two-phase setting.
Then we explain the surface MBO with and without the area preservation condition in the multiphase setting.
The computational results of our surface MBO in the multiphase setting are presented in Section \ref{results_MBO}.

We use the CPM as an approximation for the surface partial differential equation on the surface.
To reduce the computational cost of the CPM, we limit the computations to a small tubular region around $S$, as described in Section \ref{アルゴリズム}.
This region is defined by accumulating points from the computation grid that are located within a constant distance $\lambda$ from $S$ as $\Omega_\lambda$ (see equation \eqref{getband_cont}).
We determine the value of $\lambda$ in the same way as in equation (\ref{getband}).
We use a natural number $N$ and the final time $T>0$ to define the time step $\tau=T/N$.
\subsubsection{Surface MBO for two-phase regions}\label{algo-2MBO}
Below, we explain a method for approximating the MCF motion of an interface in a 2-phase region. Here, 
$\gamma^S_n$ represents the shape of the curve at time $n\tau$, where $\tau=T/N$ is the time step size, $T>0$ is the final time, and $\gamma^S_0$ is the initial curve.
Let $d_n(\boldsymbol{x})$ be the signed distance function from the point $\boldsymbol{x}$ on the surface to the interface $\gamma^S_n$.
That is, $d_n(\boldsymbol{x})$ is defined as follows:
\begin{align}\label{sd}
    \displaystyle
    d_n(\boldsymbol{x})=\begin{cases}
        \displaystyle{\inf_{\boldsymbol{y}\in\partial P^n}}||\boldsymbol{x}-\boldsymbol{y}||_S,\quad  & \boldsymbol{x}\in P^n, \\
        \displaystyle{-\inf_{\boldsymbol{y}\in\partial P^n}}||\boldsymbol{x}-\boldsymbol{y}||_S,\quad & \text{otherwise}.
    \end{cases}
\end{align}
where $P^n$ is the region occupied by phase $n$. 

The surface MBO for a 2-phase domain is as follows:
\begin{description}
    \setlength{\leftskip}{1.0cm}

    \item[1.] Create $d_0$ using Equation \eqref{sd} from the initial curve $\gamma^S_0$.
    \item[2.] Extend $d_0$ to $\Omega_\lambda$:
        $$
            d_0^\lambda(\boldsymbol{x})=d_0(C_S(\boldsymbol{x})), \quad\boldsymbol{x}\in\Omega_\lambda
        $$
    \item[3.] Repeat the following for $n=0,1,\cdots,N-1$:
        \begin{enumerate}[a.]
            \setlength{\leftskip}{1.0cm}
            \item For $t\in[0,\tau)$, solve the heat equation in $\Omega_\lambda$:
                  \begin{align}\label{algo-2MBO-heat}
                      \begin{cases}
                          u_{t}=\alpha\Delta u(\boldsymbol{x},t), \\
                          u(\boldsymbol{x},0)=d^\lambda_{n}(\boldsymbol{x}),
                      \end{cases}
                      \quad \boldsymbol{x}\in\Omega_\lambda,\quad \tau>t>0
                  \end{align}
                  where $\alpha>0$ is a constant representing the diffusion coefficient.
            \item Define the new curve $\gamma^S_{n+1}$ as the zero level set of the solution $u(\boldsymbol{x},\tau)$ of Equation \eqref{algo-2MBO-heat} on the surface $S$:
                  $$
                      \gamma^S_{n+1}=\{\boldsymbol{x}\in\Omega_\lambda|S\cap\{u(\boldsymbol{x},\tau)=0\}\}
                  $$
            \item Create $d_{n+1}$ from $\gamma^S_{n+1}$ using Equation \eqref{sd}.
            \item Extend $d_{n+1}$ to $\Omega_\lambda$ and define $d^\lambda_{n+1}$ as follows:
                  $$
                      d_{n+1}^\lambda(\boldsymbol{x})=d_{n+1}(C_S(\boldsymbol{x})),\quad\boldsymbol{x}\in\Omega_\lambda
                  $$
        \end{enumerate}
\end{description}
The results of numerical error analysis using the proposed algorithm are shown in Section \ref{2相領域に対するSurface MBOアルゴリズムの数値誤差解析}.
In the above algorithm, the sign of the signed distance function is used to extract the interface.
In the next section, we will implement the surface MBO for multiphase regions using the signed distance vector field instead of the signed distance function.
\subsubsection{Surface MBO for multiphase regions}\label{algo-mMBO}
Here, we explain a method for approximating surface MCF on interfaces consisting of $K\ge 2$ multiphase regions on the surface $S$.
Let $P^i_n$ represent the region of phase $i$ at time $n\tau$, and let $\boldsymbol{P}_n=\bigcup_{i=1}^K{P^i_n}$.
We provide initial regions $P_0^i$ for each phase $i$ as initial conditions.
Additionally, we denote the vector given to phase $i$ as $\boldsymbol{p}_i$, as explained in Section \ref{曲面版符号付き距離ベクトル場}.
The surface SDVF, written $z_n^\epsilon(\boldsymbol{x})$, is obtain from \eqref{sdvf} using $\boldsymbol{P}_n$.
The surface MBO for multiphase regions described as follows.
\begin{description}
    \setlength{\leftskip}{1.0cm}
    \item[1.] Using equation \eqref{sdvf}, create the surface SDVF $z_0^\epsilon$ from the initial domain $\boldsymbol{P}_0$.
    \item[2.] Extend $z_0^\epsilon$ to $\Omega_\lambda$ and denote it as $z_0^{\epsilon,\lambda}$ as follows:
        $$
            z_0^{\epsilon,\lambda}(\boldsymbol{x})=z_0^{\epsilon}(C_S(\boldsymbol{x})), \quad \boldsymbol{x}\in\Omega_\lambda
        $$
    \item[3.] Repeat the following for $n=0,1,\cdots,N-1$:
        \begin{enumerate}[a.]
            \setlength{\leftskip}{1.0cm}
            \item  Solve the following vector-valued heat equation for $t\in[0,\tau)$:
                  \begin{align}\label{heat-vec}
                      \begin{cases}
                          \boldsymbol{u}_{t}=\alpha\Delta \boldsymbol{u}(\boldsymbol{x},t),          \\
                          \boldsymbol{u}(\boldsymbol{x},0)=z^{\epsilon,\lambda}_{n}(\boldsymbol{x}), \\
                      \end{cases}
                      \quad \boldsymbol{x}\in\Omega_\lambda,\quad \tau>t>0
                  \end{align}
                  Here, $\alpha>0$ is a constant representing the diffusion coefficient.
            \item Extract the solution $\boldsymbol{u}(\boldsymbol{x},\tau)$ of equation \eqref{heat-vec} on the surface $S$ and denote it as $\boldsymbol{u}^S$ as follows:
                  $$
                      \boldsymbol{u}^S(\boldsymbol{x})=\boldsymbol{u}(\boldsymbol{x},\tau),\quad \boldsymbol{x}\in S\cap\Omega_\lambda
                  $$
            \item Obtain $\boldsymbol{P}_{n+1}$ on the surface $S$ using $\boldsymbol{u}^S$ as follows:
                  $$
                      \boldsymbol{P}_{n+1}=\bigcup_{i=1}^K\{P^i_{n+1}\} \quad P^i_{n+1} =\left\{ \boldsymbol{x}\in S;\boldsymbol{u}^S(\boldsymbol{x})\cdot \boldsymbol{p}_i\geq\boldsymbol{u}^S(\boldsymbol{x})\cdot \boldsymbol{p}_j,\text{ for all }j\in\{1,\cdots,K\}\right\}
                  $$
            \item Update the surface SDVF $z_{n+1}^\epsilon$ from $\boldsymbol{P}_{n+1}$ using Equation \eqref{sdvf}.
            \item Extend $z_{n+1}^\epsilon$ to $\Omega_\lambda$:
                  \begin{align*}
                      z_{n+1}^{\epsilon,\lambda}(\boldsymbol{x})=z_{n+1}^\epsilon(C_S(\boldsymbol{x})), \quad\boldsymbol{x}\in\Omega_\lambda
                  \end{align*}
        \end{enumerate}
\end{description}
Next, we will implement the surface MBO with a prescribed area-constraint for multiphase regions by combining MM with the above algorithm.

\subsubsection{Multiphase surface MBO with area constraints}\label{algo-smboc}
Here, we will explain our method for approximating multiphase surface MCF under $K$ area constraints on the surface $S$.
The initial vector field is constructed from the regions $P_0^i$ for each phase $i$, where the vector for each phase $i$, denoted by $\boldsymbol{p}_i$, is prescribed as in Section \ref{曲面版符号付き距離ベクトル場}.
The sign-distance vector field $z_n^\epsilon(\boldsymbol{x})$ is then obtained via equation \eqref{sdvf} from $\boldsymbol{P}_n$.
When approximating constrained interfacial motions, MM are often used when treating the case of area-preservation \cite{multi_volume}.

In particular, MM can be used in combination with penalty for each area constraint.

We take a sufficiently large positive constant $M$ and set $h=\tau/M$.
Given $\boldsymbol{w}_{m-1}$, we define the functional \eqref{mboF} used in the MM as follows:
\begin{align}\label{mboF}
    \mathcal{F}_m(\boldsymbol{w})=\int_{\Omega_\lambda}\left(\frac{|\boldsymbol{w}-\boldsymbol{w}_{m-1}|^2}{2h}+\alpha\frac{|\nabla\boldsymbol{w}|^2}{2}\right)d\boldsymbol{x}
    +\rho\sum^K_{i=1}|A^i-V^i_{\boldsymbol{w}}|^2,
\end{align}
where we use the extension \eqref{getband_cont} of the surface SDVF to $\Omega_\lambda$.
In \eqref{mboF}, $\alpha>0$ and $\rho>0$ are constants, and
\begin{equation}
    \begin{gathered}
        \label{cons_param_mbo}
        A^i=V^i_{z_0^{\epsilon,\lambda}},\quad V^i_{\boldsymbol{w}}=\int_{\Omega_\lambda}H^{\epsilon}(\phi^i_{\boldsymbol{w}}(\boldsymbol{x}))d\boldsymbol{x},\\
        \phi^i_{\boldsymbol{w}}(\boldsymbol{x})=\begin{cases}
            \displaystyle{\inf_{\boldsymbol{y}\in\partial Q^i_{\boldsymbol{w}}}}||\boldsymbol{x}-\boldsymbol{y}||_S, \quad \boldsymbol{x}\in Q^i_{\boldsymbol{w}}, \\
            -\displaystyle{\inf_{\boldsymbol{y}\in\partial Q^i_{\boldsymbol{w}}}}||\boldsymbol{x}-\boldsymbol{y}||_S, \quad \text{otherwise},
        \end{cases}
        \quad
        H^\epsilon(u)=\begin{cases}
            1,\quad u>\epsilon                                                                                         \\
            \frac{1}{2}+\frac{u}{2\epsilon}+\frac{1}{2\pi}\sin \frac{\pi u}{\epsilon}, \quad -\epsilon\le u\le\epsilon \\
            0,\quad u<-\epsilon                                                                                        \\
        \end{cases}\\
        Q^i_{\boldsymbol{w}}=\left\{ \boldsymbol{x}\in \Omega_\lambda|\boldsymbol{w}(\boldsymbol{x})\cdot \boldsymbol{p}_i\geq\boldsymbol{w}(\boldsymbol{x})\cdot \boldsymbol{p}_j,\text{ for all }j\in\{1,\cdots,K\}\right\},
    \end{gathered}
\end{equation}
where $z^{\epsilon,\lambda}_0$ is the extension of $z^\epsilon_0$ obtained from the initial region $\boldsymbol{P}_0$ on the surface $S$ to $\Omega_\lambda$.
The parameter $\rho$ controls the strength of the penalty. 
Note that when $\rho=0$, equation \eqref{mboF} takes the same form as the vectorized version of equation \eqref{mmheat_flat}, and applying the following algorithm results in a numerical solution of surface MCF without area preservation.

The surface MBO for realizing multiphase area-preserving surface MCF is as follows:
\begin{description}
    \setlength{\leftskip}{1.0cm}

    \item[1.] Using equation \eqref{sdvf}, create the surface SDVF $z_0^\epsilon$ from the initial domain $\boldsymbol{P}_0$.
    \item[2.] Extend $z_0^\epsilon$ to $\Omega_\lambda$ and define:
        $$
            z_0^{\epsilon,\lambda}(\boldsymbol{x})=z_0^{\epsilon}(C_S(\boldsymbol{x})), \quad \boldsymbol{x}\in\Omega_\lambda
        $$
    \item[3.] Using \eqref{cons_param_mbo}, obtain $A^i$ from $z_0^{\epsilon,\lambda}$.
    \item[4.] Repeat the following for $n=1,2,\cdots,N-1$:
        \begin{enumerate}[a.]
            \setlength{\leftskip}{1.0cm}
            \item Set $\boldsymbol{w}_0 = z^{\epsilon,\lambda}_{n-1}(\boldsymbol{x})$.
            \item For $m=1,\cdots,M$, find the minimizer $\boldsymbol{w}$ of the functional $\mathcal{F}_m(\boldsymbol{w})$ (refer to equation \eqref{mboF}). Denote the minimizer by $\boldsymbol{w}_{M}$.
            \item 
                  Extract $\boldsymbol{w}_M$ on the surface $S$ and denote it by $\boldsymbol{w}^S$:

                  $$ \boldsymbol{w}^S(\boldsymbol{x})=\boldsymbol{w}_M(\boldsymbol{x}),\quad\boldsymbol{x}\in S\cap\Omega_\lambda
                  $$
            \item 
                  Obtain $\boldsymbol{P}_{n+1}$ on the surface $S$ using $\boldsymbol{w}^S$ as follows:
                  \begin{align}
                      \boldsymbol{P}_{n+1}&=\bigcup_{i=1}^K\{P^i_{n+1}\},\notag\\
                      P^i_{n+1} &=\left\{ \boldsymbol{x}\in S|\boldsymbol{w}^S(\boldsymbol{x})\cdot \boldsymbol{p}_i\geq\boldsymbol{w}^S(\boldsymbol{x})\cdot \boldsymbol{p}_j,\text{ for all }j\in\{1,\cdots,K\}\right\}\notag
                  \end{align}
            \item Create the surface SDVF $z_{n+1}^\epsilon$ from $\boldsymbol{P}_{n+1}$ using equation \eqref{sdvf}.
            \item Extend $z_{n+1}^\epsilon$ to $\Omega_\lambda$ and denote it as $z_{n+1}^{\epsilon,\lambda}$ as follows:
                  \begin{align*}
                      z_{n+1}^{\epsilon,\lambda}(\boldsymbol{x})=z_{n+1}^\epsilon(C_S(\boldsymbol{x})), \quad\boldsymbol{x}\in\Omega_\lambda.
                  \end{align*}
        \end{enumerate}
\end{description}
The numerical examples using the approximation methods presented in this section are shown in Section \ref{results_MBO}. 
In the next section, we will explain an approximation methods for hyperbolic mean curvature flow on surfaces.
\subsection{Surface HMBO}\label{algo-HMBO}
In this section, we discuss an approximation method for surface HMCF on a surface $S$. 
First, we explain the surface HMBO for two-phase regions. Then, we describe its multiphase counterpart without area preservation. Afterwards, the surface HMBO for multiphase regions with area constraints is presented. 
The numerical results using the surface HMBO for multiphase regions is presented in Section \ref{results_HMBO}.

As before, let $\Omega_\lambda$ be a tubular region of $S$ where the distance from $S$ is within a constant value $\lambda$ as in Eq. \eqref{getband_cont}. 
The value of $\lambda$ is determined using the same method as described at  \eqref{getband}. 
The time step used in the implemented algorithm is denoted by $\tau=T/N$, where $N$ is a natural number representing the number of steps and $T>0$ is the final time.
\subsubsection{Surface HMBO for two-phase regions}\label{algo-2HMBO}
Below, we explain the method for approximating the motion of interfaces in a 2-phase region evolving by surface HMCF. 

Given an initial curve $\gamma^S_0$ and an initial velocity $v_0$, let  $\gamma^S_n$ represent the shape of the curve at time $n\tau$. 
Let $d_n(\boldsymbol{x})$ be a signed distance function on the surface to the interface $\gamma^S_n$. 
That is, if we let $P_n$ be the region enclosed by $\gamma_n^S$, $d_n(\boldsymbol{x})$ can be expressed by \eqref{sd}.
The surface HMBO for a 2-phase region is then described as follows:
\begin{description}
    \setlength{\leftskip}{1.0cm}

    \item[1.] Define $\gamma^S_1=\gamma^S_0+\tau v_0$ using the initial curve $\gamma^S_0$ and the initial velocity $v_0$.
    \item[2.] Create $d_0$ and $d_1$ from $\gamma^S_0$ and $\gamma^S_1$ using equation \eqref{sd}.
    \item[3.] Extend $d_0$ and $d_1$ to $\Omega_\lambda$ and denote them by $d^\lambda_0$ and $d^\lambda_1$ respectively, as follows:
        $$        d_l^\lambda(\boldsymbol{x})=d_l(C_S(\boldsymbol{x})),\quad\boldsymbol{x}\in\Omega_\lambda, \hspace{10pt}l=0,1
        $$
    \item[4.] Repeat the following for $n=1,2,\cdots,N-1$:
        \begin{enumerate}[a.]
            \setlength{\leftskip}{1.0cm}
            \item Solve the following wave equation: 
                  \begin{align}\label{algo-2HMBO-wave}
                      \begin{cases}
                          u_{tt}=\alpha\Delta u,                                                              \\
                          u(\boldsymbol{x},0)=2d^\lambda_{n}(\boldsymbol{x})-d^\lambda_{n-1}(\boldsymbol{x}), \\
                          u_t(\boldsymbol{x},0)=0
                      \end{cases}
                      \quad \boldsymbol{x}\in\Omega_\lambda,\quad \tau>t>0
                  \end{align}
                  where $\alpha>0$ is a constant.
            \item Define the new curve $\gamma^S_{n+1}$ as the zero level set of the solution $u(\boldsymbol{x},\tau)$ of equation \eqref{algo-2HMBO-wave} on the surface $S$:
                  $$
                      \gamma^S_{n+1}=\{\boldsymbol{x}\in\Omega_\lambda|S\cap\{u(\boldsymbol{x},\tau)=0\}\}
                  $$
            \item Create $d_{n+1}$ from $\gamma^S_{n+1}$ using equation \eqref{sd}.
            \item Extend $d_{n+1}$ to $\Omega_\lambda$ and define $d^\lambda_{n+1}$ as follows:
                  $$
                      d_{n+1}^\lambda(\boldsymbol{x})=d_{n+1}(C_S(\boldsymbol{x})), \quad\boldsymbol{x}\in\Omega_\lambda.
                  $$
        \end{enumerate}
\end{description}
The results of the numerical error analysis using the above algorithm are presented in Section \ref{2相領域に対するSurface HMBOアルゴリズムの数値誤差解析}. 
In the above algorithm, the signed distance function is used to detect the location of the interface. 
As before, by using the signed distance vector field instead of the signed distance function, we can implement the surface HMBO in the multiphase setting. This will be explained in the next section.
\subsubsection{Surface HMBO for multiphase regions}\label{algo-mHMBO}
In the following, we describe a method for approximating surface HMCF of  interfaces separating $K$ multiphase regions on a surface $S$. 
To this end, let $P^i_n$ represent the region of phase $i$ at time $n\tau$, and $\boldsymbol{P}_n=\bigcup_{i=1}^K{P^i_n}$. 
In addition to the initial regions $P_0^i$, we also provide the initial velocities $\boldsymbol{v}_0^i$ for each phase $i$. 
We denote the vector given to phase $i$ as $\boldsymbol{p}_i$, as explained in Section \ref{曲面版符号付き距離ベクトル場}. 
Again, $z_n^\epsilon(\boldsymbol{x})$ denotes surface SDVF defined by equation \eqref{sdvf} and which is constructed from $\boldsymbol{P}_n$. 

The surface HMBO algorithm for multiphase regions is as follows:
\begin{description}
    \setlength{\leftskip}{1.0cm}

    \item[1.] Determine $\boldsymbol{P}_1$ from the initial domain $\boldsymbol{P}_0$ and the initial velocities $\boldsymbol{v}_0^i$ on $\partial P_0^i$. For details, refer to Section \ref{多相領域に対するSurface HMBOアルゴリズムと初期速度}.
    \item[2.] Using equation \eqref{sdvf}, create the surface SDVFs  $z_0^\epsilon$, and $z_1^\epsilon$, from $\boldsymbol{P}_0,\boldsymbol{P}_1$.
    \item[3.] Extend $z_0^\epsilon$, and $z_1^\epsilon$ to $\Omega_\lambda$, and denote their extensions by $z_0^{\epsilon,\lambda}$ and $z_1^{\epsilon,\lambda}$ respectively, as follows:
        $$
            z_l^{\epsilon,\lambda}(\boldsymbol{x})=z_l^{\epsilon}(C_S(\boldsymbol{x})), \quad \boldsymbol{x}\in\Omega_\lambda,\quad l=0,1
        $$
    \item[4.] Repeat the following for $n=1,2,\cdots,N-1$:
        \begin{enumerate}[a.]
            \setlength{\leftskip}{1.0cm}
            \item Solve the vectorial wave equation: 
                  \begin{align}\label{wave-vec}
                      \begin{cases}
                          \boldsymbol{u}_{tt}=\alpha\Delta \boldsymbol{u},           
                          \\
                          \boldsymbol{u}(\boldsymbol{x},0)=2z^{\epsilon,\lambda}_{n}(\boldsymbol{x})-z^{\epsilon,\lambda}_{n-1}(\boldsymbol{x}), \\
                          \boldsymbol{u}_t(\boldsymbol{x},0)=\boldsymbol{0},
                      \end{cases}
                      \quad \boldsymbol{x}\in\Omega_\lambda,\quad \tau>t>0
                  \end{align}
                  where, $\alpha>0$ is a constant.

            \item Extract the solution $\boldsymbol{u}(\boldsymbol{x},\tau)$ of equation \eqref{wave-vec} on the surface $S$ and denote it by  $\boldsymbol{u}^S$:

                  $$
                      \boldsymbol{u}^S(\boldsymbol{x})=\boldsymbol{u}(\boldsymbol{x},\tau),\quad \boldsymbol{x}\in S\cap\Omega_\lambda.
                  $$
            \item Obtain $\boldsymbol{P}_{n+1}$ on the surface $S$ using $\boldsymbol{u}^S$ as follows:
                  \begin{align}
                      \boldsymbol{P}_{n+1}&=\bigcup_{i=1}^K\{P^i_{n+1}\},\notag\\
                      P^i_{n+1} &=\left\{ \boldsymbol{x}\in S|\boldsymbol{u}^S(\boldsymbol{x})\cdot \boldsymbol{p}_i\geq\boldsymbol{u}^S(\boldsymbol{x})\cdot \boldsymbol{p}_j,\text{ for all }j\in\{1,\cdots,K\}\right\}\notag
                  \end{align}
            \item Create the surface SDVF $z_{n+1}^\epsilon$ from $\boldsymbol{P}_{n+1}$ using equation \eqref{sdvf}.
            \item Extend $z_{n+1}^\epsilon$ to $\Omega_\lambda$ and denote it by $z_{n+1}^{\epsilon,\lambda}$ as follows:
                  \begin{align*}
                      z_{n+1}^{\epsilon,\lambda}(\boldsymbol{x})=z_{n+1}^\epsilon(C_S(\boldsymbol{x})), \quad\boldsymbol{x}\in\Omega_\lambda
                  \end{align*}
        \end{enumerate}
\end{description}
In the next section, we show how MM can be combined with the above surface HMBO algorithm to realize multiphase surface HMCF of interfaces with optional area-preserving conditions.
\subsubsection{Surface HMBO for multiphase  area-preserving motions}\label{algo-shmboc}
In the following, we describe our method for approximating multiphase surface HMCF, where the area of each domain's preserved. 
We assume that the initial conditions for each phase $i$ are given by the initial domain $P_0^i$ and the initial velocity $\boldsymbol{v}_0^i$. 
As in section \ref{曲面版符号付き距離ベクトル場}, we denote the vector given to phase $i$ by $\boldsymbol{p}_i$.

Similar to section \ref{algo-smboc}, MM are used to incorporate the area-preserving conditions. We take a sufficiently large integer $M$ and set $h=\tau/M$. 
Let the functions $\boldsymbol{w}_{m-1}$ and $\boldsymbol{w}_{m-2}$ be defined using the surface SDVF extended to $\Omega_\lambda$ and the MM functional be as follows:

\begin{align}\label{F}   \mathcal{F}_m(\boldsymbol{w})=\int_{\Omega_\lambda}\left(\frac{|\boldsymbol{w}-2\boldsymbol{w}_{m-1}+\boldsymbol{w}_{m-2}|^2}{2h^2}+\alpha\frac{|\nabla\boldsymbol{w}|^2}{2}\right)d\boldsymbol{x}
    +\rho\sum^K_{i=1}|A^i-V^i_{\boldsymbol{w}}|^2
\end{align}
Here, $\alpha>0$ and $\rho>0$ are constants, and $A^i$ and $V_{\boldsymbol{w}}^i$ are defined by \eqref{cons_param_mbo}.

The definitions of $\boldsymbol{w}_1$ and $\boldsymbol{w}_{0}$ are as described before. Namely, we use $z^{\epsilon,\lambda}_0$ (the extension of $z^\epsilon_0$ obtained from the initial domain $\boldsymbol{P}_0$) together with the initial velocities $\boldsymbol{v}^i_0$  on the surface $S$. 
Here, $\rho>0$ is a sufficiently large penalty parameter. 
Note that when $\rho=0$, the functional \eqref{F} is the same as the vectorized extension of equation \eqref{mmwave_flat}, and that applying the algorithm below would result in a numerical approximation of surface HMCF without area preservation. 

The Surface HMBO that realizes area preservation for multiphase domains is as follows.

\begin{description}
    \setlength{\leftskip}{1.0cm}
     \item[1.] 
        Determine $\boldsymbol{P}_1$ from the initial domain $\boldsymbol{P}_0$ and initial velocities $\boldsymbol{v}_0^i$ on $\partial P_0^i$. For details, refer to Section \ref{多相領域に対するSurface HMBOアルゴリズムと初期速度}.
    \item[2.] 
        Using equation \eqref{sdvf}, create the surface SDVF $z_0^\epsilon$ and $z_1^\epsilon$ from $\boldsymbol{P}_0$ and $\boldsymbol{P}_1$.

    \item[3.] Define the extensions of $z_0^\epsilon$ and $z_1^\epsilon$ to $\Omega_\lambda$  as follows:
        $$
            z_l^{\epsilon,\lambda}(\boldsymbol{x})=z_l^{\epsilon}(C_S(\boldsymbol{x})), \quad \boldsymbol{x}\in\Omega_\lambda,\quad l=0,1
        $$
    \item[4.] 
    Using equation \eqref{cons_param_mbo}, compute each $A^i$ from $z_0^{\epsilon,\lambda}$.
    \item[5.] Repeat the following for $n=1,2,\cdots,N-1$:
        \begin{enumerate}[a.]
            \setlength{\leftskip}{1.0cm}
            \item Set $\boldsymbol{w}_0=\boldsymbol{w}_1=2z^{\epsilon,\lambda}_{n}(\boldsymbol{x})-z^{\epsilon,\lambda}_{n-1}(\boldsymbol{x})$.
            \item For each $m=2,\cdots,M$, minimize $\mathcal{F}_m(\boldsymbol{w})$ (given by \eqref{F}) and denote each minimizer by $\boldsymbol{w}_{m}$.
                  
            \item Let $\boldsymbol{w}^S$ denote the restriction of  $\boldsymbol{w}_M$ to the surface $S$ as follows:

                  $$             \boldsymbol{w}^S(\boldsymbol{x})=\boldsymbol{w}_M(\boldsymbol{x}),\quad\boldsymbol{x}\in S\cap\Omega_\lambda
                  $$
            \item Obtain $\boldsymbol{P}_{n+1}$ on the surface $S$ using $\boldsymbol{w}^S$ as follows:
                  \begin{align}
                      \boldsymbol{P}_{n+1}&=\bigcup_{i=1}^K\{P^i_{n+1}\}\notag\\ P^i_{n+1} &=\left\{ \boldsymbol{x}\in S|\boldsymbol{w}^S(\boldsymbol{x})\cdot \boldsymbol{p}_i\geq\boldsymbol{w}^S(\boldsymbol{x})\cdot \boldsymbol{p}_j,\text{ for all }j\in\{1,\cdots,K\}\right\}\notag
                  \end{align}
            \item Create Surface SDVF $z_{n+1}^\epsilon$ from $\boldsymbol{P}_{n+1}$ using Equation \eqref{sdvf}.
            \item Let $z_{n+1}^{\epsilon,\lambda}$ be the extension of  $z_{n+1}^\epsilon$ to $\Omega_\lambda$:
                  \begin{align*}
                      z_{n+1}^{\epsilon,\lambda}(\boldsymbol{x})=z_{n+1}^\epsilon(C_S(\boldsymbol{x})), \quad\boldsymbol{x}\in\Omega_\lambda.
                  \end{align*}
        \end{enumerate}
\end{description}
Numerical examples using the method described in this section are presented in section \ref{results_HMBO}.
\section{Numerical results and considerations}\label{数値計算結果と考察}

In this section, we use the approximation methods presented in sections \ref{algo-MBO} and \ref{algo-HMBO} to numerically solve the mean curvature flow and hyperbolic mean curvature flow on surfaces under various conditions. 

The discrete approximation of $\Omega_\lambda$ uses a uniformly spaced orthogonal grid $\Omega_\lambda^D$ with a spacing of $\Delta x$ in all three directions.
Note that $\Omega_\lambda^D$ is obtained using the same method as in section \ref{アルゴリズム}. 
Details regarding the approximation of the MM functional values are explained in section \ref{実装方法について}. 
In all cases, the interpolation parameter $\epsilon$ used to construct the surface SDVF represented by Equation \eqref{sdvf} is set to $\epsilon=0.03$.\\
{\bf{Remark:}}
The interpolation parameter $\epsilon$ needs to be appropriately selected depending on the discretization of the surface $S$. 
In the case that the surface is discretized by a point cloud, it was found from the numerical investigations in this study that a value of 3-5 times the average distance to neighboring points within the point cloud is appropriate for the interpolation parameter $\epsilon$.
\subsection{Regarding the initial conditions}\label{sp_ini}
Boundaries between regions determine the shape of the interface and hence the initial conditions used in the numerical calculations. 
The following two types of initial conditions for the numerical calculations on the unit sphere were used. In both cases, points with different colors indicate different phases.
\begin{enumerate}
    \item Two-phases on the unit sphere
          \begin{figure}[H]
              \centering
              \includegraphics[bb=0cm 0cm 10cm 5cm,scale=0.7]{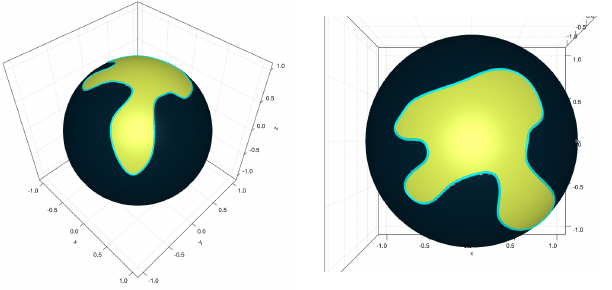}
              \caption{Initial Condition (Two-Phase)}
              \label{ini_sp_phase2}
              \centerline{(Left: Figure viewed from an oblique angle. Right: Figure viewed from directly above.)}
          \end{figure}
    \item Four-phases on the unit sphere
          \begin{figure}[H]
              \centering
              \includegraphics[bb=0cm 0cm 10cm 5cm,scale=0.7]{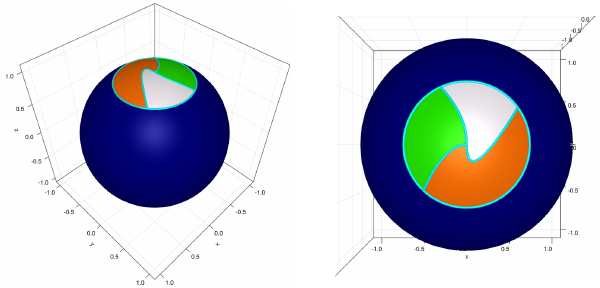}
              \caption{Initial Condition (Four-Phase)}
              \label{ini_sp_phase4}
          \end{figure}
          \vspace{-12pt}
          \centerline{(Left: Figure viewed from an oblique angle. Right: Figure viewed from directly above.)}
\end{enumerate}
\subsection{Computational details regarding surface mean curvature flow}\label{results_MBO}
In the following, we introduce the parameters used in our computations of surface mean curvature flow. 
We performed computations of the surface  mean curvature flow involving interfaces in two and  four phase environments on the unit sphere, with and without area preservation.

\begin{figure}[H]
    \centering
    \fbox{
        \begin{tabular}{cc}
            \begin{minipage}[t]{0.25\columnwidth}
                \centering
                \includegraphics[height=0.15\vsize]{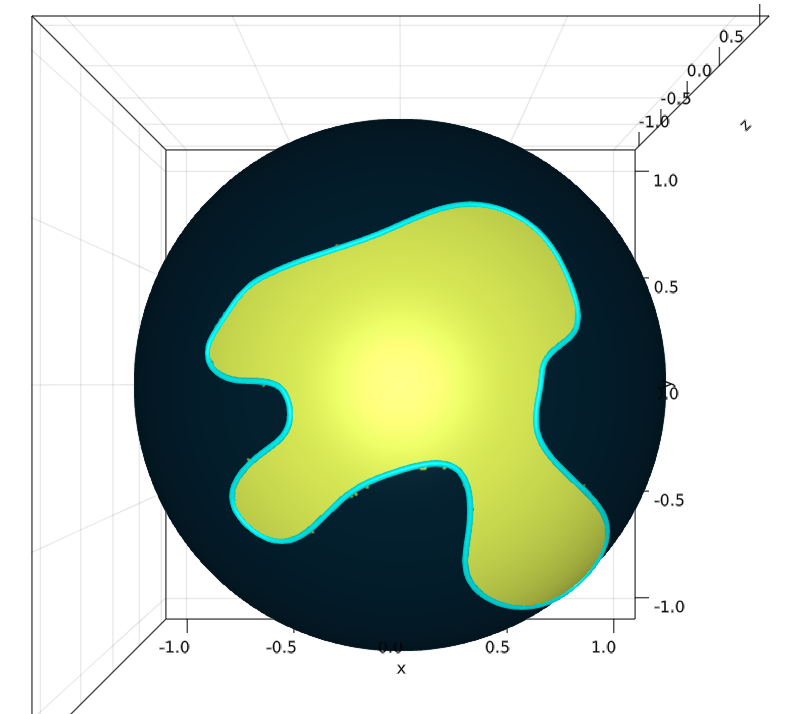}
                \subcaption{}

            \end{minipage} &
            \begin{minipage}[t]{0.25\columnwidth}
                \centering
                \includegraphics[height=0.15\vsize]{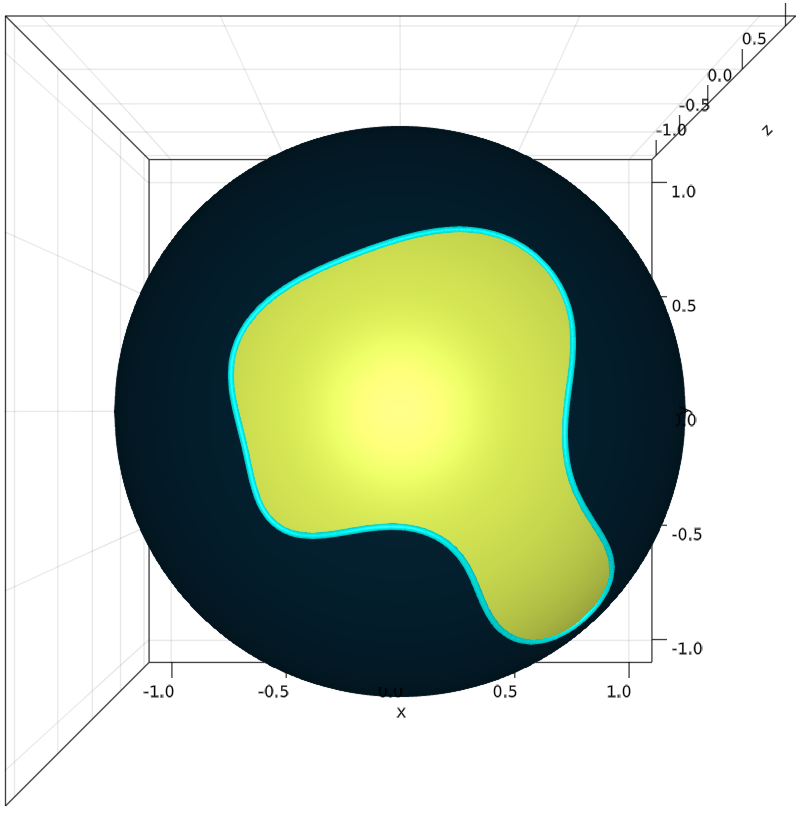}
                \subcaption{}
                \label{}
            \end{minipage}
            \\
            \begin{minipage}[t]{0.25\columnwidth}
                \centering
                \includegraphics[height=0.15\vsize]{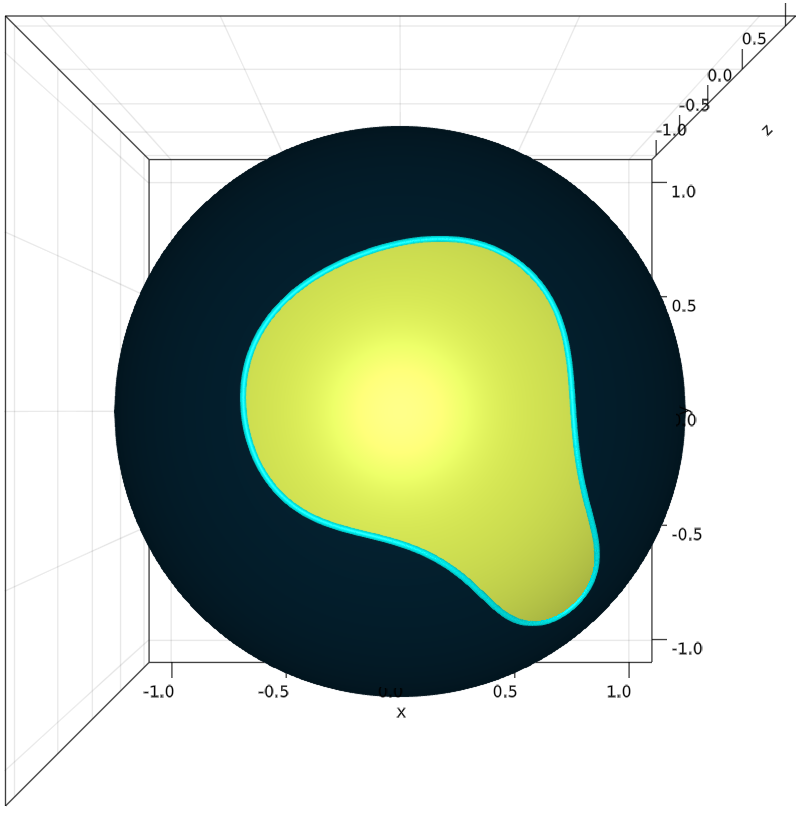}
                \subcaption{}

            \end{minipage} &
            \begin{minipage}[t]{0.25\columnwidth}
                \centering
                \includegraphics[height=0.15\vsize]{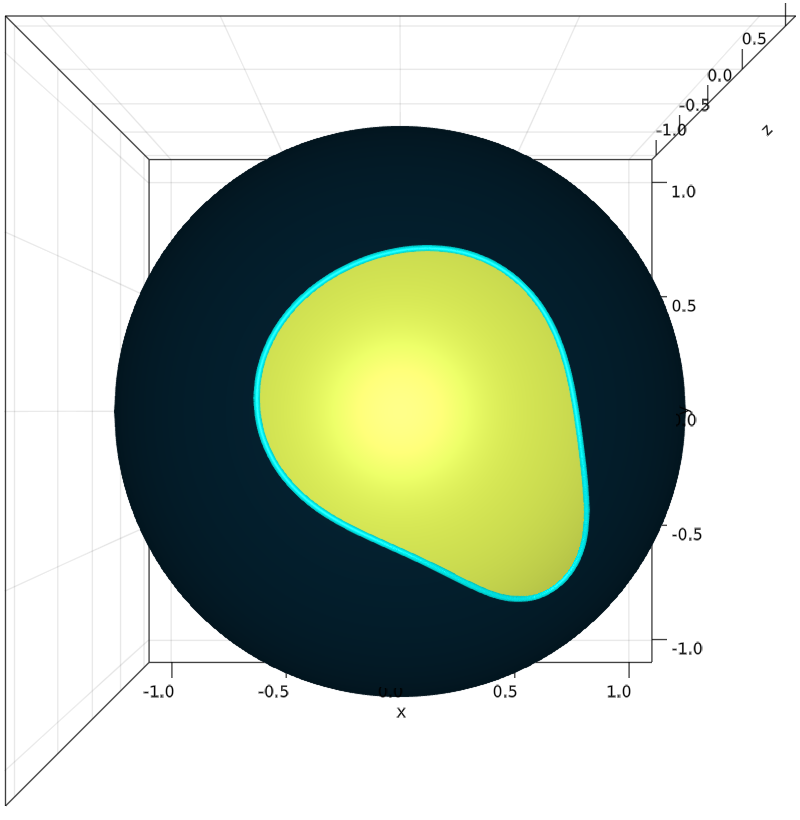}
                \subcaption{}

            \end{minipage}
            \\
            \begin{minipage}[t]{0.25\columnwidth}
                \centering
                \includegraphics[height=0.15\vsize]{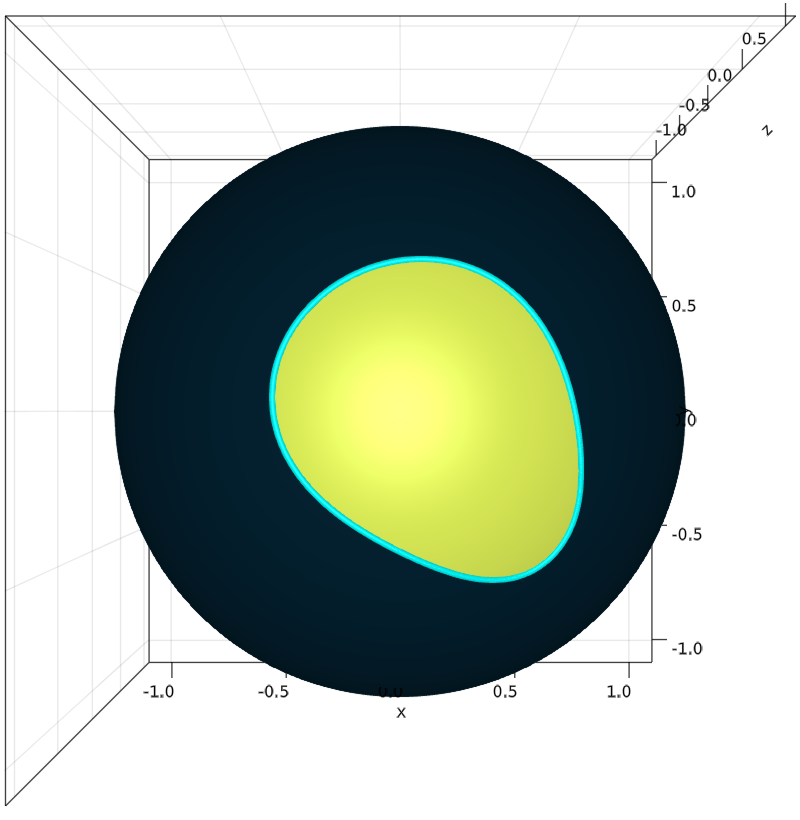}
                \subcaption{}

            \end{minipage} &
            \begin{minipage}[t]{0.25\columnwidth}
                \centering
                \includegraphics[height=0.15\vsize]{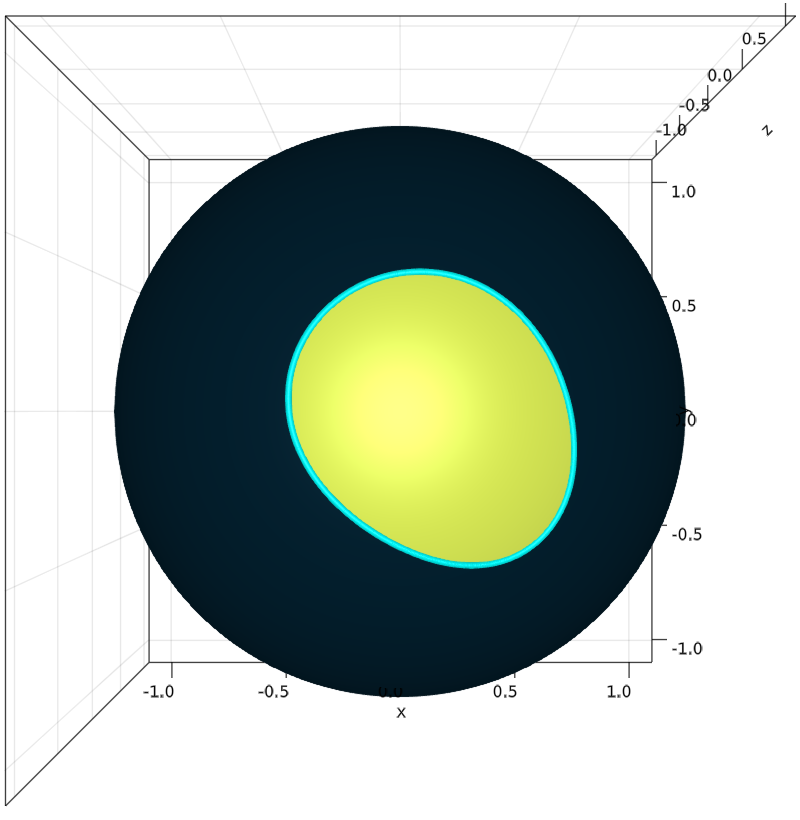}
                \subcaption{}

            \end{minipage}
            \\
            \begin{minipage}[t]{0.25\columnwidth}
                \centering
                \includegraphics[height=0.15\vsize]{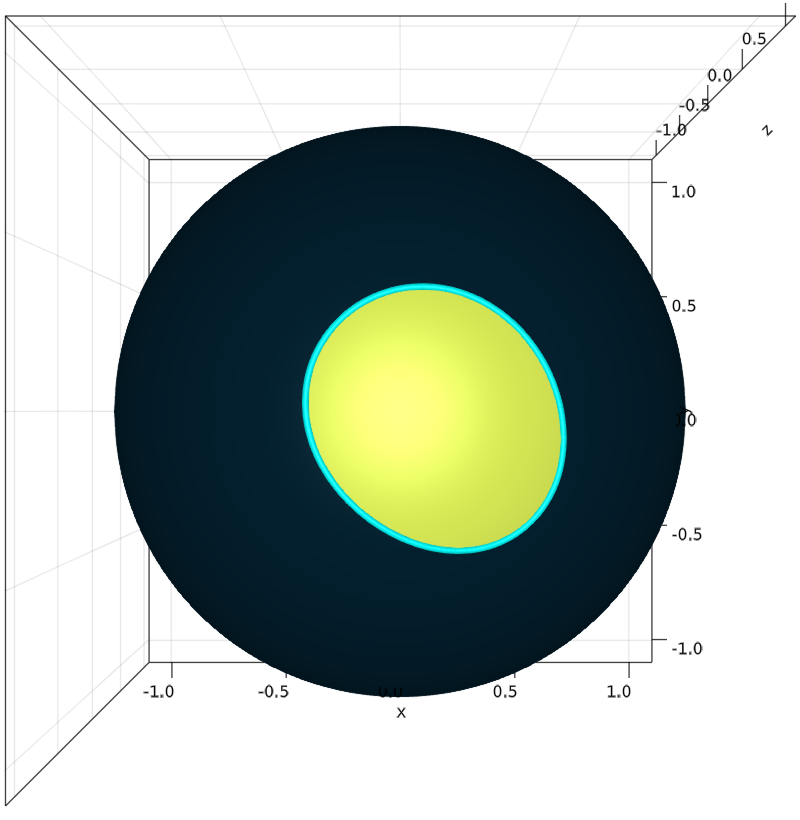}
                \subcaption{}

            \end{minipage} &
            \begin{minipage}[t]{0.25\columnwidth}
                \centering
                \includegraphics[height=0.15\vsize]{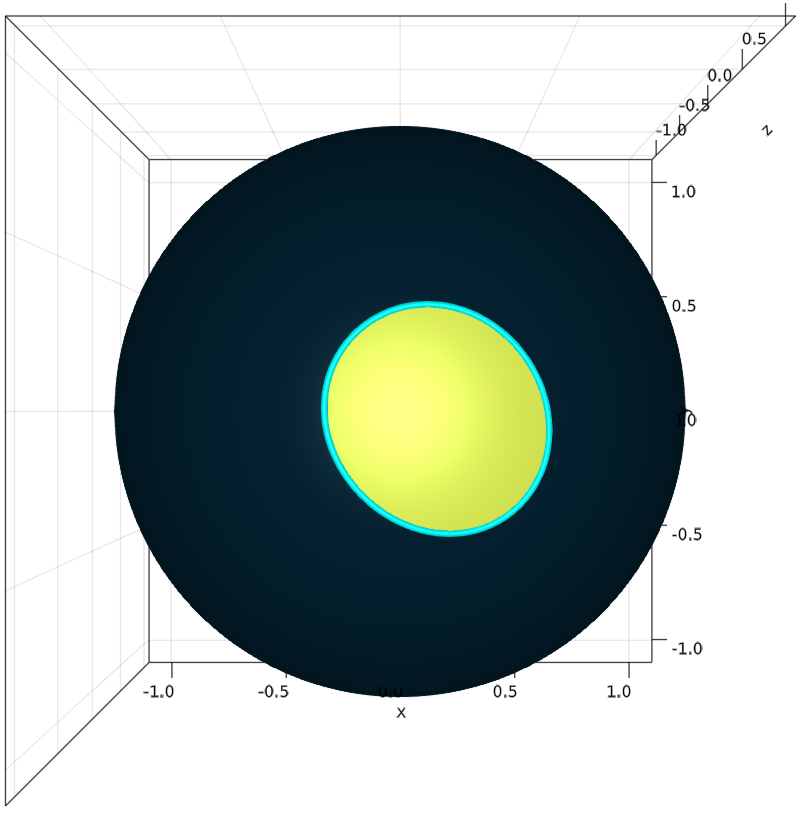}
                \subcaption{}

            \end{minipage}
            \\
            \begin{minipage}[t]{0.25\columnwidth}
                \centering
                \includegraphics[height=0.15\vsize]{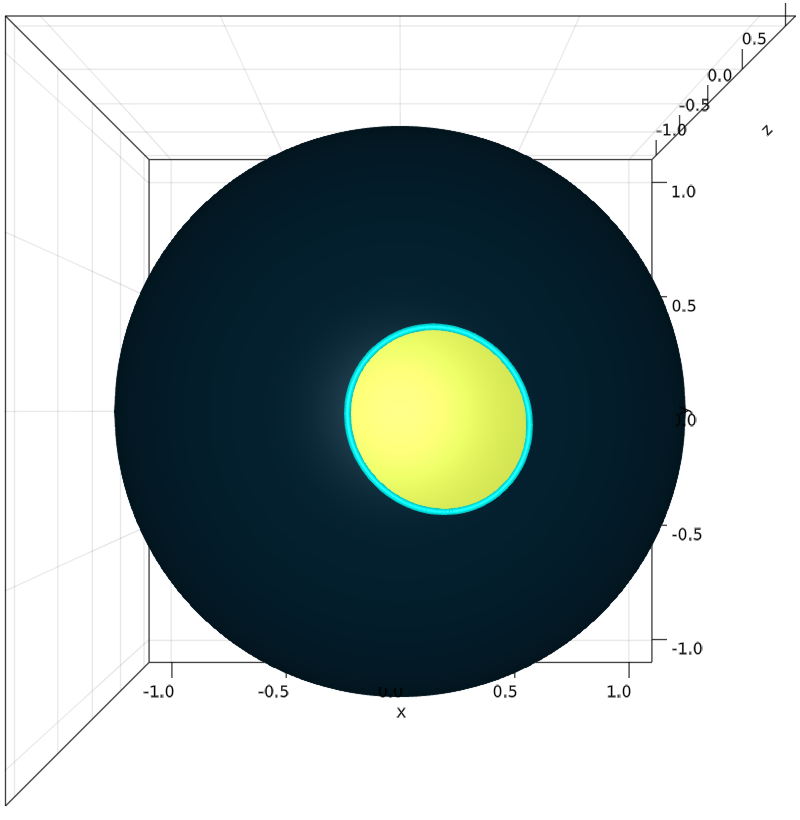}
                \subcaption{}

            \end{minipage} &
            \begin{minipage}[t]{0.25\columnwidth}
                \centering
                \includegraphics[height=0.15\vsize]{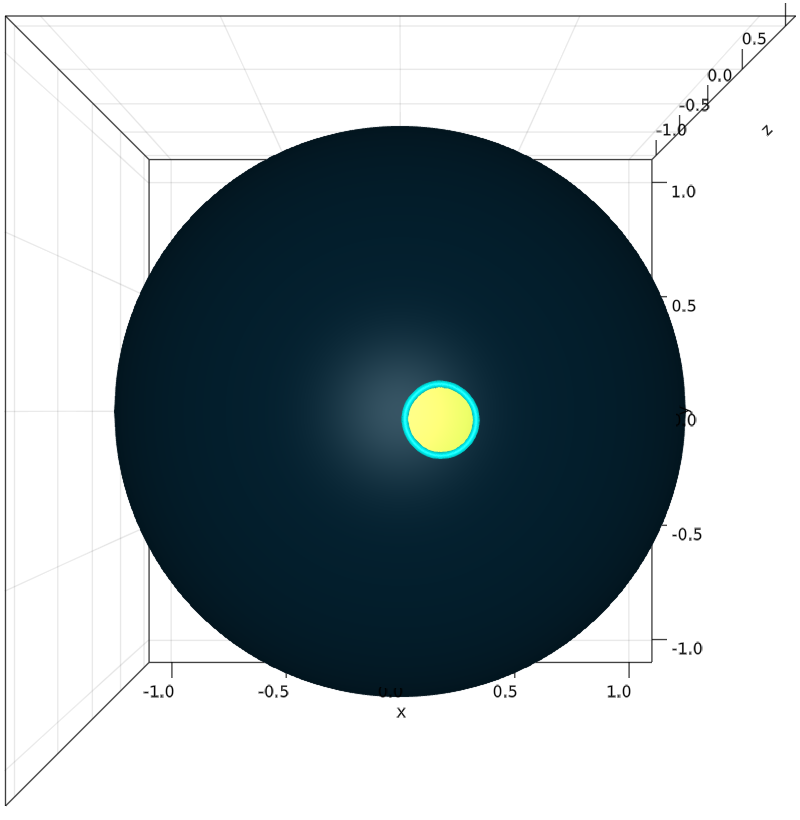}
                \subcaption{}

            \end{minipage}
        \end{tabular}
    }
    \caption{Result of two-phase MCF on the unit sphere(without area preservation). Arranged in alphabetical order with equal intervals from the initial time to the time the interface disappears.}
    \label{res:sp2_mcf}
\end{figure}

\begin{figure}[H]
    \centering
    \fbox{
        \begin{tabular}{cc}
            \begin{minipage}[t]{0.25\columnwidth}
                \centering
                \includegraphics[height=0.15\vsize]{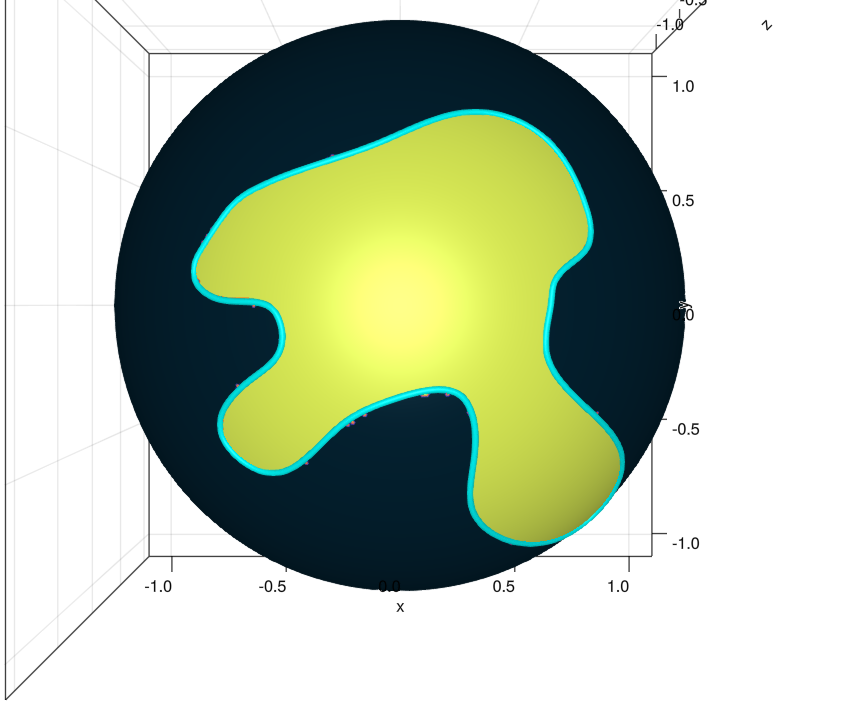}
                \subcaption{}

            \end{minipage} &
            \begin{minipage}[t]{0.25\columnwidth}
                \centering
                \includegraphics[height=0.15\vsize]{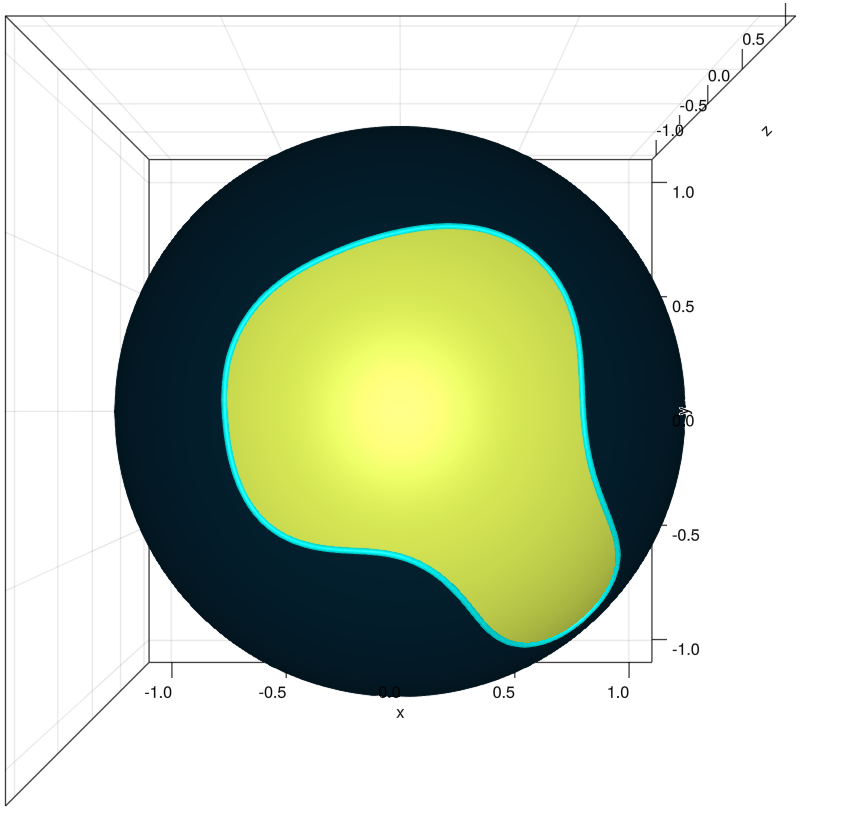}
                \subcaption{}
                \label{}
            \end{minipage}
            \\
            \begin{minipage}[t]{0.25\columnwidth}
                \centering
                \includegraphics[height=0.15\vsize]{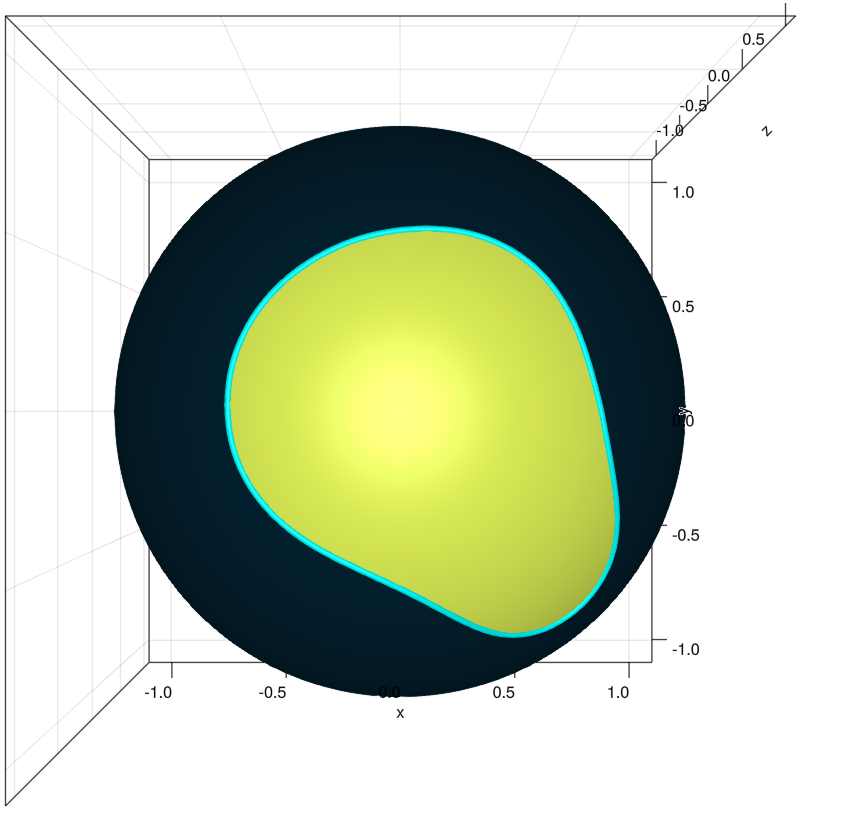}
                \subcaption{}

            \end{minipage} &
            \begin{minipage}[t]{0.25\columnwidth}
                \centering
                \includegraphics[height=0.15\vsize]{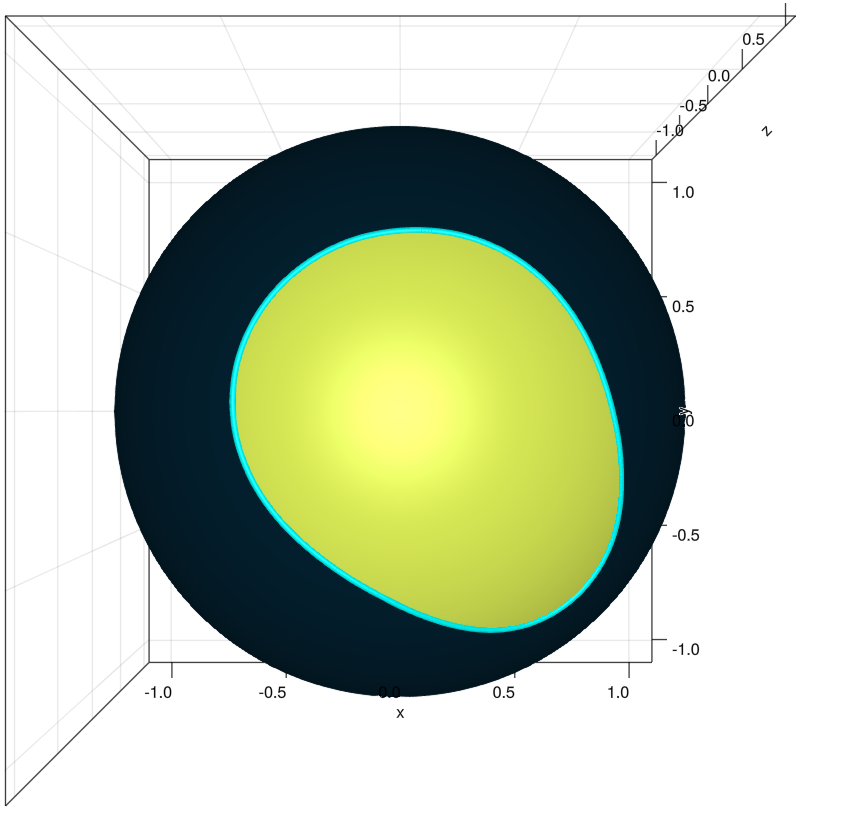}
                \subcaption{}

            \end{minipage}
            \\
            \begin{minipage}[t]{0.25\columnwidth}
                \centering
                \includegraphics[height=0.15\vsize]{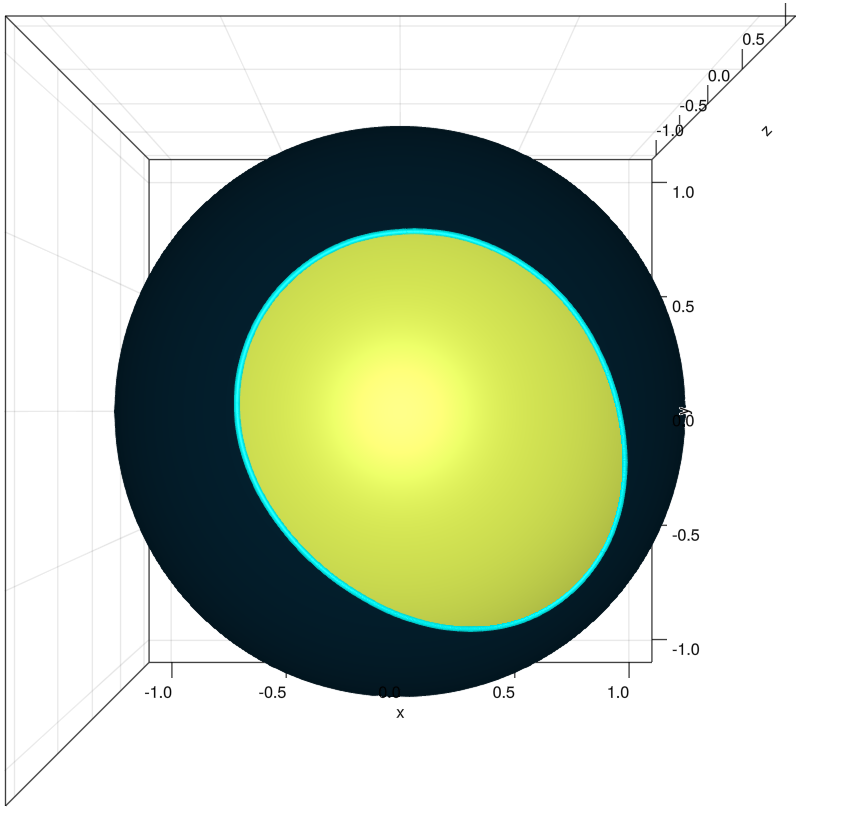}
                \subcaption{}

            \end{minipage} &
            \begin{minipage}[t]{0.25\columnwidth}
                \centering
                \includegraphics[height=0.15\vsize]{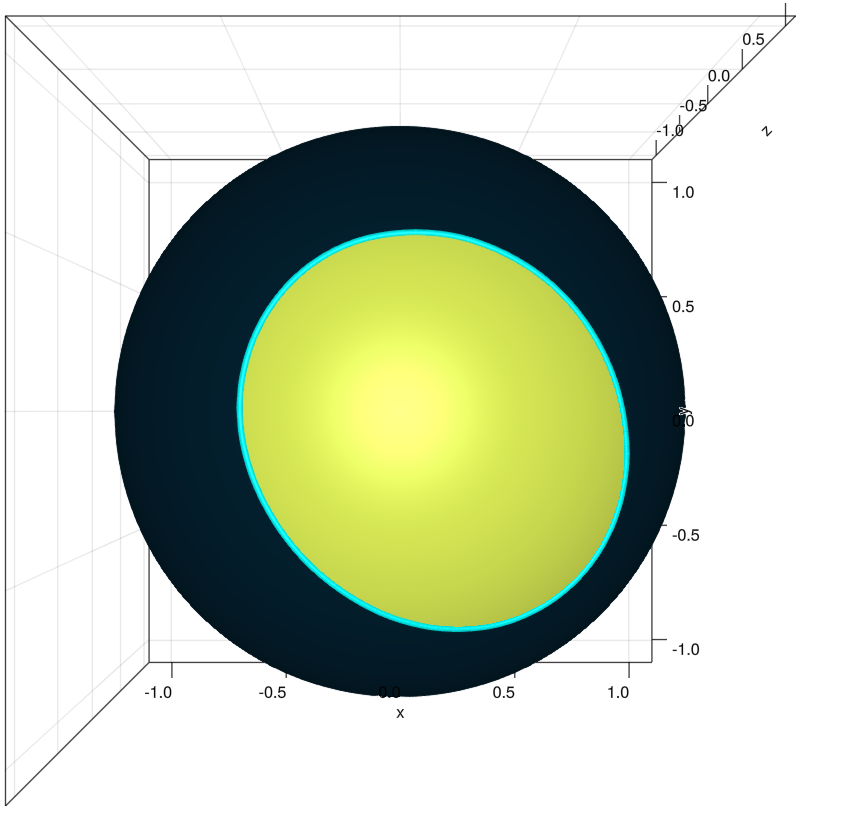}
                \subcaption{}

            \end{minipage}
            \\
            \begin{minipage}[t]{0.25\columnwidth}
                \centering
                \includegraphics[height=0.15\vsize]{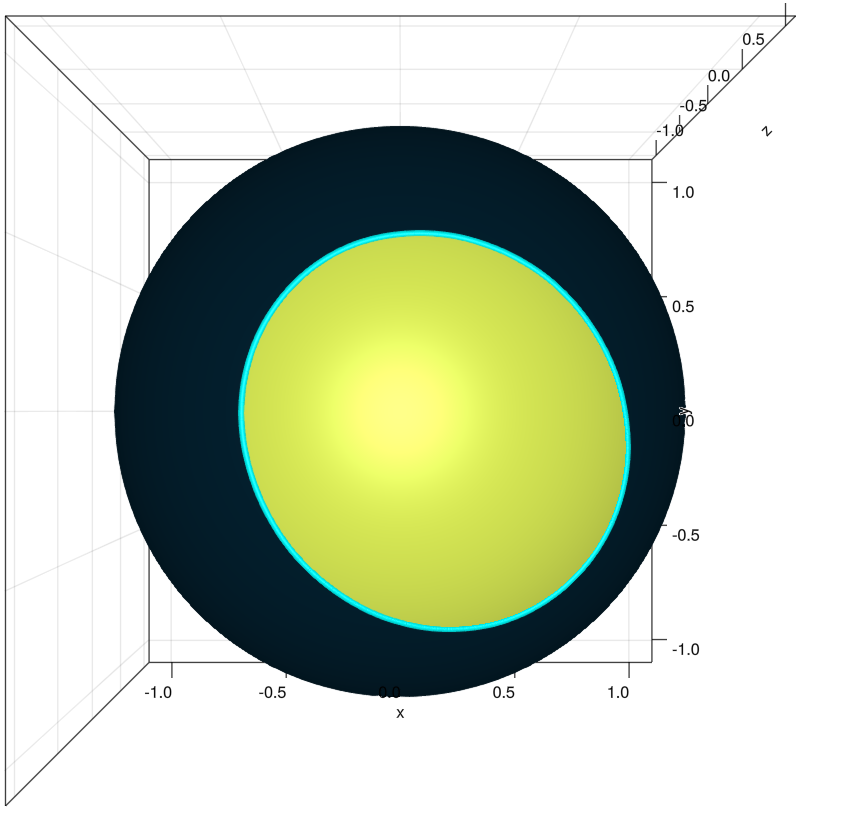}
                \subcaption{}

            \end{minipage} &
            \begin{minipage}[t]{0.25\columnwidth}
                \centering
                \includegraphics[height=0.15\vsize]{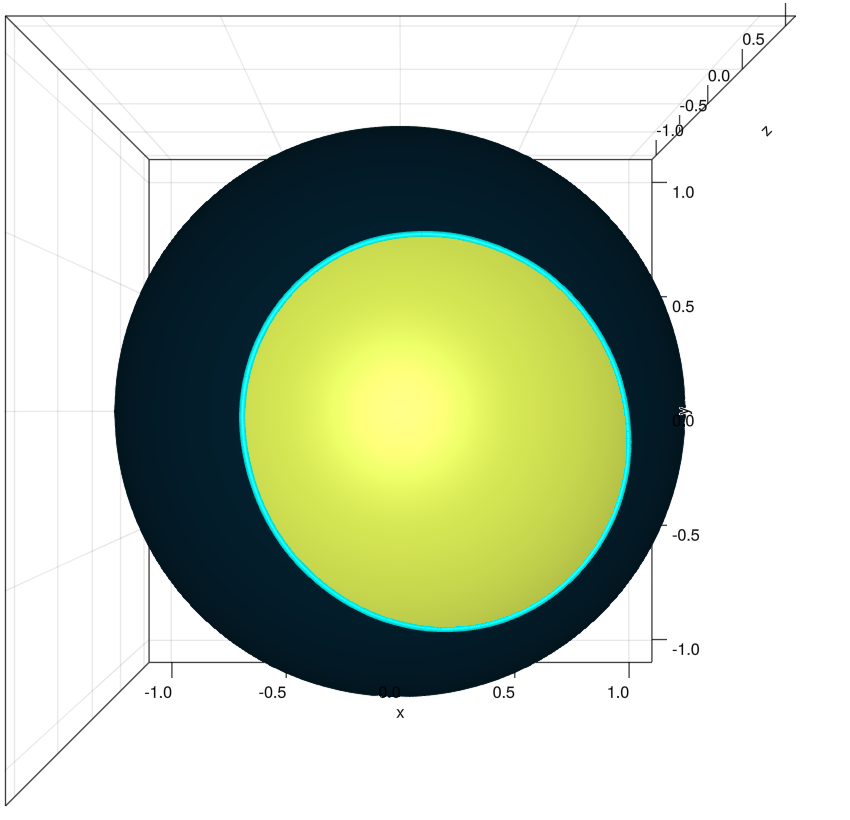}
                \subcaption{}

            \end{minipage}
            \\
            \begin{minipage}[t]{0.25\columnwidth}
                \centering
                \includegraphics[height=0.15\vsize]{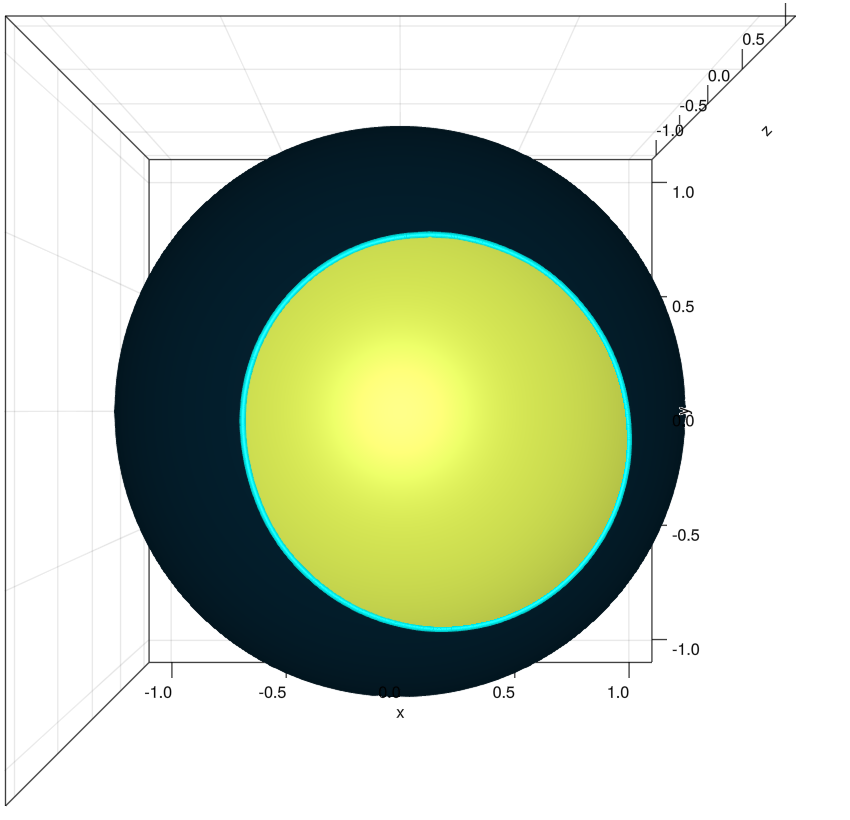}
                \subcaption{}

            \end{minipage} &
            \begin{minipage}[t]{0.25\columnwidth}
                \centering
                \includegraphics[height=0.15\vsize]{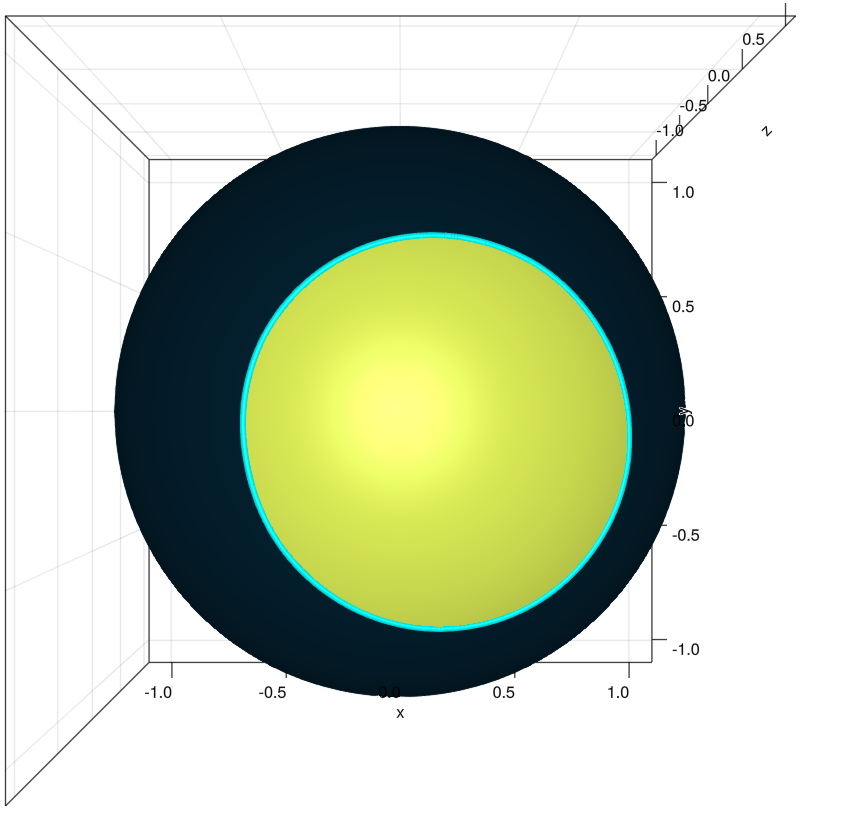}
                \subcaption{}

            \end{minipage}
        \end{tabular}
    }
    \caption{Result of two-phase MCF on the unit sphere (with area preservation). Arranged in alphabetical order. The time intervals from the initial time to the time that the interface reaches a nearly stationary state are equally spaced.}
    \label{res:sp2_mcf_cons}
\end{figure}

\begin{figure}[H]
    \centering
    \fbox{
        \begin{tabular}{cc}
            \begin{minipage}[t]{0.25\columnwidth}
                \centering
                \includegraphics[height=0.15\vsize]{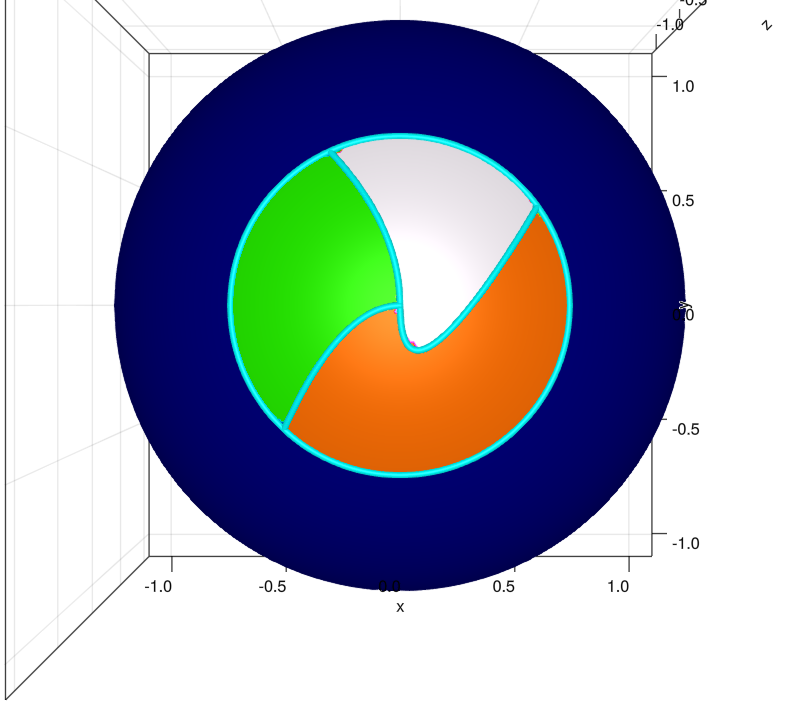}
                \subcaption{}

            \end{minipage} &
            \begin{minipage}[t]{0.25\columnwidth}
                \centering
                \includegraphics[height=0.15\vsize]{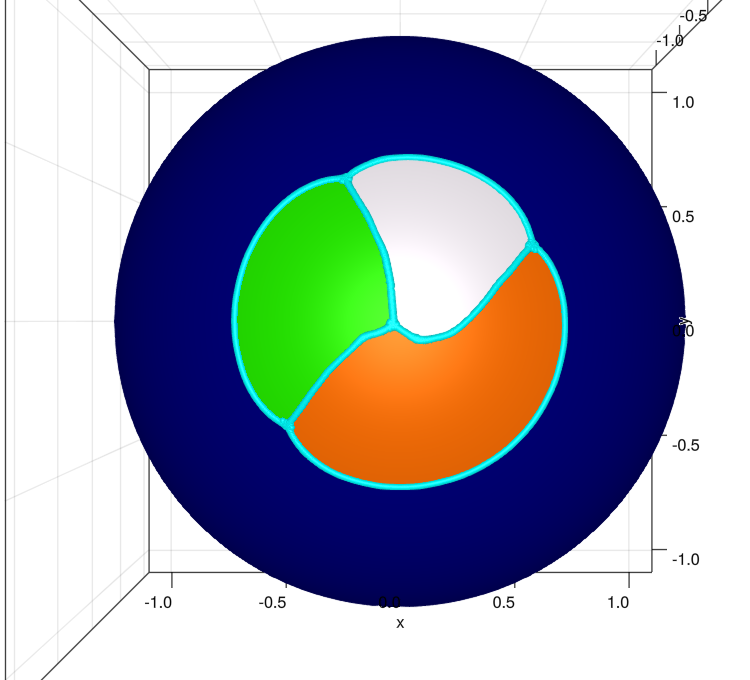}
                \subcaption{}
                \label{}
            \end{minipage}
            \\
            \begin{minipage}[t]{0.25\columnwidth}
                \centering
                \includegraphics[height=0.15\vsize]{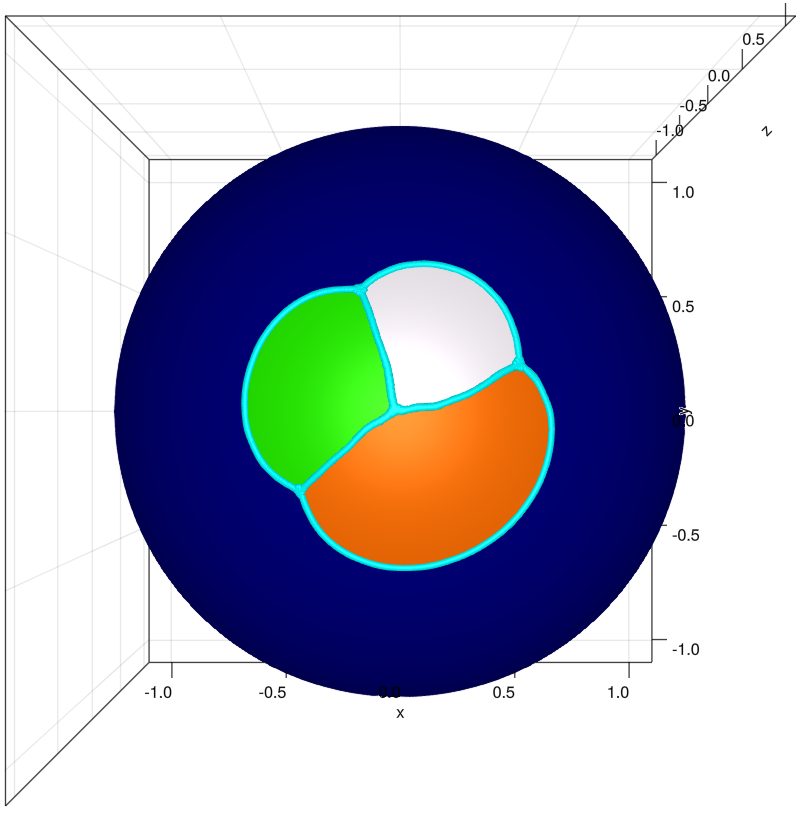}
                \subcaption{}

            \end{minipage} &
            \begin{minipage}[t]{0.25\columnwidth}
                \centering
                \includegraphics[height=0.15\vsize]{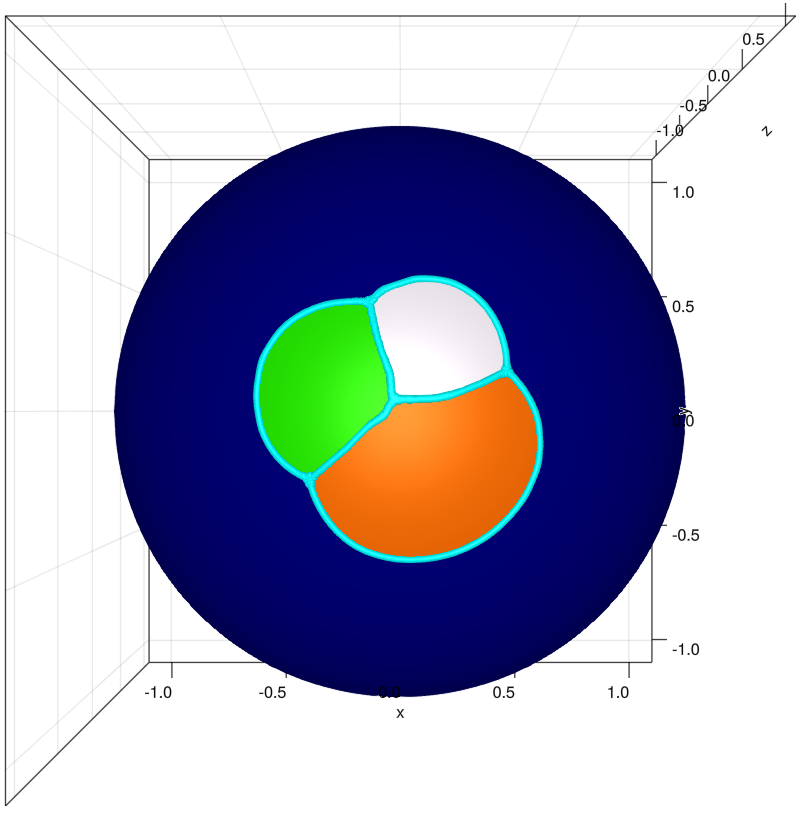}
                \subcaption{}

            \end{minipage}
            \\
            \begin{minipage}[t]{0.25\columnwidth}
                \centering
                \includegraphics[height=0.15\vsize]{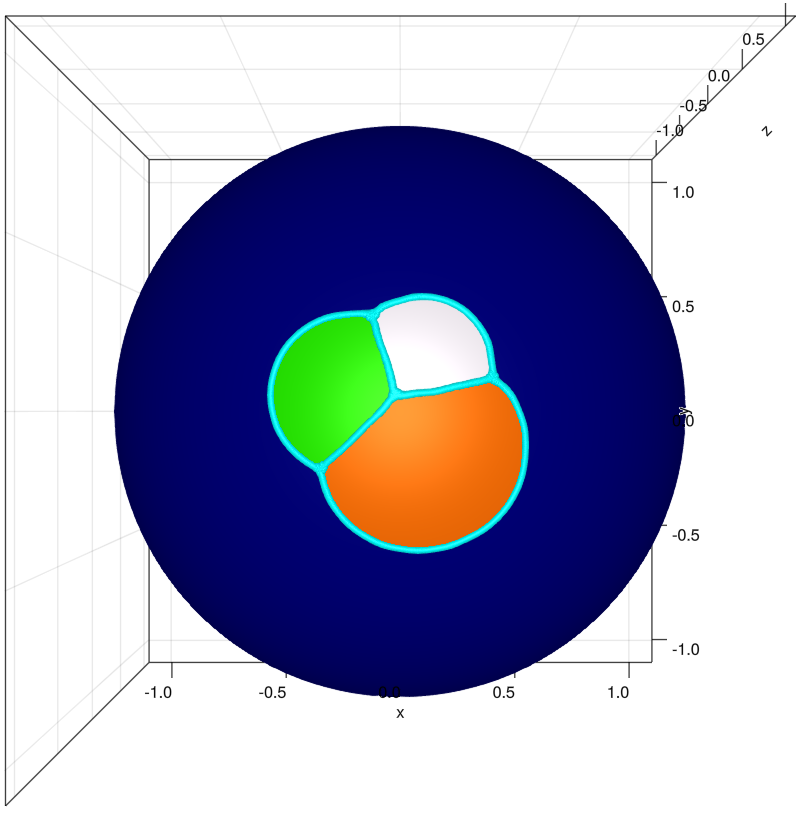}
                \subcaption{}

            \end{minipage} &
            \begin{minipage}[t]{0.25\columnwidth}
                \centering
                \includegraphics[height=0.15\vsize]{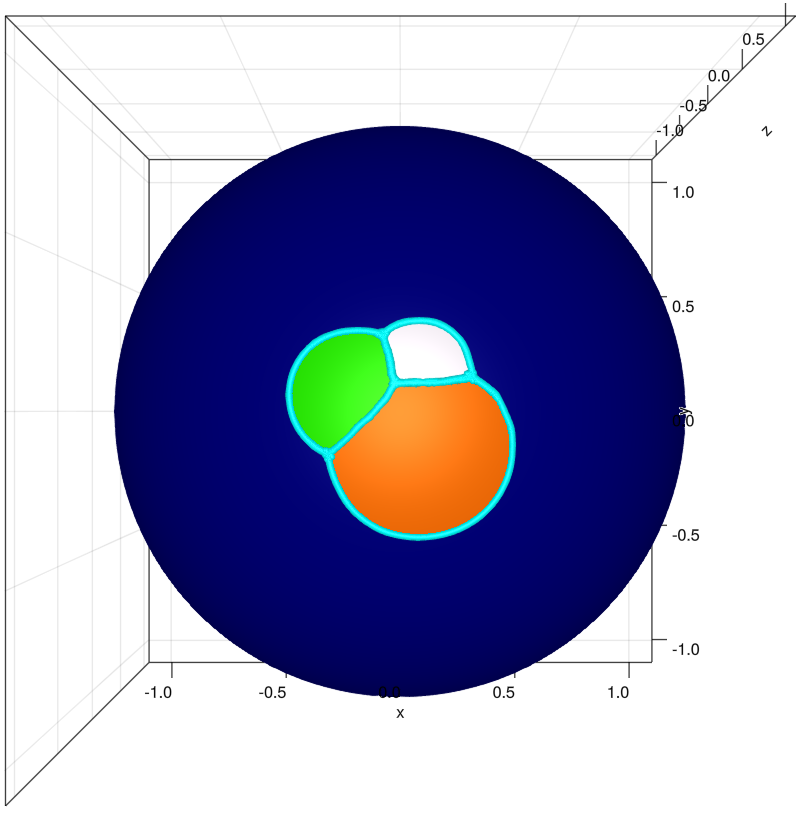}
                \subcaption{}

            \end{minipage}
            \\
            \begin{minipage}[t]{0.25\columnwidth}
                \centering
                \includegraphics[height=0.15\vsize]{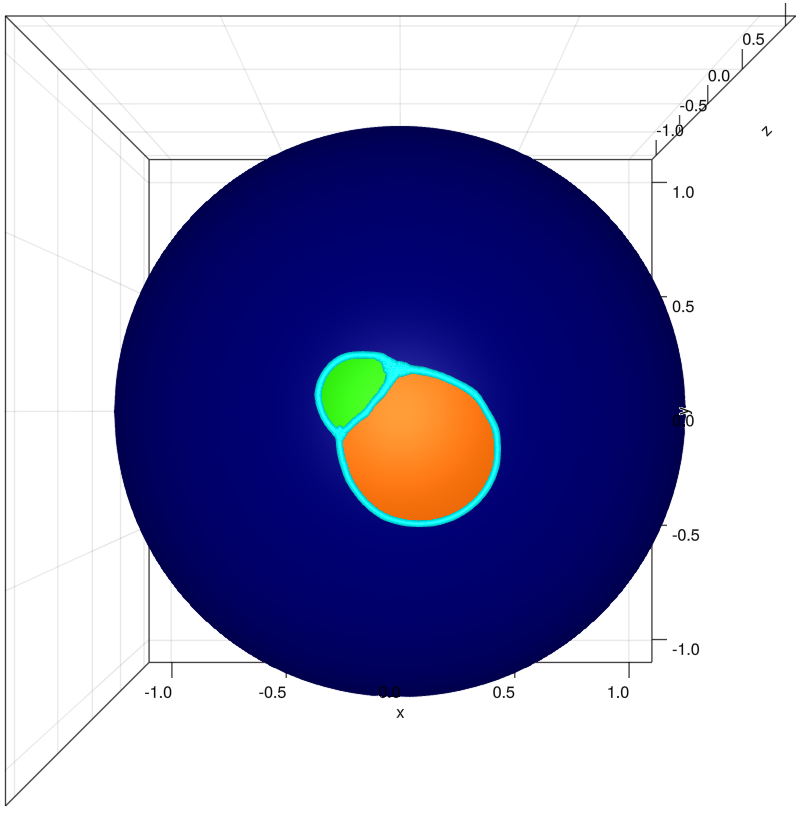}
                \subcaption{}

            \end{minipage} &
            \begin{minipage}[t]{0.25\columnwidth}
                \centering
                \includegraphics[height=0.15\vsize]{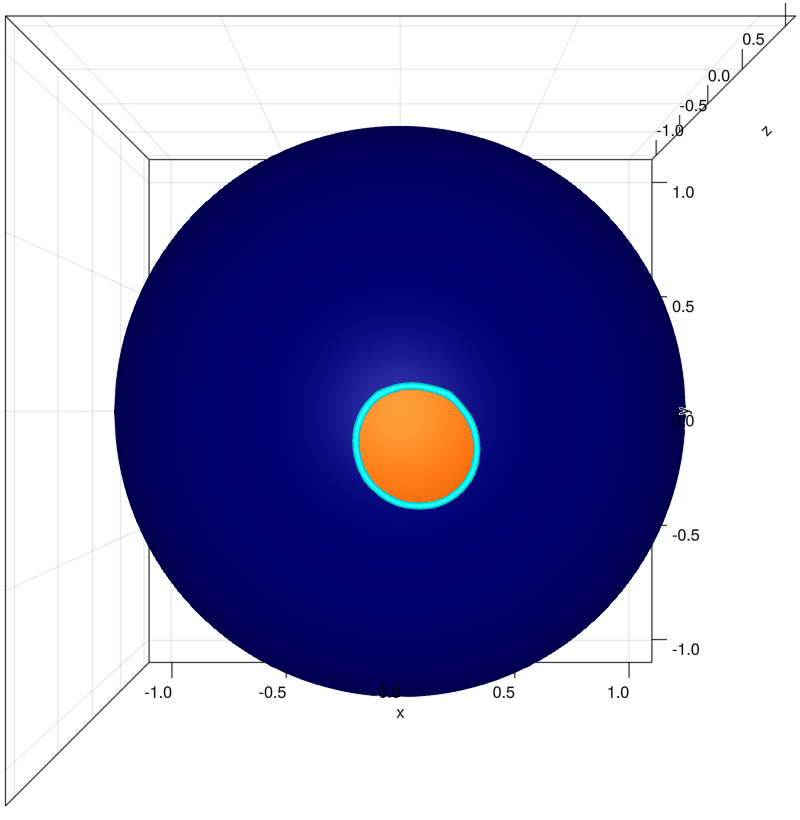}
                \subcaption{}

            \end{minipage}
            \\
            \begin{minipage}[t]{0.25\columnwidth}
                \centering
                \includegraphics[height=0.15\vsize]{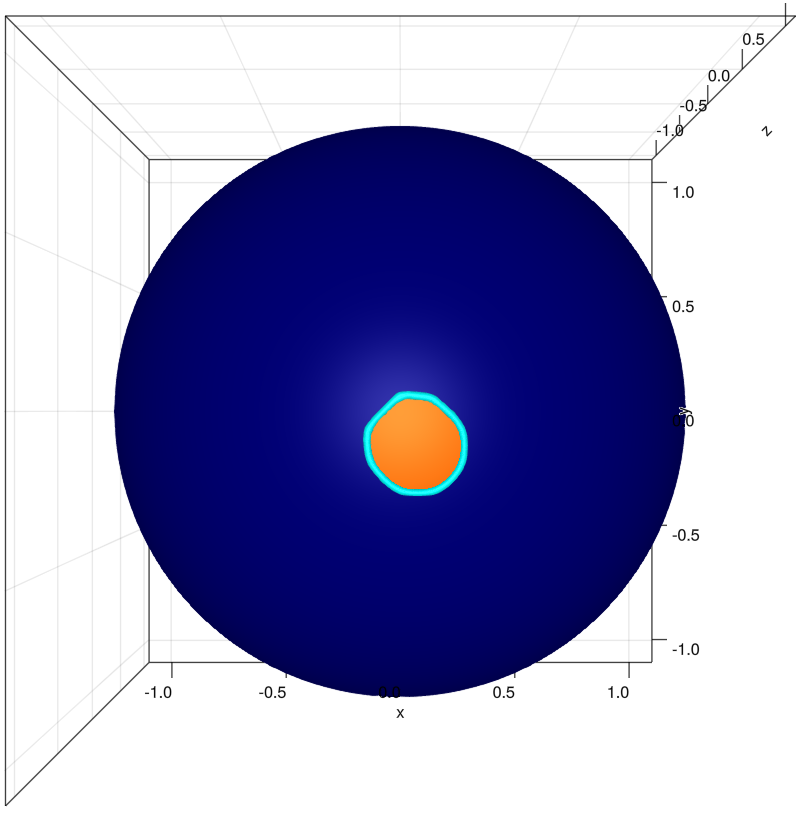}
                \subcaption{}

            \end{minipage} &
            \begin{minipage}[t]{0.25\columnwidth}
                \centering
                \includegraphics[height=0.15\vsize]{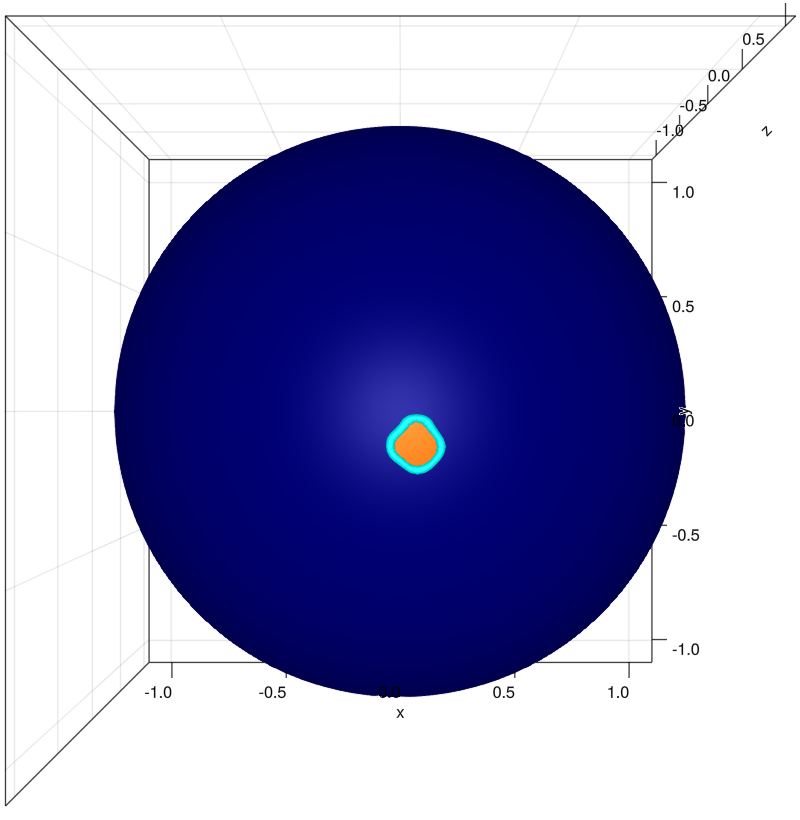}
                \subcaption{}

            \end{minipage}
        \end{tabular}
    }
    \caption{Result of four-phase MCF on the unit sphere (without area preservation). Arranged in alphabetical order with equal intervals between the initial time and the time the interface disappears.}
    \label{res:sp4_mcf}
\end{figure}

\begin{figure}[H]
    \centering
    \fbox{
        \begin{tabular}{cc}
            \begin{minipage}[t]{0.25\columnwidth}
                \centering
                \includegraphics[height=0.15\vsize]{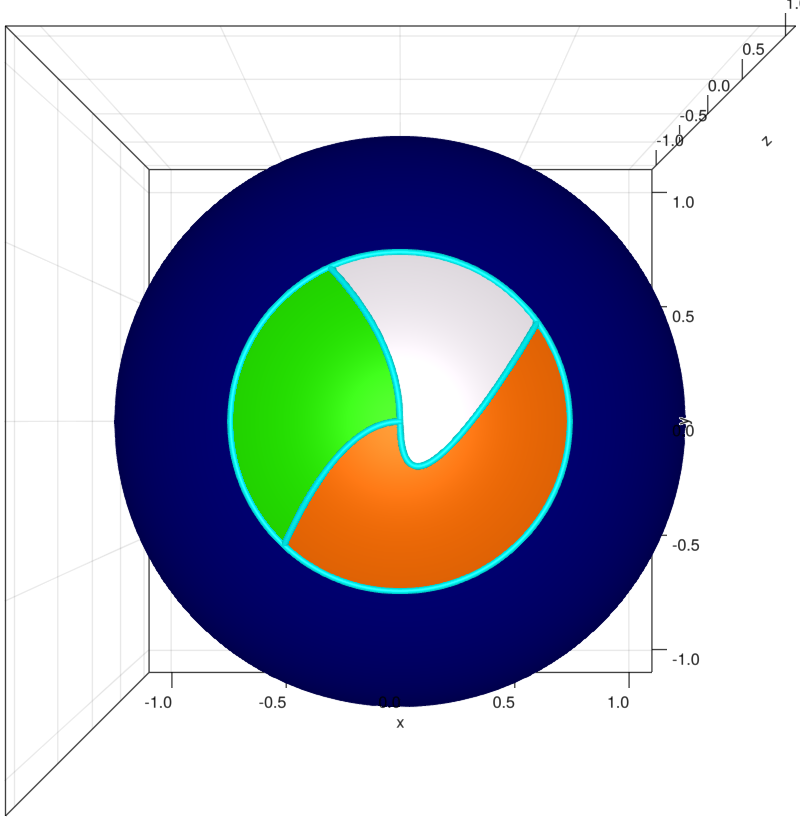}
                \subcaption{}

            \end{minipage} &
            \begin{minipage}[t]{0.25\columnwidth}
                \centering
                \includegraphics[height=0.15\vsize]{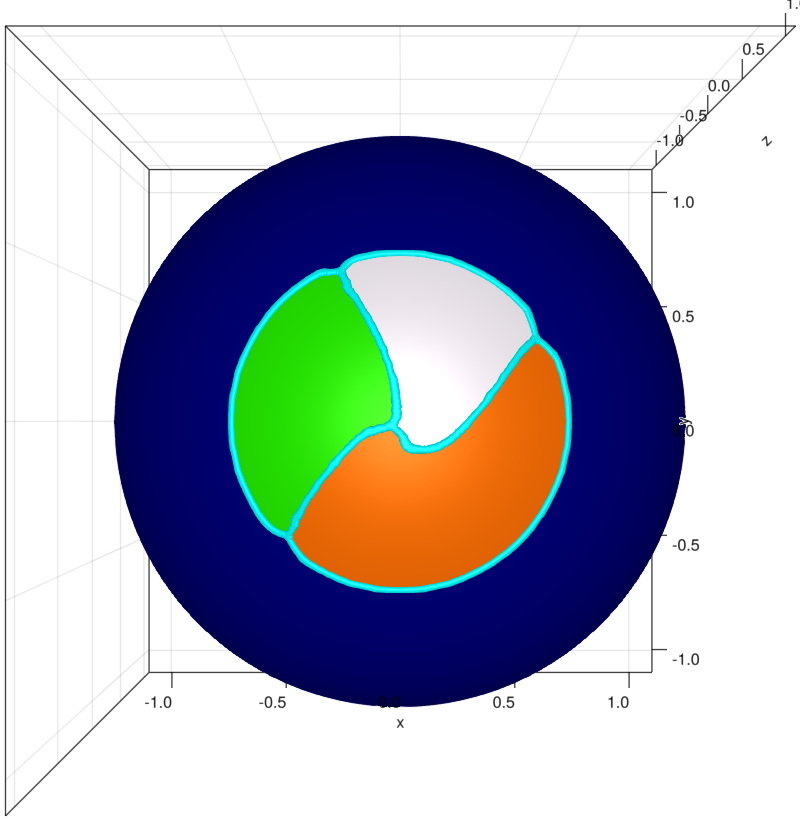}
                \subcaption{}
                \label{}
            \end{minipage}
            \\
            \begin{minipage}[t]{0.25\columnwidth}
                \centering
                \includegraphics[height=0.15\vsize]{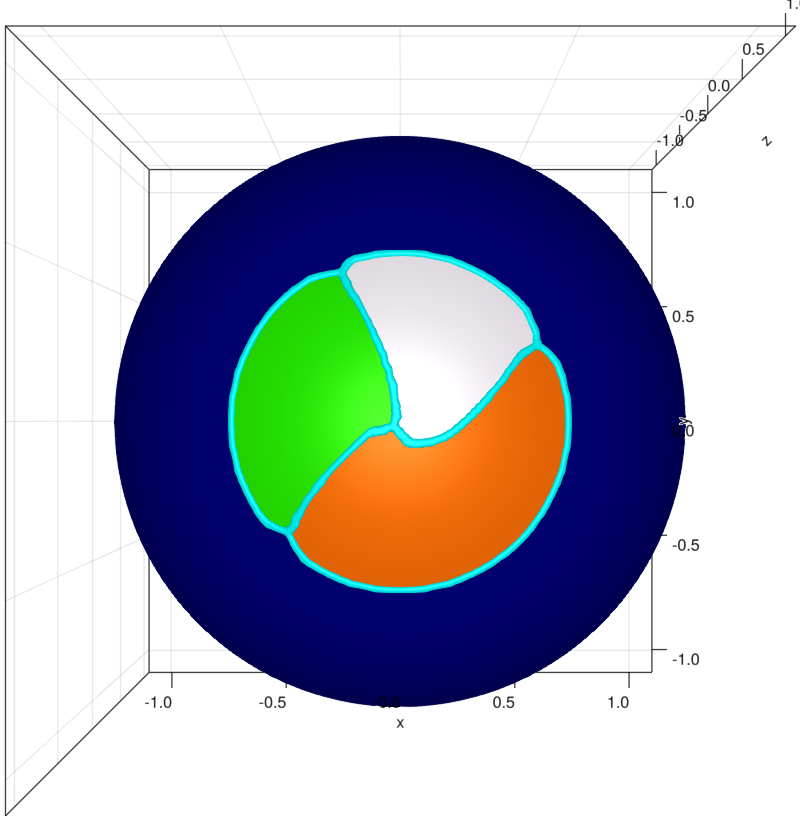}
                \subcaption{}

            \end{minipage} &
            \begin{minipage}[t]{0.25\columnwidth}
                \centering
                \includegraphics[height=0.15\vsize]{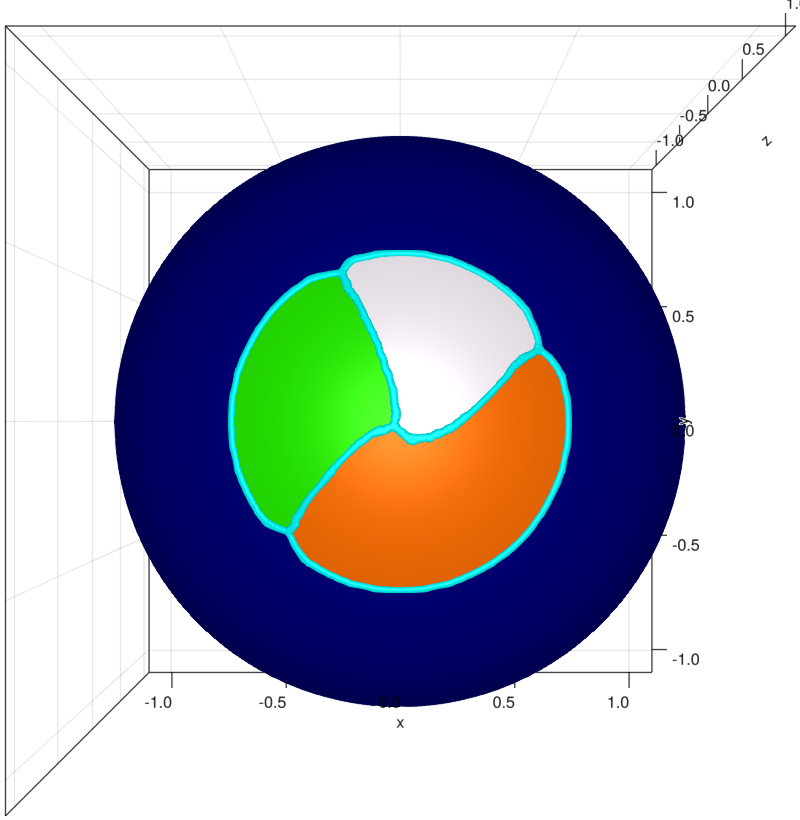}
                \subcaption{}

            \end{minipage}
            \\
            \begin{minipage}[t]{0.25\columnwidth}
                \centering
                \includegraphics[height=0.15\vsize]{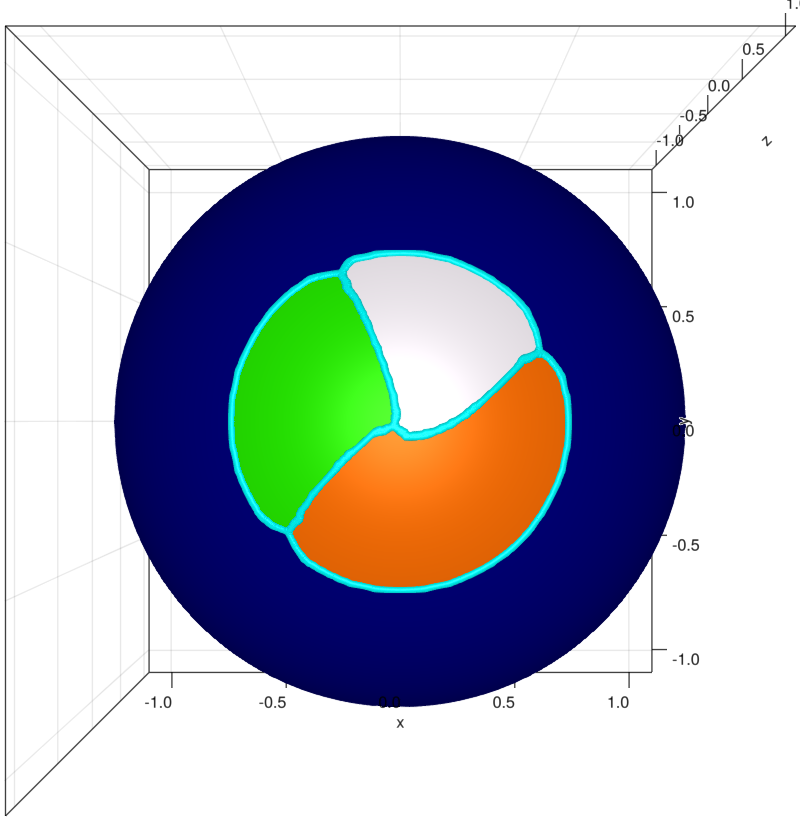}
                \subcaption{}

            \end{minipage} &
            \begin{minipage}[t]{0.25\columnwidth}
                \centering
                \includegraphics[height=0.15\vsize]{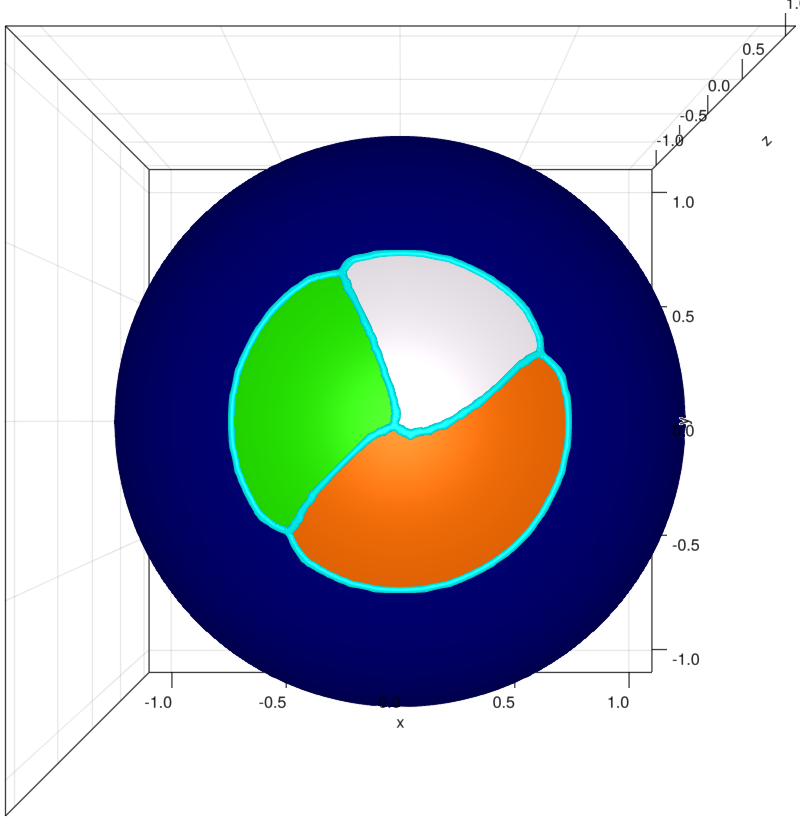}
                \subcaption{}

            \end{minipage}
            \\
            \begin{minipage}[t]{0.25\columnwidth}
                \centering
                \includegraphics[height=0.15\vsize]{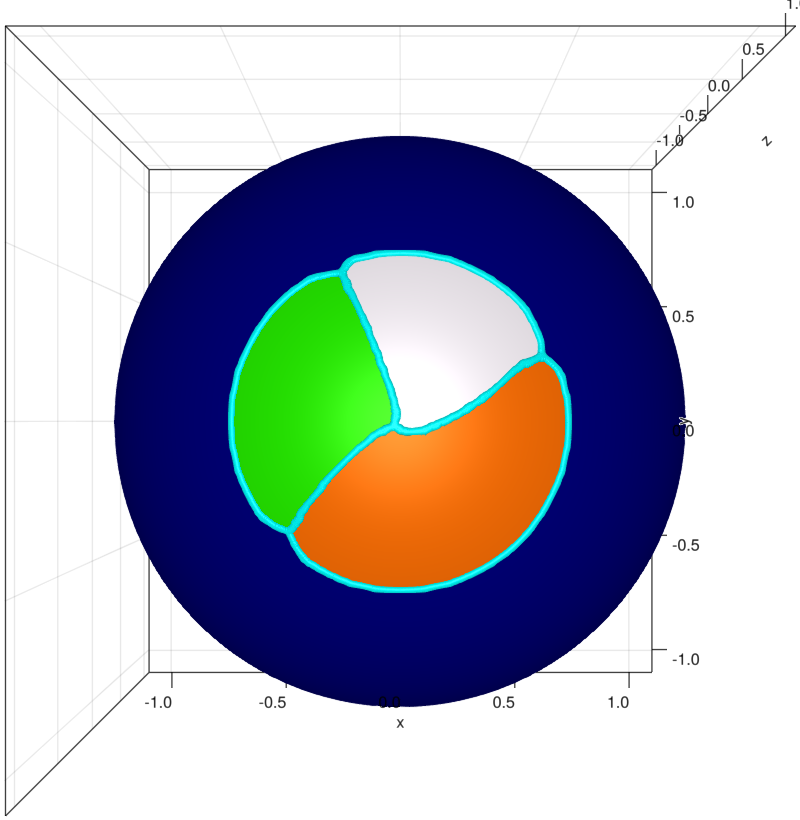}
                \subcaption{}

            \end{minipage} &
            \begin{minipage}[t]{0.25\columnwidth}
                \centering
                \includegraphics[height=0.15\vsize]{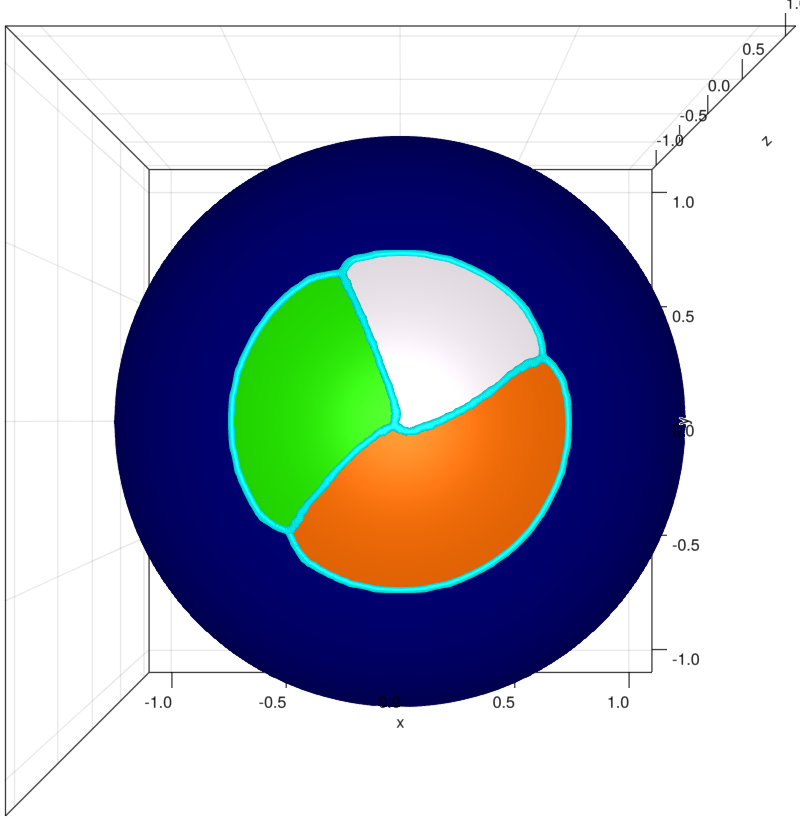}
                \subcaption{}

            \end{minipage}
            \\
            \begin{minipage}[t]{0.25\columnwidth}
                \centering
                \includegraphics[height=0.15\vsize]{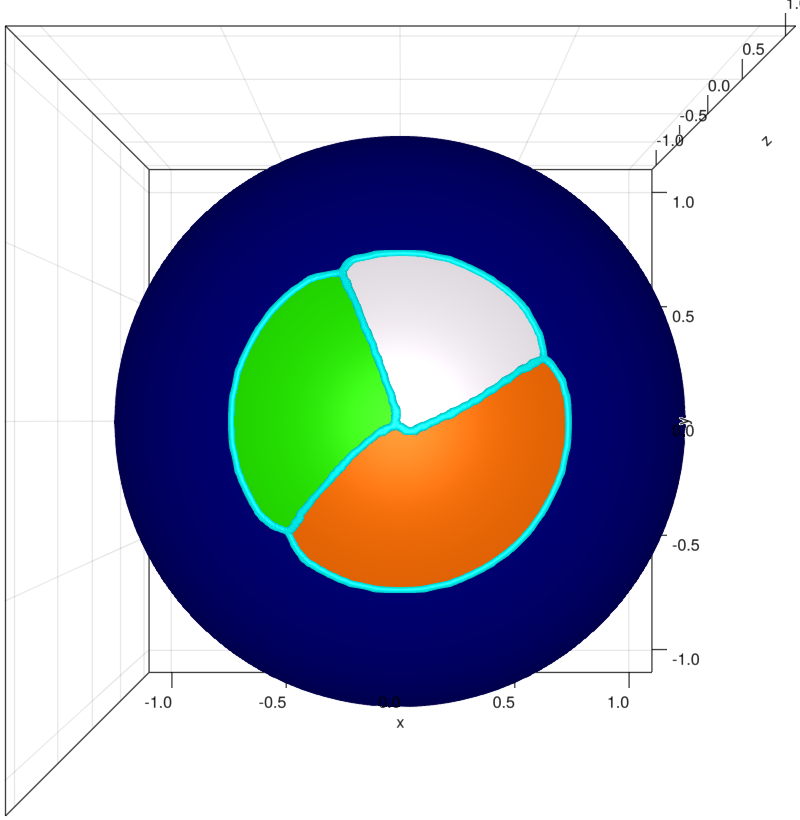}
                \subcaption{}

            \end{minipage} &
            \begin{minipage}[t]{0.25\columnwidth}
                \centering
                \includegraphics[height=0.15\vsize]{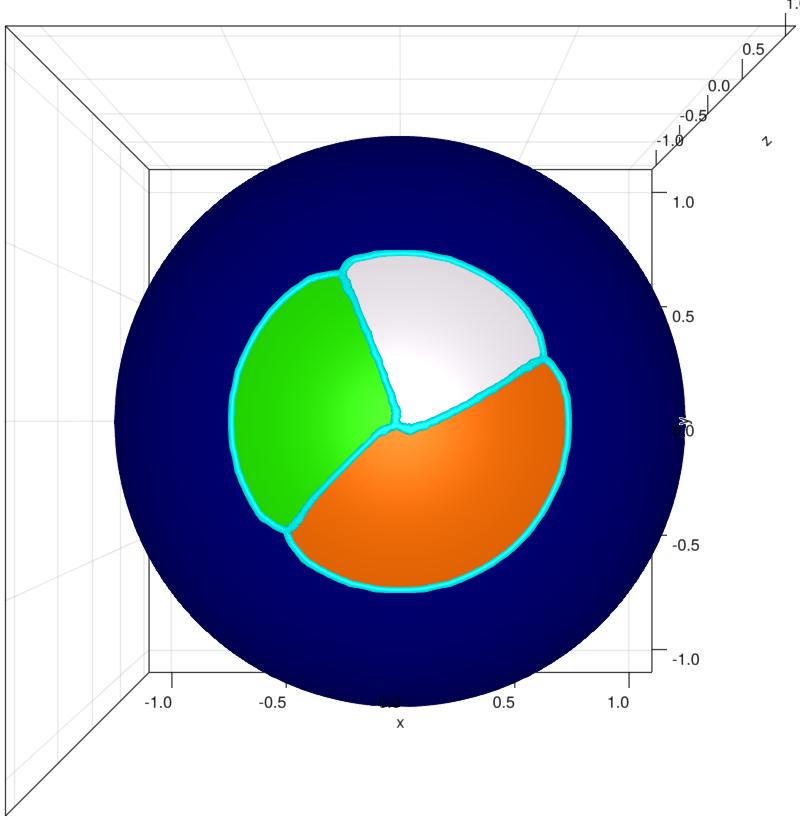}
                \subcaption{}

            \end{minipage}
        \end{tabular}
    }
    \caption{Result 1 of four-phase MCF on the unit sphere (with area preservation), the time intervals from the initial time to the time the interface reaches a nearly stationary state are equally spaced. }
    \label{res:sp4_mcf_cons}
\end{figure}

\begin{figure}[H]
    \centering
    \fbox{
        \begin{tabular}{cc}
            \begin{minipage}[t]{0.25\columnwidth}
                \centering
                \includegraphics[height=0.15\vsize]{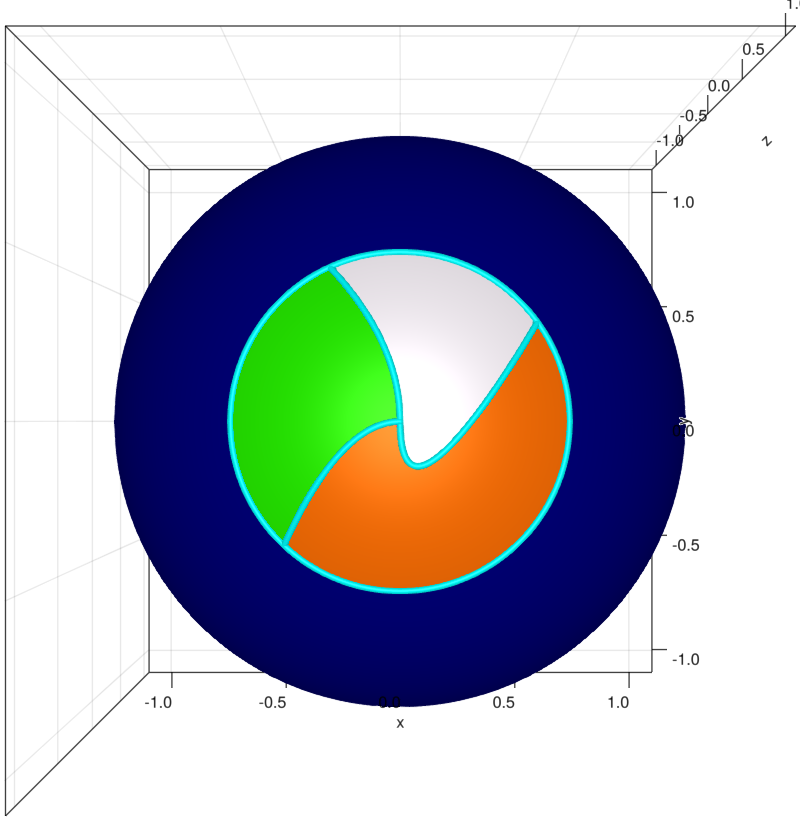}
                \subcaption{}
            \end{minipage} &
            \begin{minipage}[t]{0.25\columnwidth}
                \centering
                \includegraphics[height=0.15\vsize]{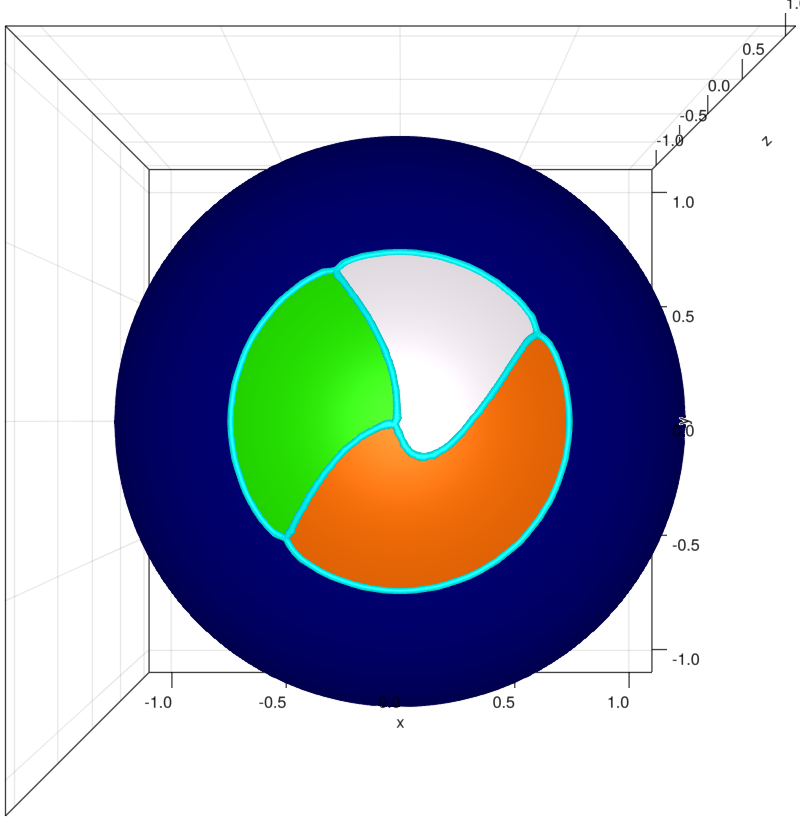}
                \subcaption{}
                \label{}
            \end{minipage}
            \\
            \begin{minipage}[t]{0.25\columnwidth}
                \centering
                \includegraphics[height=0.15\vsize]{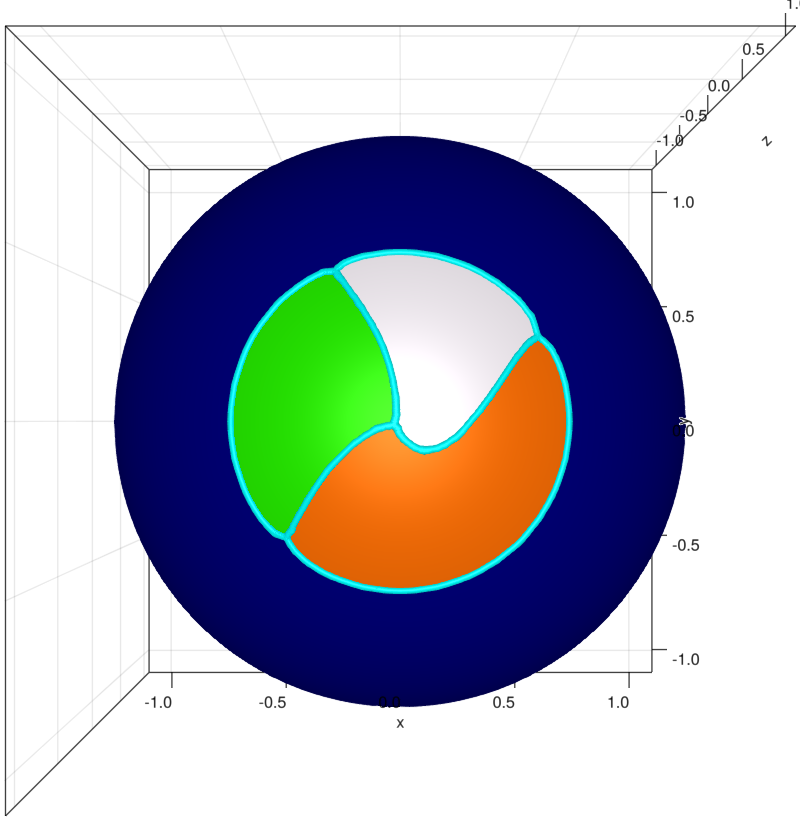}
                \subcaption{}
            \end{minipage} &
            \begin{minipage}[t]{0.25\columnwidth}
                \centering
                \includegraphics[height=0.15\vsize]{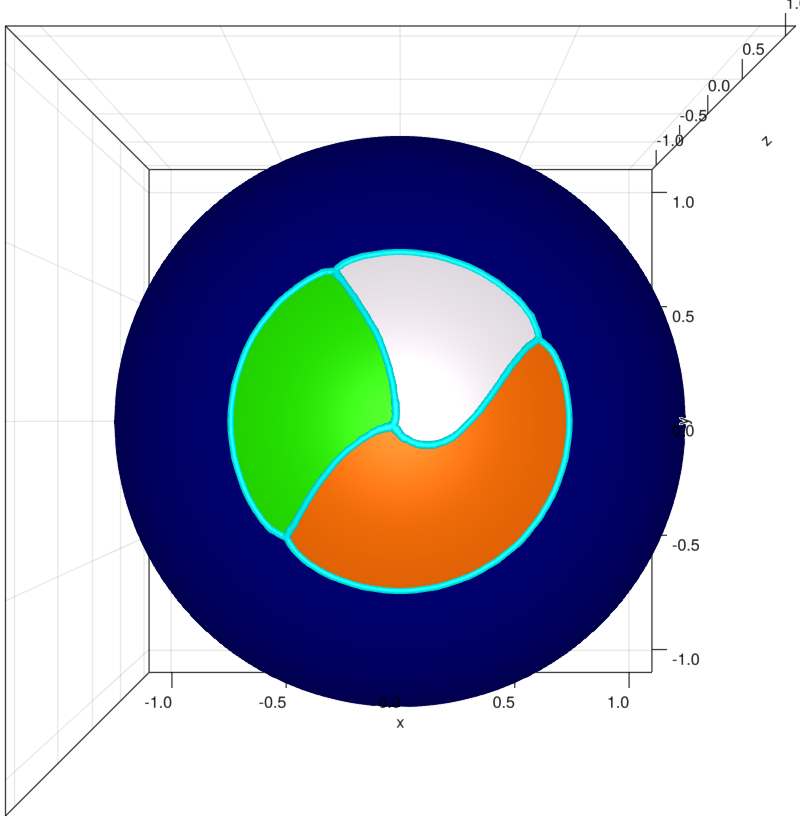}
                \subcaption{}
            \end{minipage}
            \\
            \begin{minipage}[t]{0.25\columnwidth}
                \centering
                \includegraphics[height=0.15\vsize]{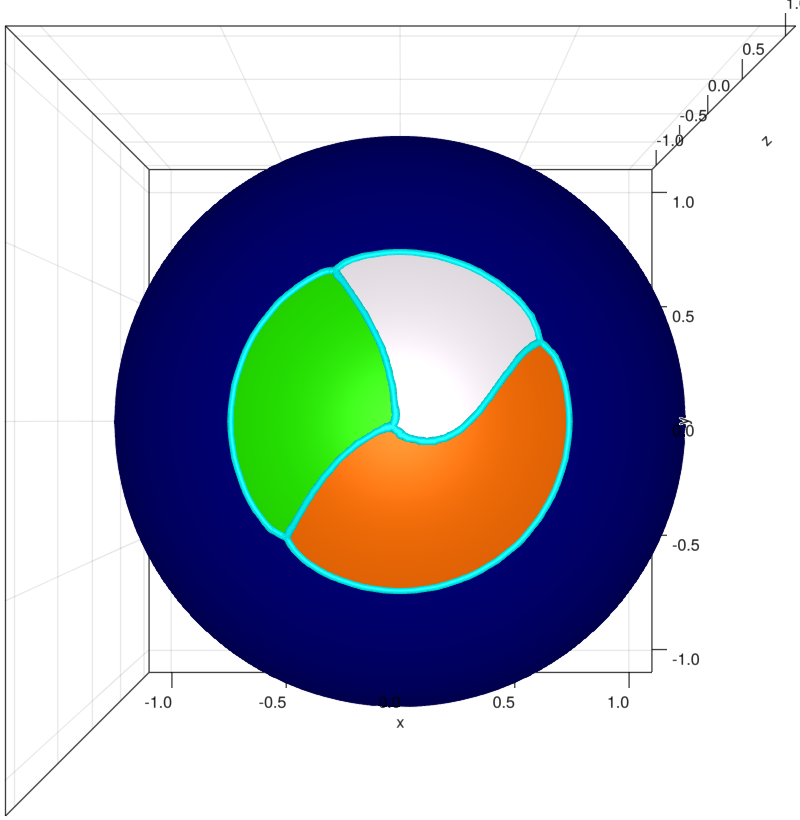}
                \subcaption{}
            \end{minipage} &
            \begin{minipage}[t]{0.25\columnwidth}
                \centering
                \includegraphics[height=0.15\vsize]{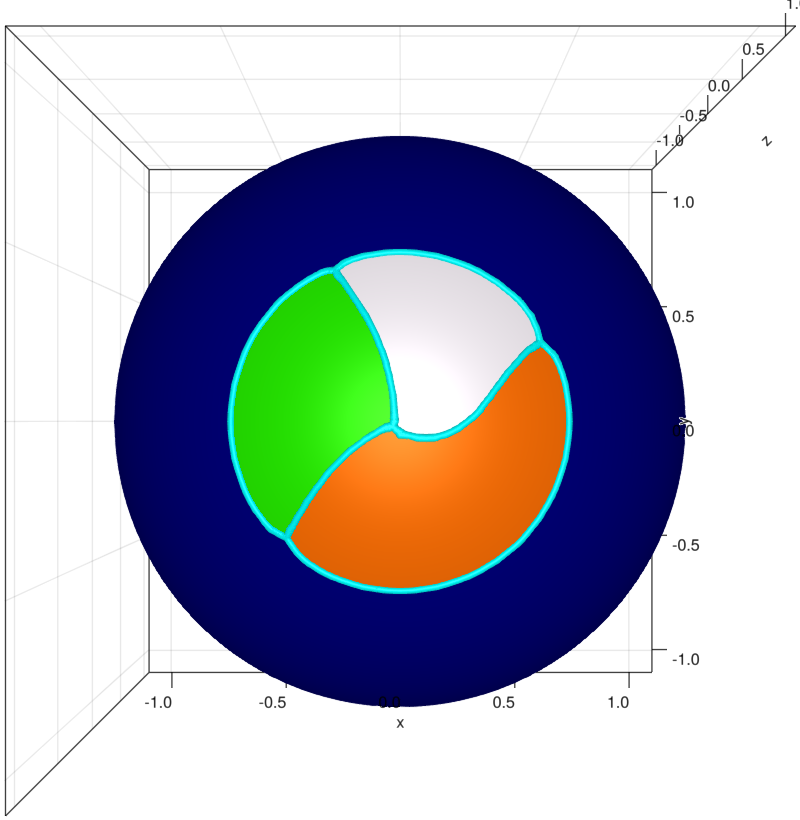}
                \subcaption{}
            \end{minipage}
            \\
            \begin{minipage}[t]{0.25\columnwidth}
                \centering
                \includegraphics[height=0.15\vsize]{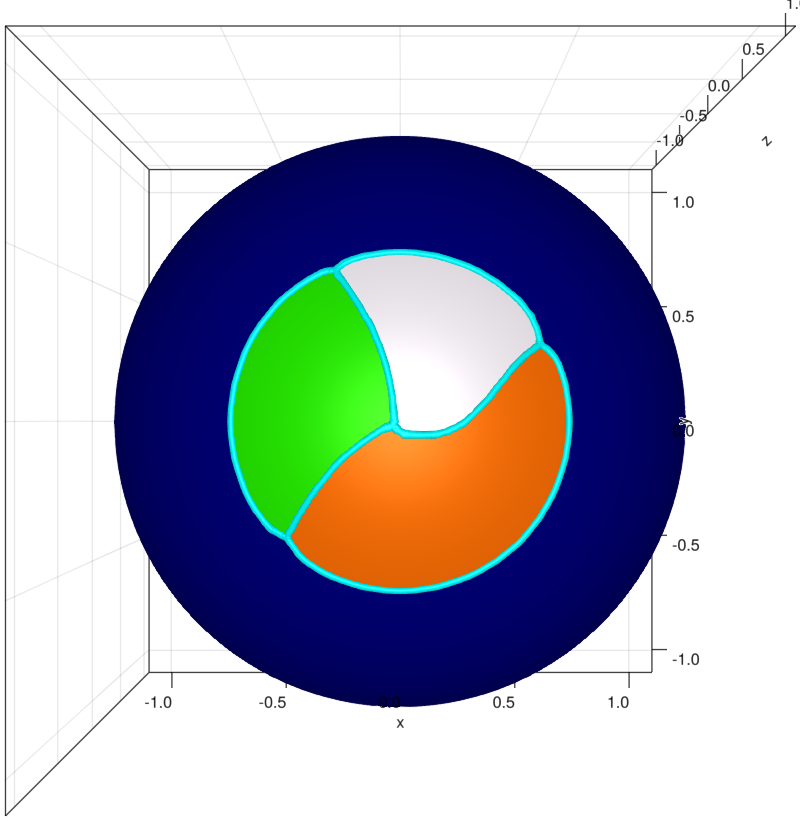}
                \subcaption{}
            \end{minipage} &
            \begin{minipage}[t]{0.25\columnwidth}
                \centering
                \includegraphics[height=0.15\vsize]{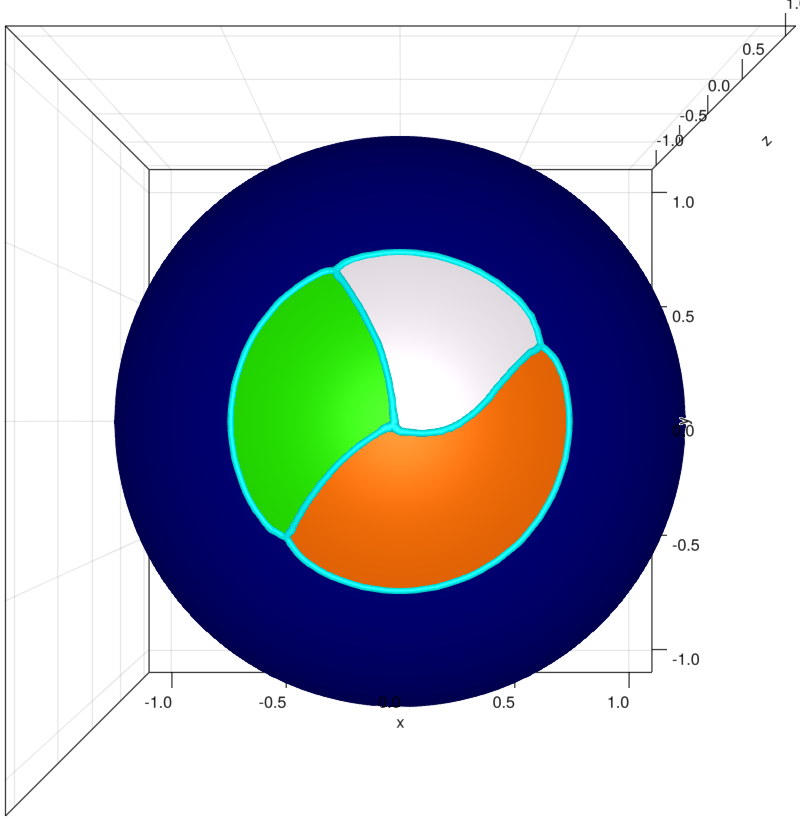}
                \subcaption{}
            \end{minipage}
            \\
            \begin{minipage}[t]{0.25\columnwidth}
                \centering
                \includegraphics[height=0.15\vsize]{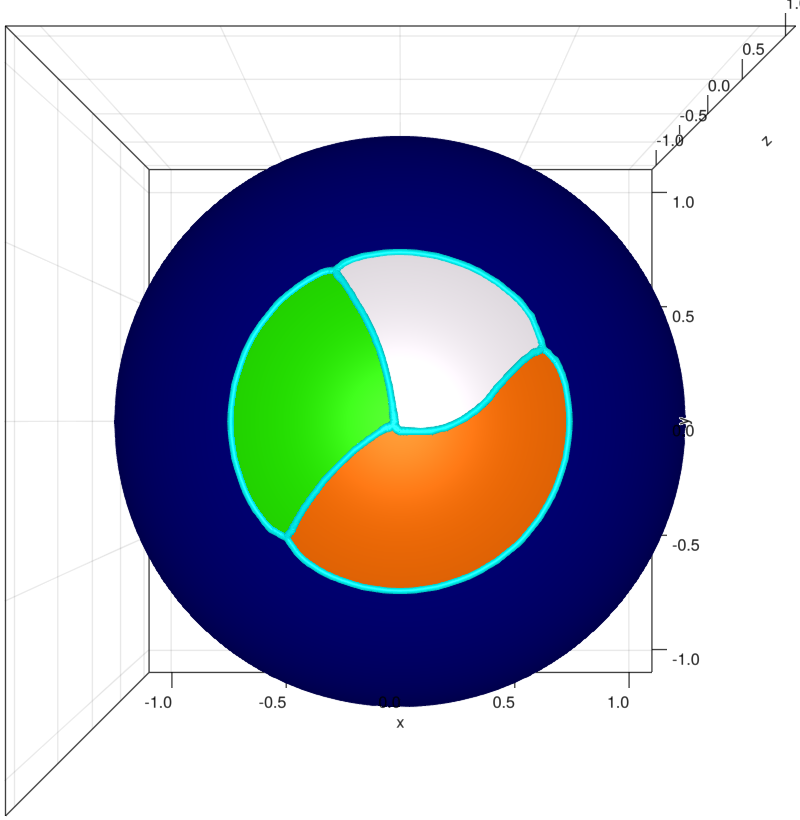}
                \subcaption{}

            \end{minipage} &
            \begin{minipage}[t]{0.25\columnwidth}
                \centering
                \includegraphics[height=0.15\vsize]{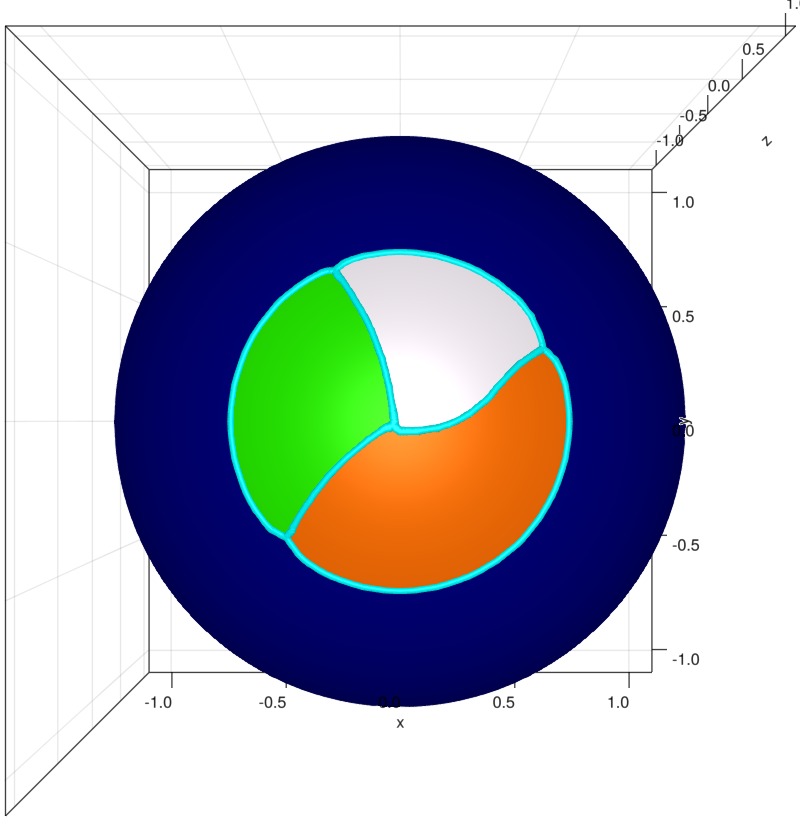}
                \subcaption{}

            \end{minipage}
        \end{tabular}
    }
    \caption{Result 2 of four-phase MCF on the unit sphere (with area preservation). Arranged in alphabetical order. The time intervals from the initial time to the time when the interface reaches a nearly stationary state are equally spaced. }
    \label{res:sp4_mcf_cons_2}
\end{figure}

Discussions of the corresponding computations are presented in section \ref{結果に対する考察}.
\begin{description}
    \item[Result of two-phase mean curvature flow on the unit sphere (without area preservation)]~\\

    We use the two-phase initial condition shown in Figure \ref{ini_sp_phase2} and the algorithm described in section \ref{algo-smboc} for the numerical calculation. 
    Parameters were set as follows:
    $$
        \alpha=1.0,\quad \Delta x=0.05, \quad h= \Delta x^2/6,\quad\tau=15h,\quad \rho=0
    $$
    The numerical result is shown in Figure \ref{res:sp2_mcf}.
    \item[Result of two-phase mean curvature on the unit sphere (with area preservation)]~\\

    We used the two-phase initial condition shown in Figure \ref{ini_sp_phase2} and the algorithm described in section \ref{algo-smboc} for the numerical calculation. 
    Parameters were set as follows:
    $$
        \alpha=0.05,\quad\Delta x=0.05, \quad h= \Delta x^2/6,\quad\tau=100h,\quad \rho=10^3
    $$
    The numerical result is shown in Figure \ref{res:sp2_mcf_cons}.
    \item[Result of four-phase mean curvature flow on the unit sphere (without area preservation)]~\\

    We used the four-phase initial condition shown in Figure \ref{ini_sp_phase4} and the algorithm described in section \ref{algo-smboc} for the numerical calculation. 
    Parameters were set as follows:
    $$
        \alpha=1.0,\quad\Delta x=0.05, \quad h= \Delta x^2/6,\quad\tau=15h,\quad \rho=0
    $$
    The numerical result is shown in Figure \ref{res:sp4_mcf}.
    \item[Result 1 of four-phase mean curvature flow on the unit sphere (with area preservation)]~\\

    We used the four-phase initial condition shown in Figure \ref{ini_sp_phase4} and the algorithm described in Section \ref{algo-smboc} for the numerical calculation. 
    Parameters were set as follows:
    $$
        \alpha=1.0,\quad\Delta x=0.015, \quad h= \Delta x^2,\quad\tau=15h,\quad \rho=10^5
    $$
    The numerical result is shown in Figure \ref{res:sp4_mcf_cons}.
    \item[Result 2 of four-phase mean curvature flow on the unit sphere (with area preservation)]~\\
    We used the four-phase initial condition shown in Figure \ref{ini_sp_phase4} and the algorithm described in Section \ref{algo-smboc} for the numerical calculation. 
    Parameter were set as follows:
    $$
        \alpha=0.1,\quad\Delta x=0.01, \quad h= \Delta x^2,\quad\tau=100h,\quad \rho=4\times10^4
    $$
    The numerical result is shown in Figure \ref{res:sp4_mcf_cons_2}.
\end{description}
\subsection{Computational details regarding surface hyperbolic mean curvature flow}\label{results_HMBO}
\begin{figure}[H]
    \centering
    \fbox{
        \begin{tabular}{cc}
            \begin{minipage}[t]{0.25\columnwidth}
                \centering
                \includegraphics[height=0.15\vsize]{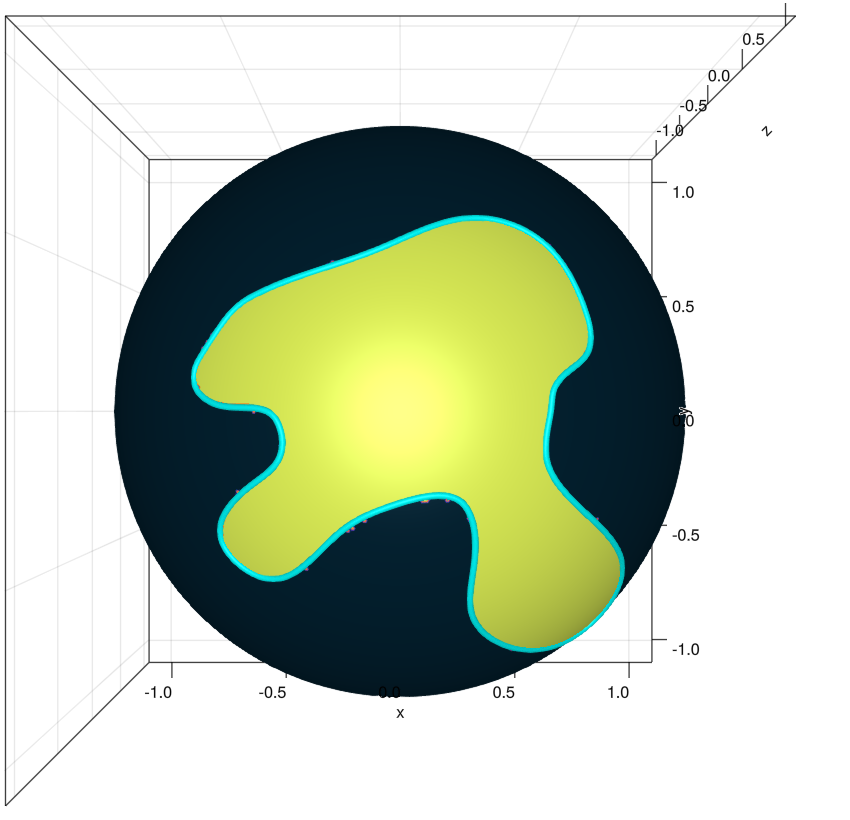}
                \subcaption{}

            \end{minipage} &
            \begin{minipage}[t]{0.25\columnwidth}
                \centering
                \includegraphics[height=0.15\vsize]{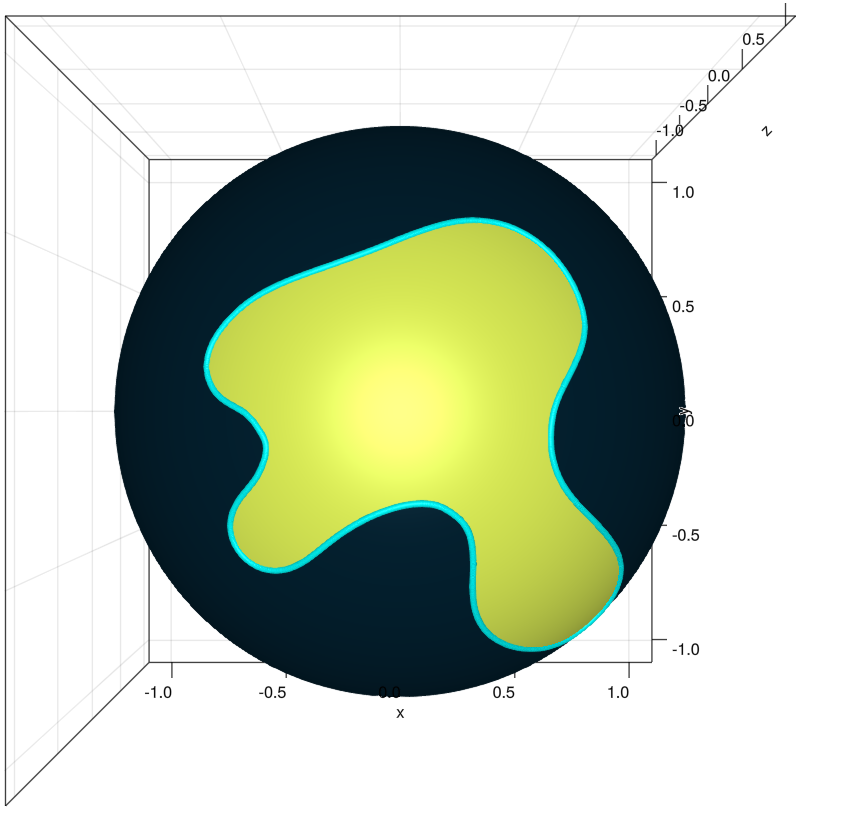}
                \subcaption{}
                \label{}
            \end{minipage}
            \\
            \begin{minipage}[t]{0.25\columnwidth}
                \centering
                \includegraphics[height=0.15\vsize]{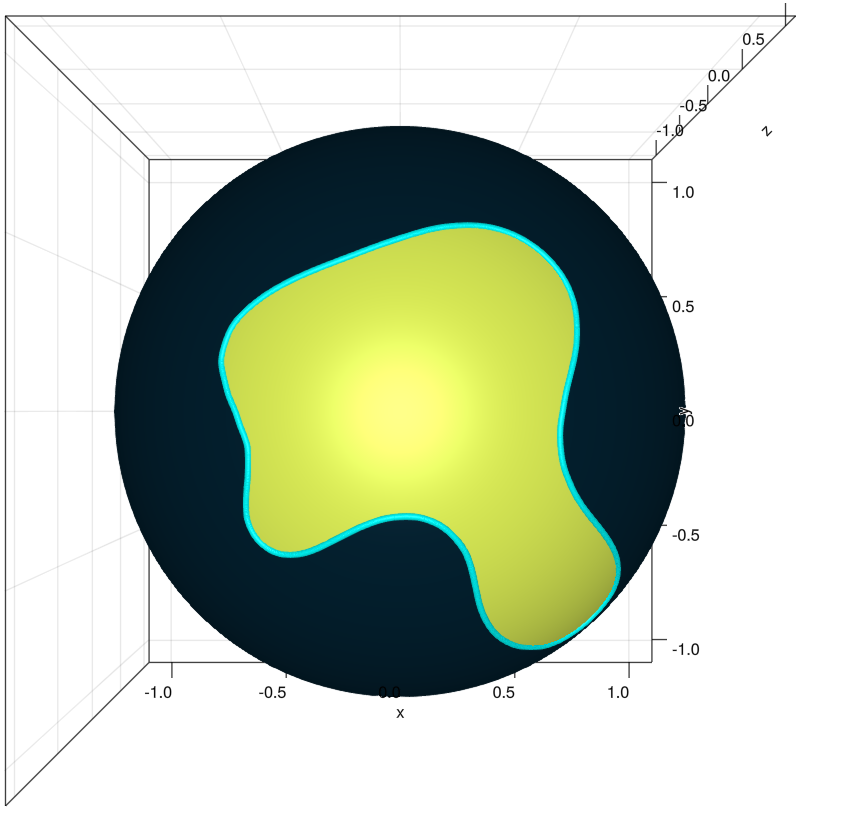}
                \subcaption{}

            \end{minipage} &
            \begin{minipage}[t]{0.25\columnwidth}
                \centering
                \includegraphics[height=0.15\vsize]{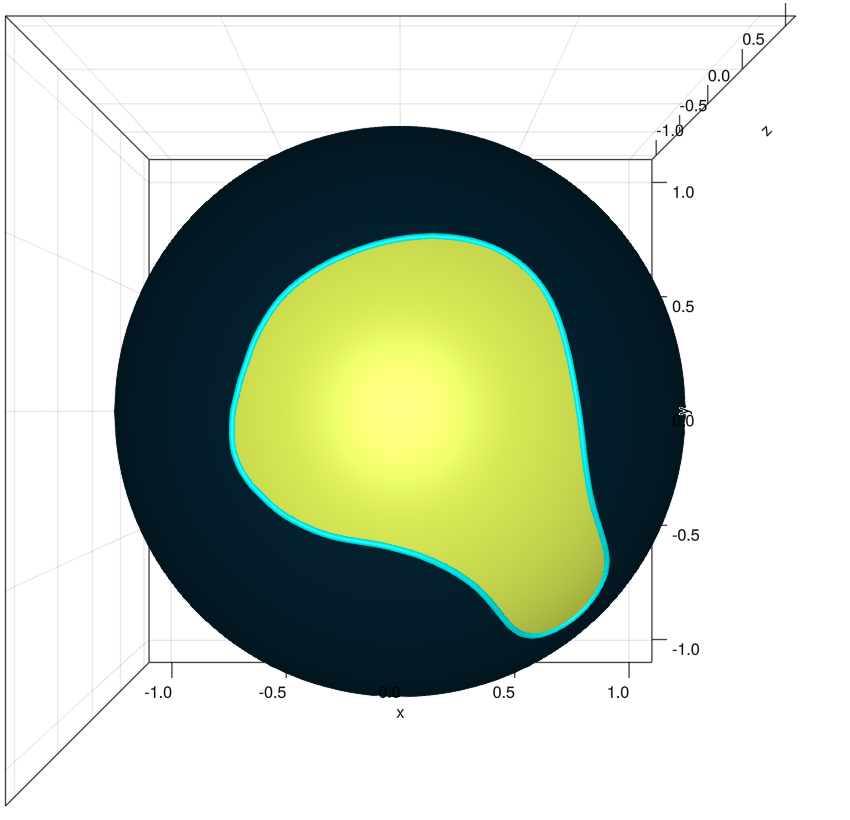}
                \subcaption{}

            \end{minipage}
            \\
            \begin{minipage}[t]{0.25\columnwidth}
                \centering
                \includegraphics[height=0.15\vsize]{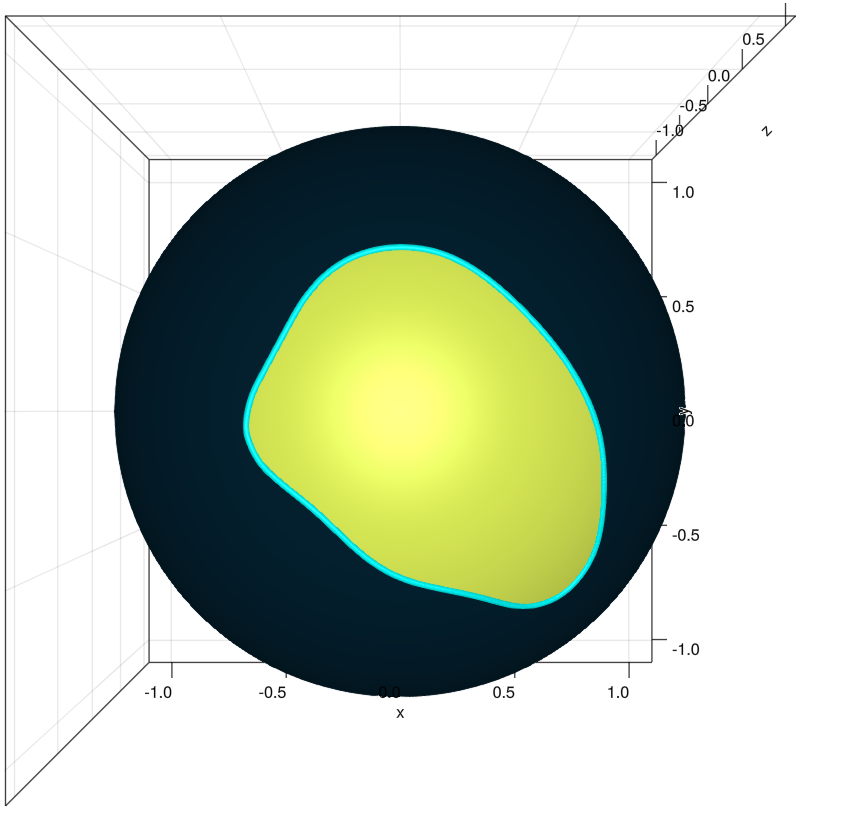}
                \subcaption{}

            \end{minipage} &
            \begin{minipage}[t]{0.25\columnwidth}
                \centering
                \includegraphics[height=0.15\vsize]{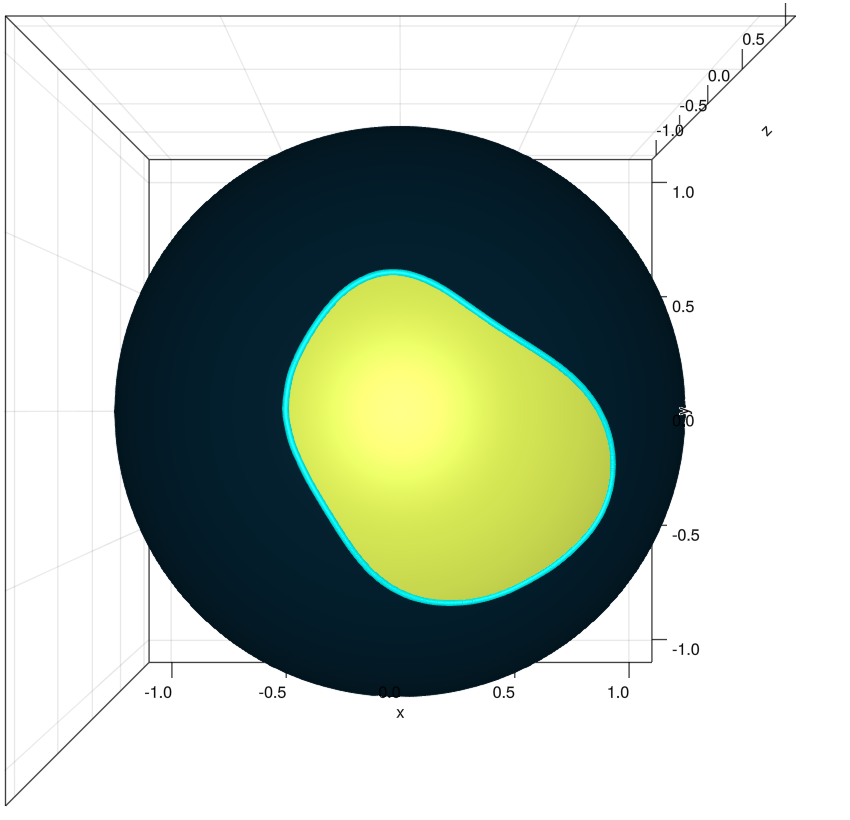}
                \subcaption{}

            \end{minipage}
            \\
            \begin{minipage}[t]{0.25\columnwidth}
                \centering
                \includegraphics[height=0.15\vsize]{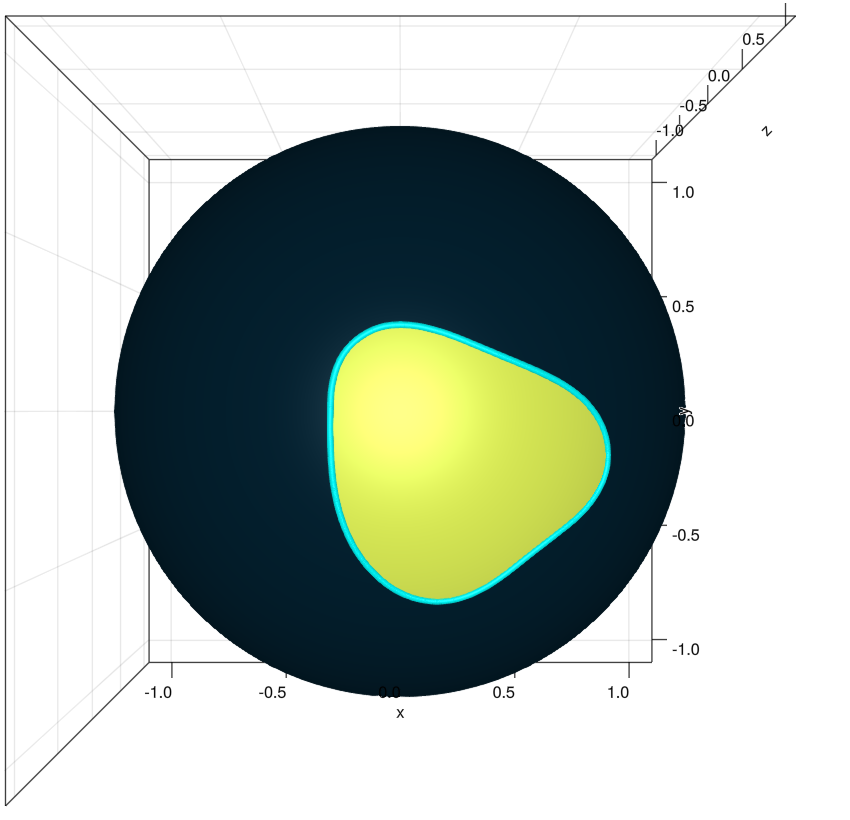}
                \subcaption{}

            \end{minipage} &
            \begin{minipage}[t]{0.25\columnwidth}
                \centering
                \includegraphics[height=0.15\vsize]{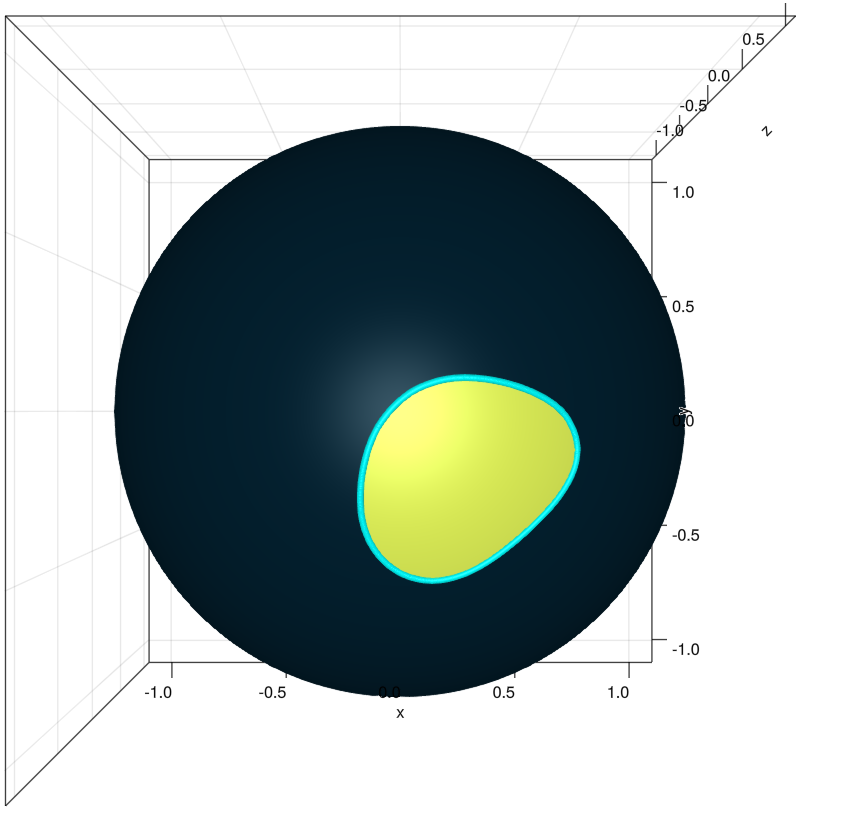}
                \subcaption{}

            \end{minipage}
            \\
            \begin{minipage}[t]{0.25\columnwidth}
                \centering
                \includegraphics[height=0.15\vsize]{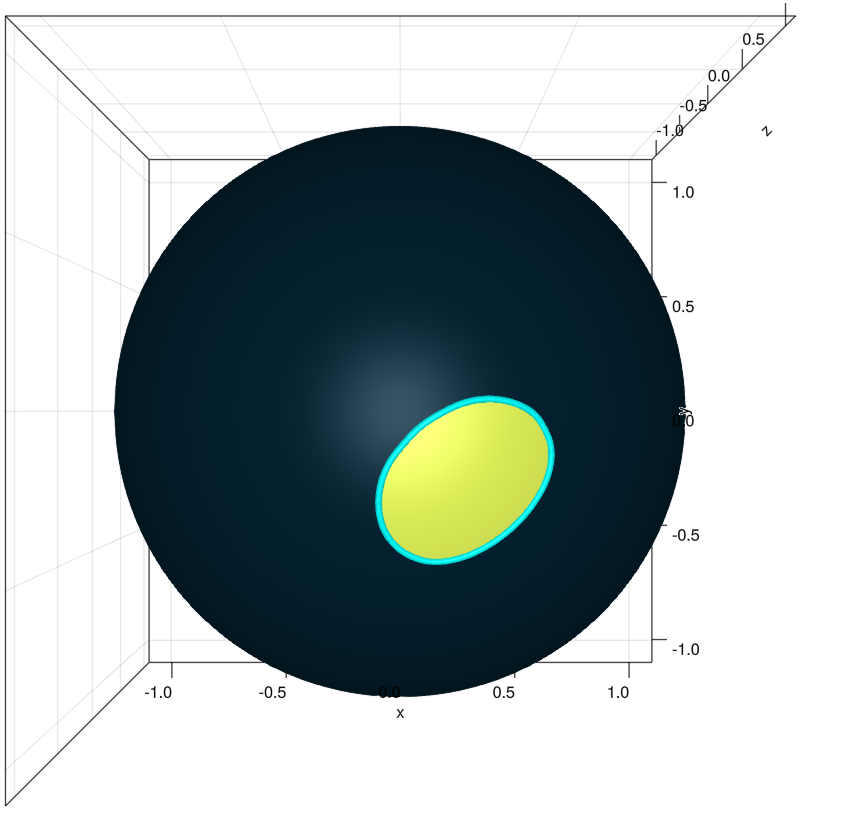}
                \subcaption{}

            \end{minipage} &
            \begin{minipage}[t]{0.25\columnwidth}
                \centering
                \includegraphics[height=0.15\vsize]{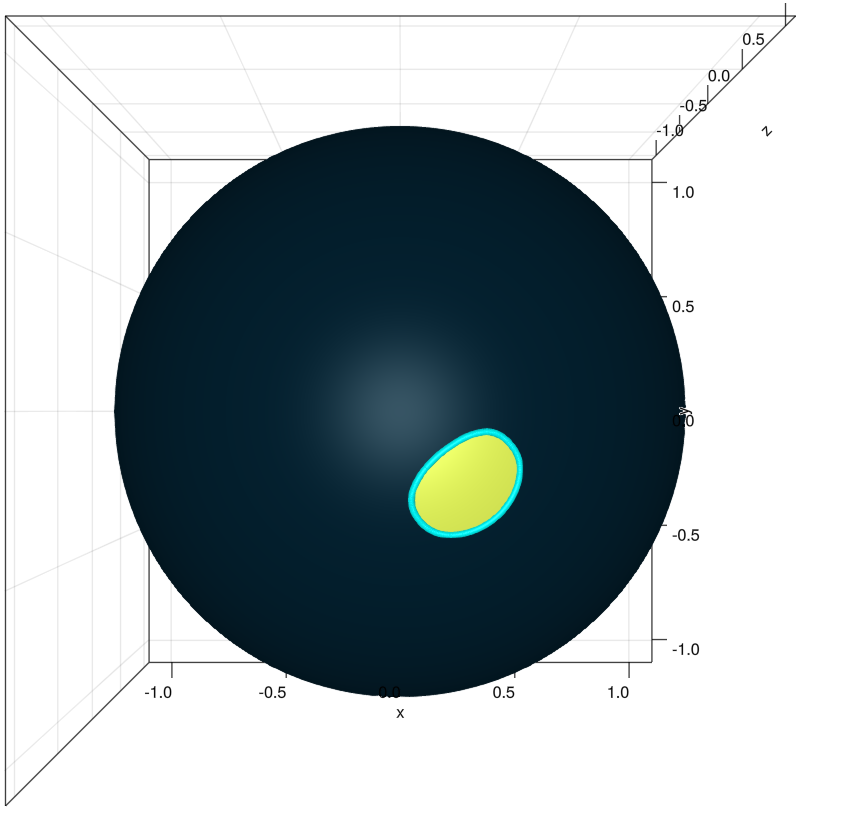}
                \subcaption{}

            \end{minipage}
        \end{tabular}
    }
    \caption{Result of two-phase hyperbolic MCF on the unit sphere (without area preservation). Arranged in alphabetical order with equal intervals from the initial time to the time the interface disappears.}
    \label{res:sp2_hmcf}
\end{figure}

\begin{figure}[H]
    \centering
    \fbox{
        \begin{tabular}{cc}
            \begin{minipage}[t]{0.25\columnwidth}
                \centering
                \includegraphics[height=0.15\vsize]{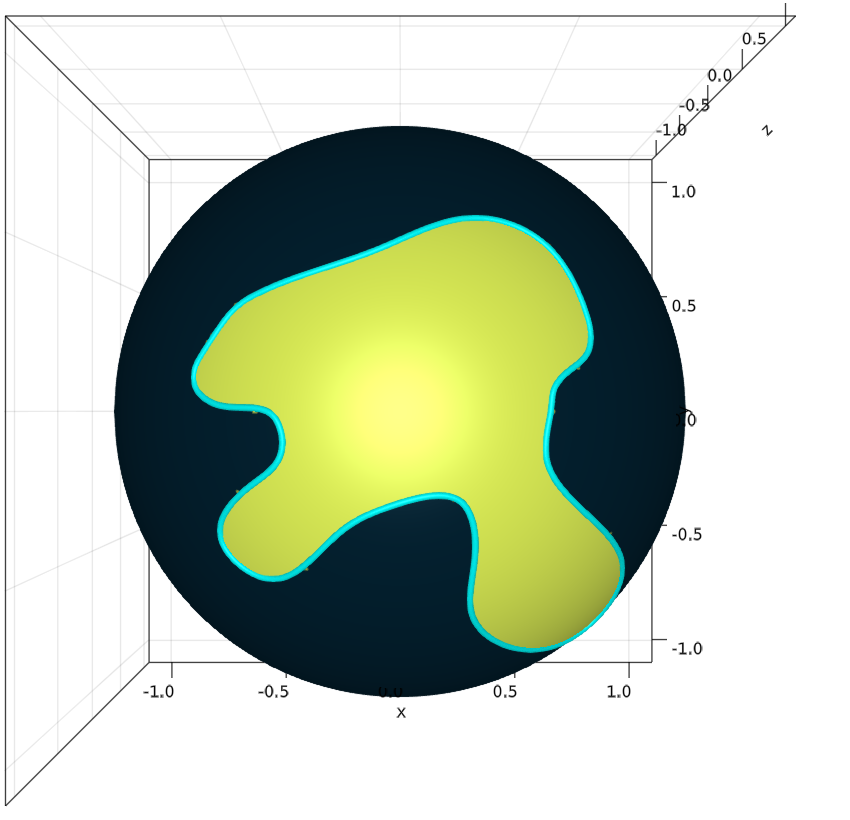}
                \subcaption{}

            \end{minipage} &
            \begin{minipage}[t]{0.25\columnwidth}
                \centering
                \includegraphics[height=0.15\vsize]{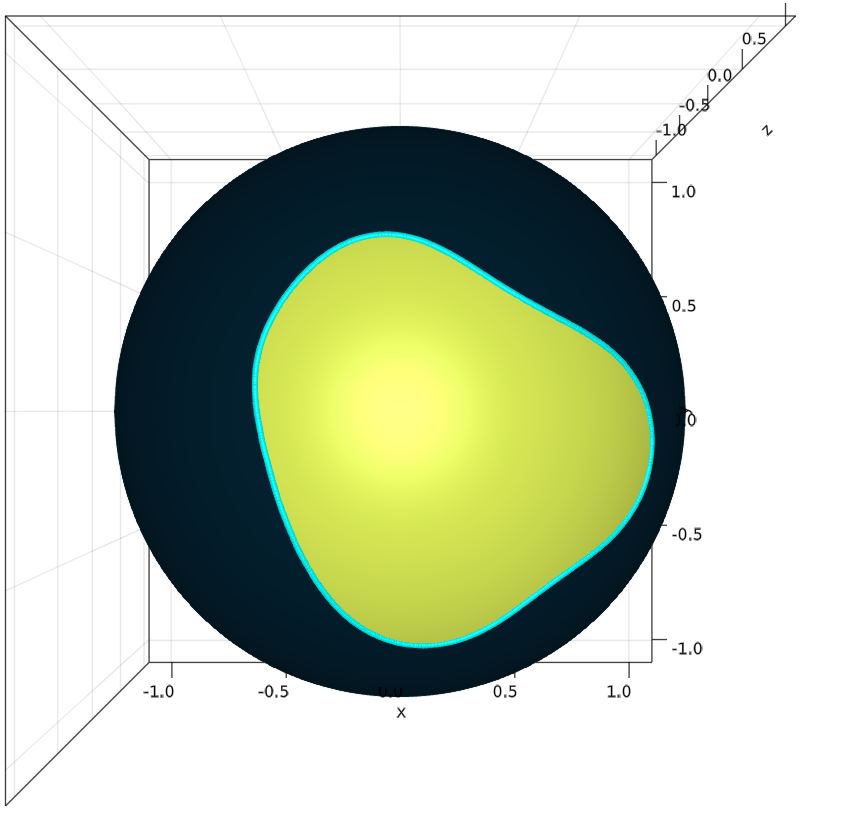}
                \subcaption{}
                \label{}
            \end{minipage}
            \\
            \begin{minipage}[t]{0.25\columnwidth}
                \centering
                \includegraphics[height=0.15\vsize]{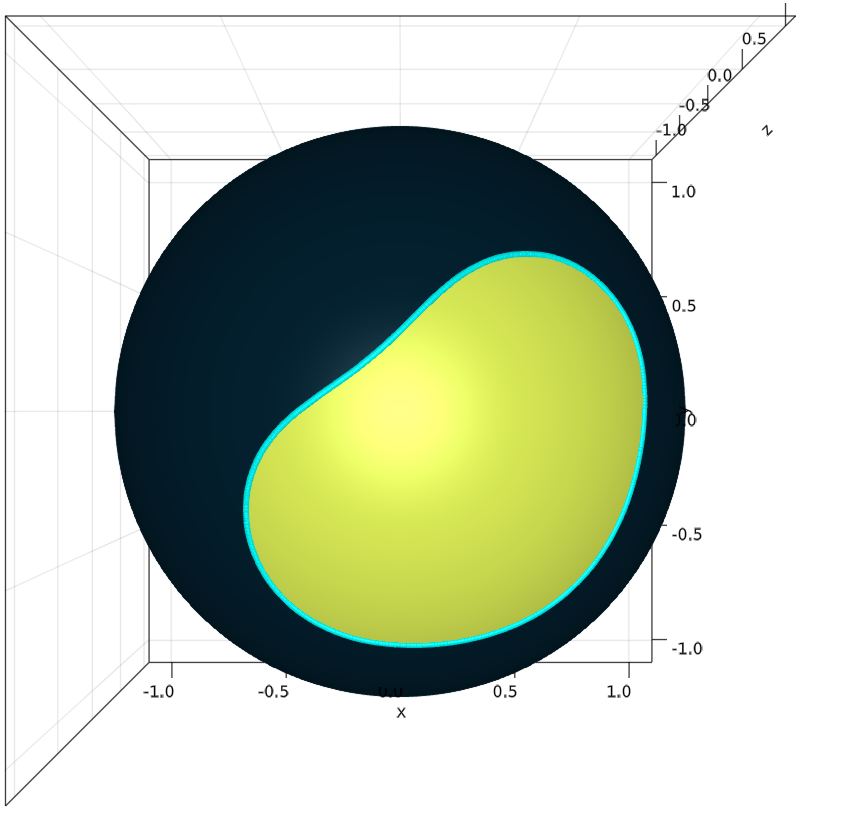}
                \subcaption{}

            \end{minipage} &
            \begin{minipage}[t]{0.25\columnwidth}
                \centering
                \includegraphics[height=0.15\vsize]{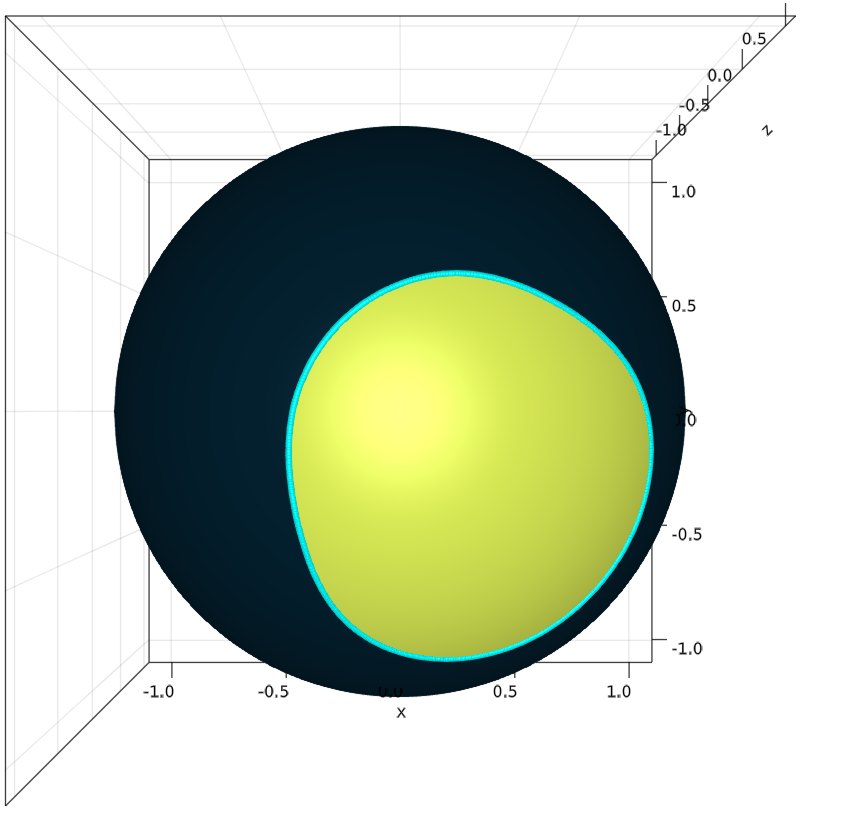}
                \subcaption{}

            \end{minipage}
            \\
            \begin{minipage}[t]{0.25\columnwidth}
                \centering
                \includegraphics[height=0.15\vsize]{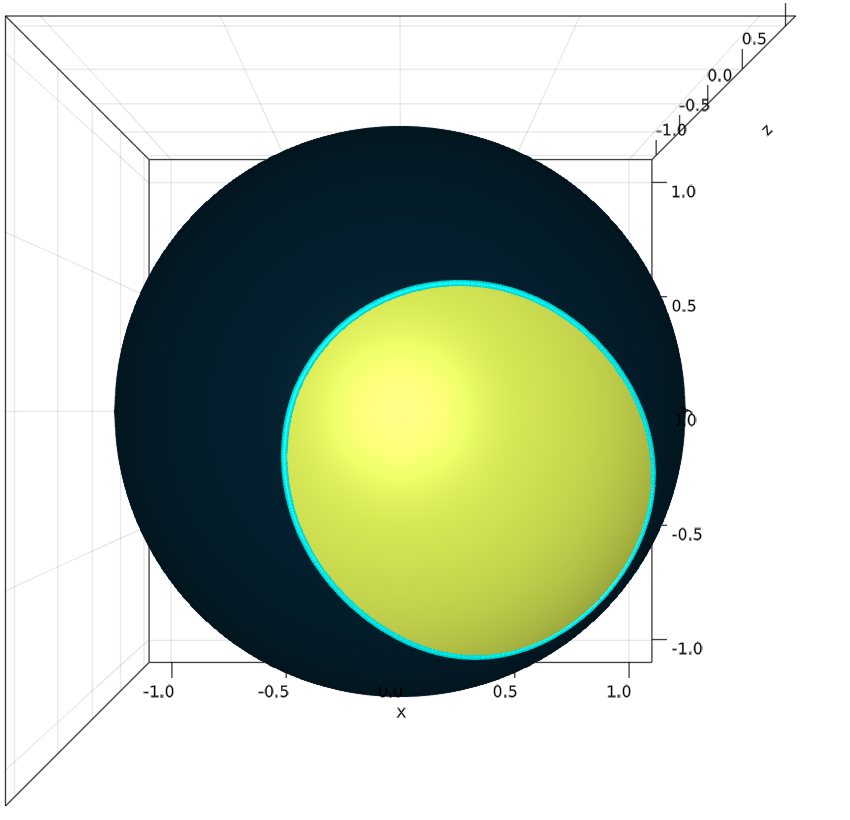}
                \subcaption{}

            \end{minipage} &
            \begin{minipage}[t]{0.25\columnwidth}
                \centering
                \includegraphics[height=0.15\vsize]{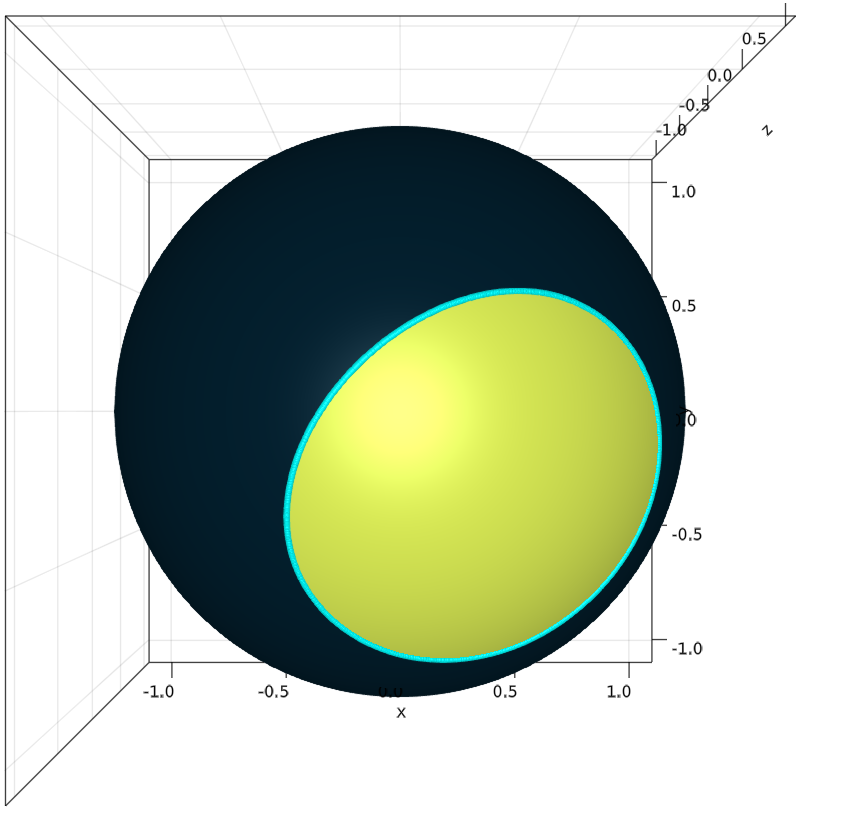}
                \subcaption{}

            \end{minipage}
            \\
            \begin{minipage}[t]{0.25\columnwidth}
                \centering
                \includegraphics[height=0.15\vsize]{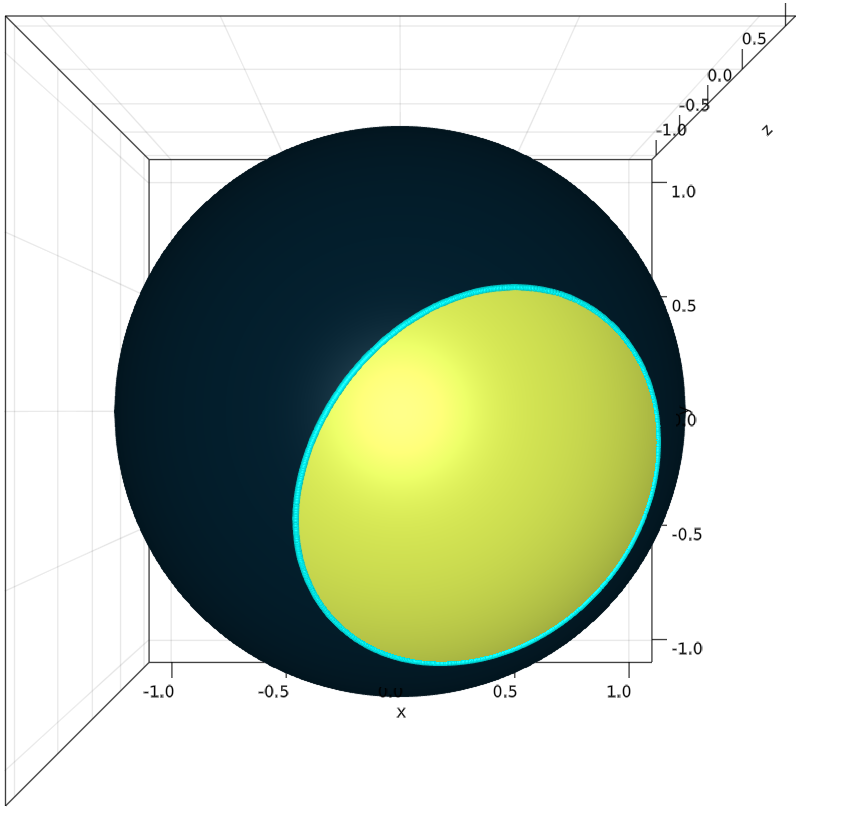}
                \subcaption{}

            \end{minipage} &
            \begin{minipage}[t]{0.25\columnwidth}
                \centering
                \includegraphics[height=0.15\vsize]{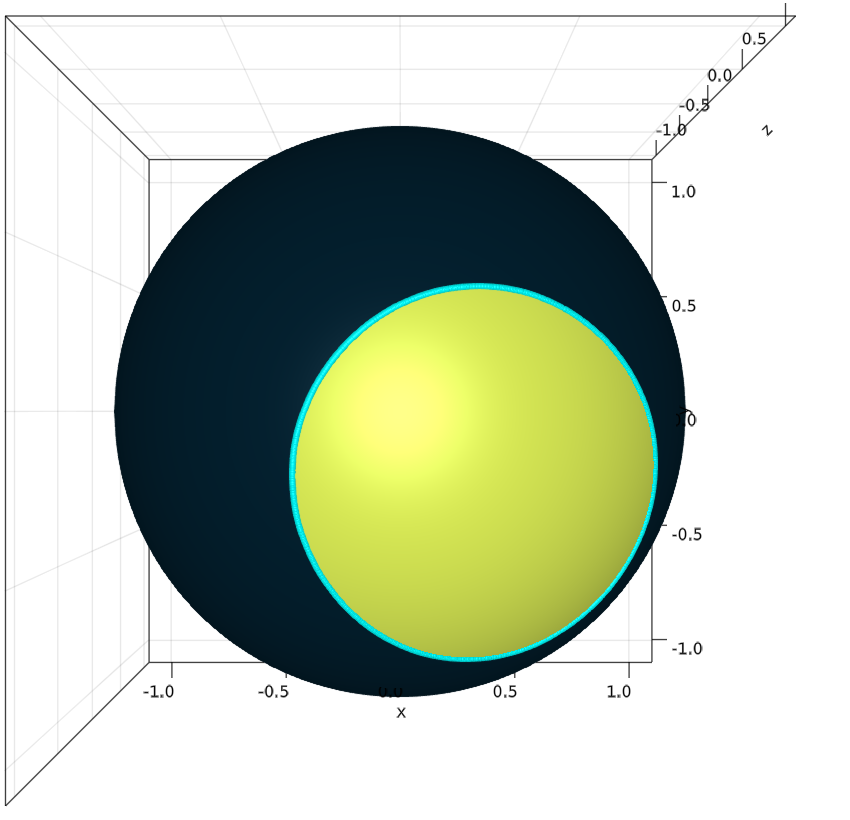}
                \subcaption{}

            \end{minipage}
            \\
            \begin{minipage}[t]{0.25\columnwidth}
                \centering
                \includegraphics[height=0.15\vsize]{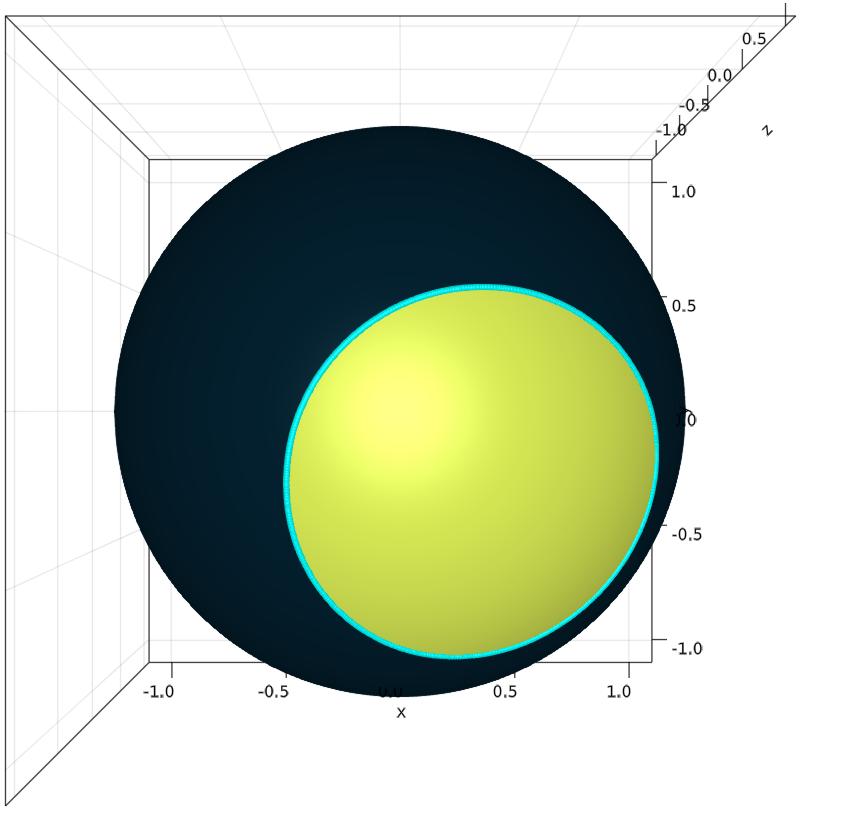}
                \subcaption{}
            \end{minipage} &
            \begin{minipage}[t]{0.25\columnwidth}
                \centering
                \includegraphics[height=0.15\vsize]{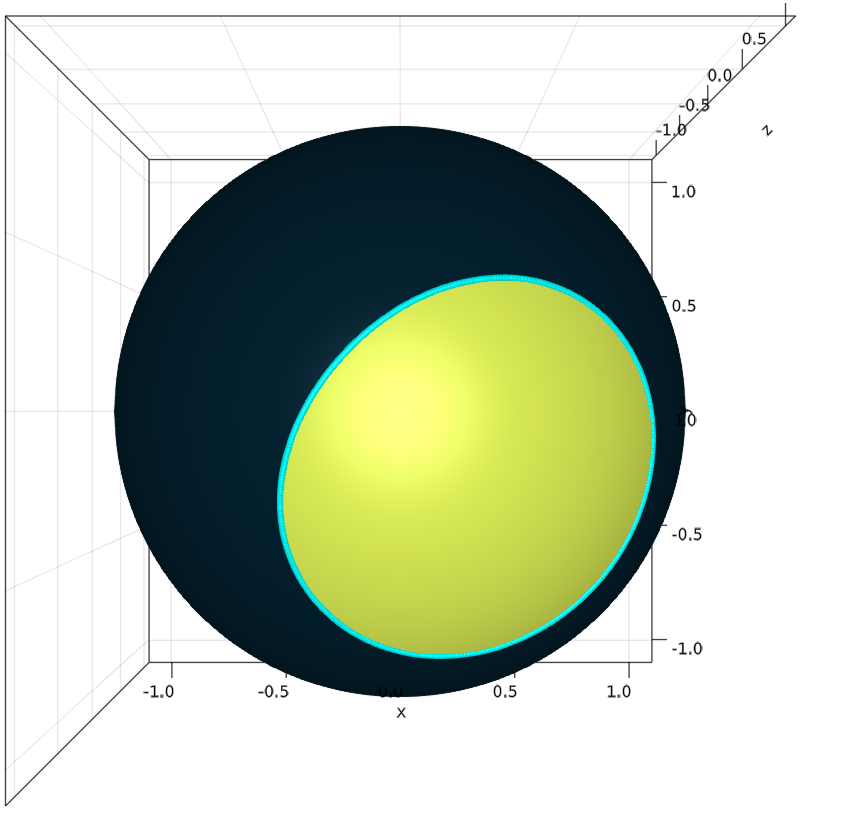}
                \subcaption{}
            \end{minipage}
        \end{tabular}
    }
    \caption{Result of two-phase hyperbolic MCF on the unit sphere (with area preservation). Arranged in alphabetical order. The time intervals from the initial time to the time the interface reaches a nearly steady state are equally spaced.}
    \label{res:sp2_hmcf_cons}
\end{figure}

\begin{figure}[H]
    \centering
    \fbox{
        \begin{tabular}{cc}
            \begin{minipage}[t]{0.25\columnwidth}
                \centering
                \includegraphics[height=0.15\vsize]{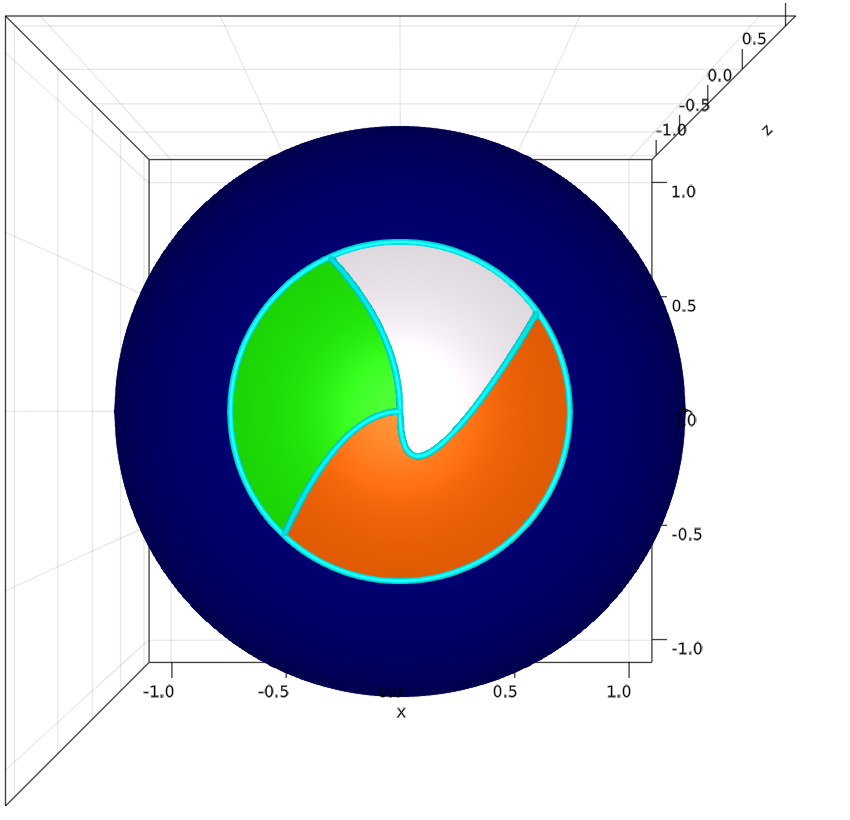}
                \subcaption{}

            \end{minipage} &
            \begin{minipage}[t]{0.25\columnwidth}
                \centering
                \includegraphics[height=0.15\vsize]{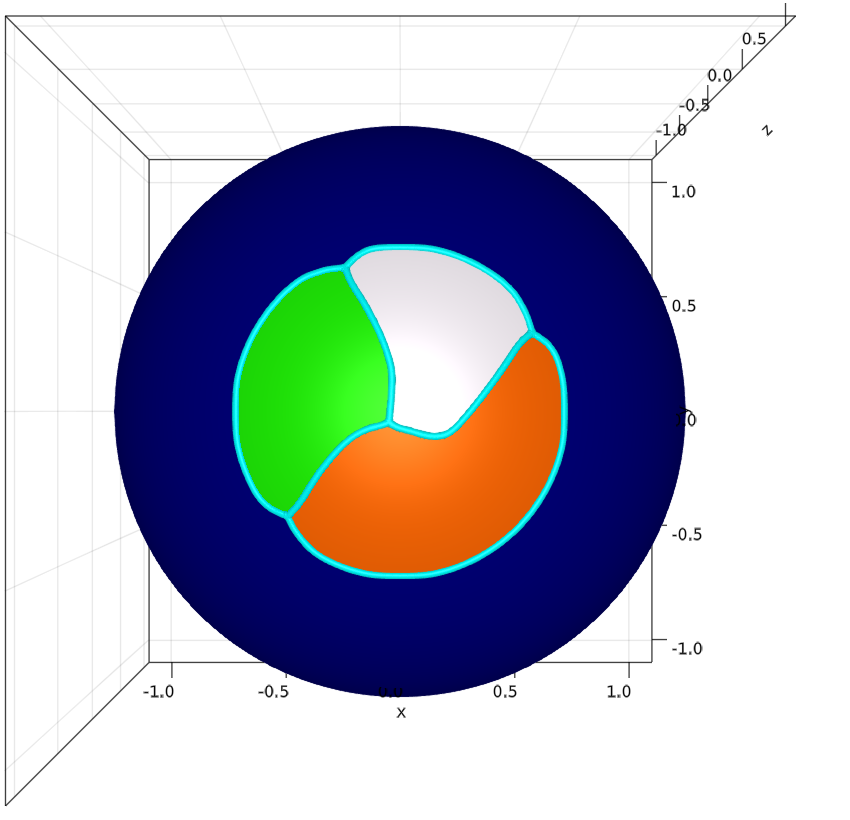}
                \subcaption{}
                \label{}
            \end{minipage}
            \\
            \begin{minipage}[t]{0.25\columnwidth}
                \centering
                \includegraphics[height=0.15\vsize]{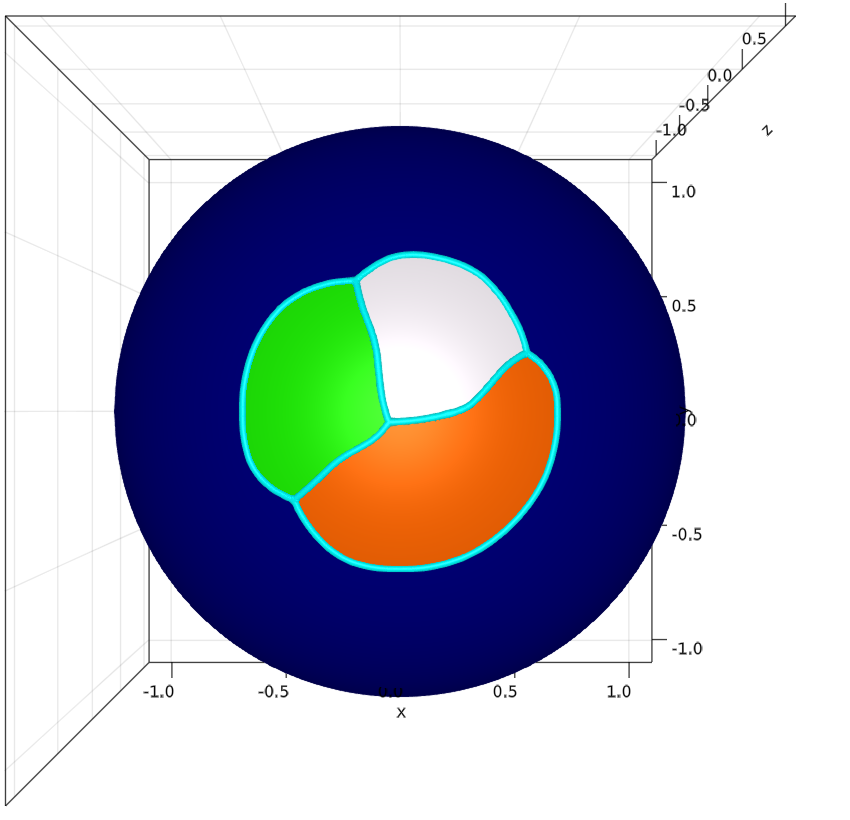}
                \subcaption{}
            \end{minipage} &
            \begin{minipage}[t]{0.25\columnwidth}
                \centering
                \includegraphics[height=0.15\vsize]{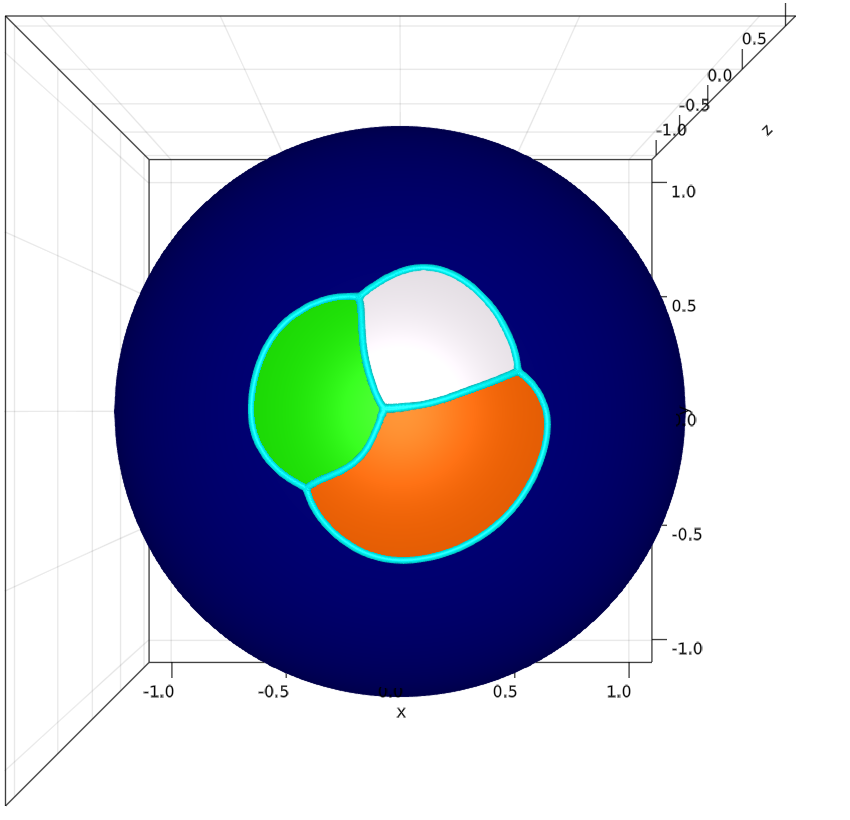}
                \subcaption{}
            \end{minipage}
            \\
            \begin{minipage}[t]{0.25\columnwidth}
                \centering
                \includegraphics[height=0.15\vsize]{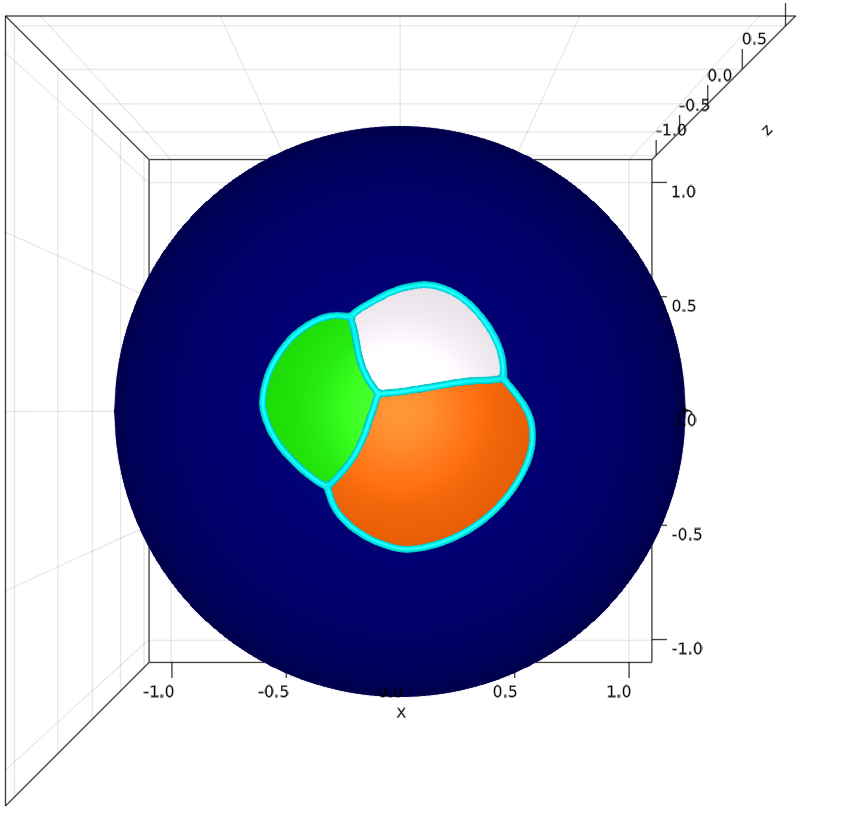}
                \subcaption{}
            \end{minipage} &
            \begin{minipage}[t]{0.25\columnwidth}
                \centering
                \includegraphics[height=0.15\vsize]{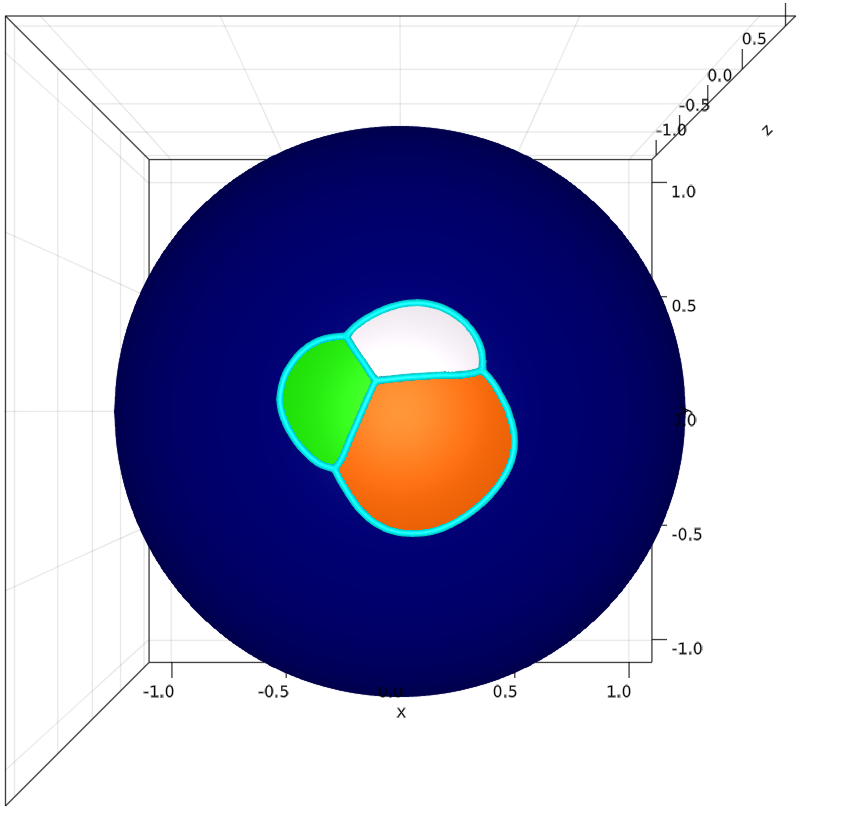}
                \subcaption{}
            \end{minipage}
            \\
            \begin{minipage}[t]{0.25\columnwidth}
                \centering
                \includegraphics[height=0.15\vsize]{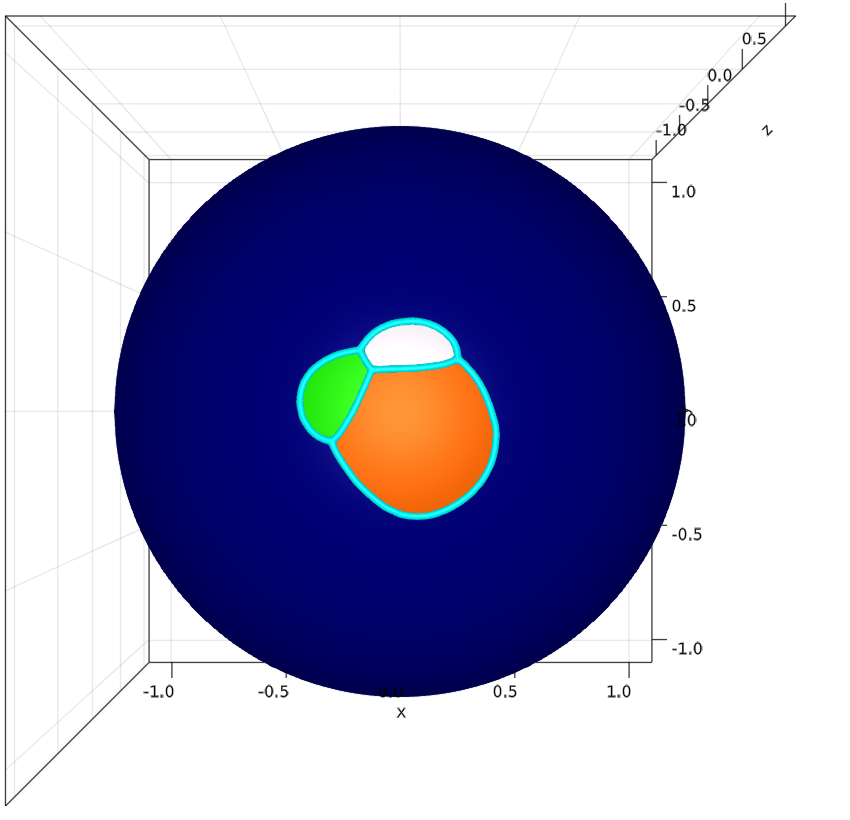}
                \subcaption{}
            \end{minipage} &
            \begin{minipage}[t]{0.25\columnwidth}
                \centering
                \includegraphics[height=0.15\vsize]{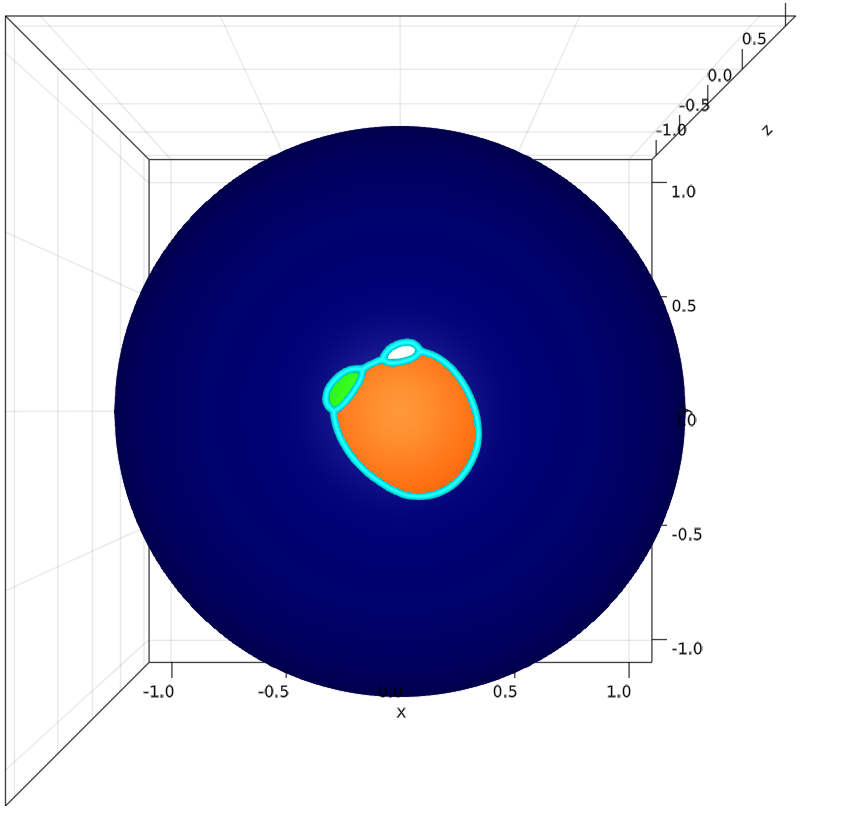}
                \subcaption{}
            \end{minipage}
            \\
            \begin{minipage}[t]{0.25\columnwidth}
                \centering
                \includegraphics[height=0.15\vsize]{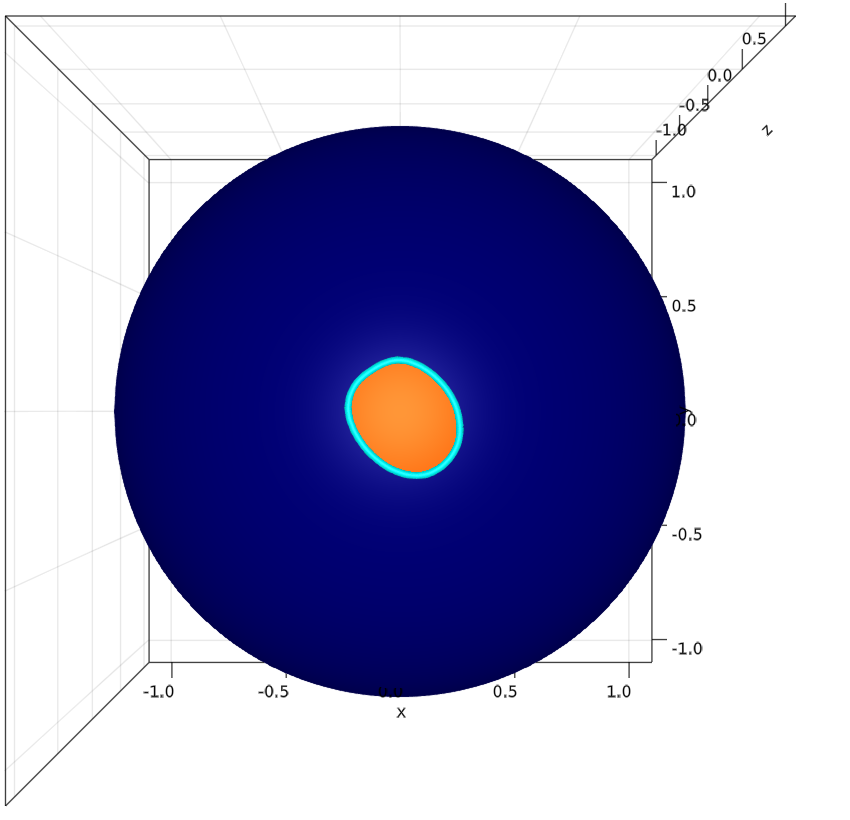}
                \subcaption{}
            \end{minipage} &
            \begin{minipage}[t]{0.25\columnwidth}
                \centering
                \includegraphics[height=0.15\vsize]{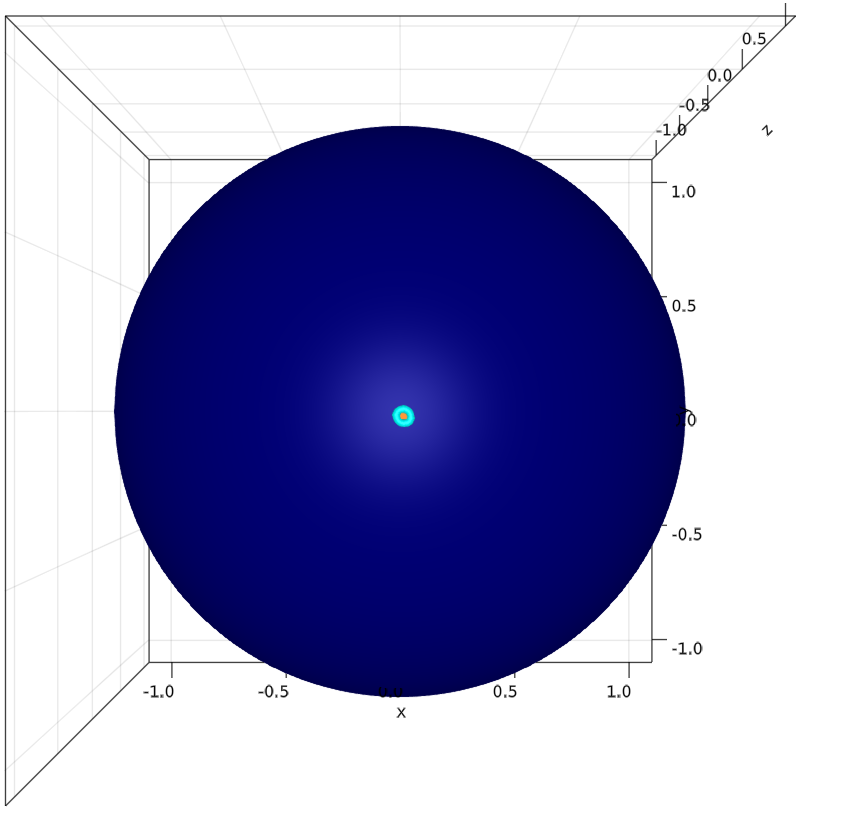}
                \subcaption{}
            \end{minipage}
        \end{tabular}
    }
    \caption{Result of four-phase hyperbolic MCF on the unit sphere (without area preservation). Arranged in alphabetical order with equal intervals from the initial time to the time the interface disappears.}
    \label{res:sp4_hmcf}
\end{figure}

\begin{figure}[H]
    \centering
    \fbox{
        \begin{tabular}{cc}
            \begin{minipage}[t]{0.25\columnwidth}
                \centering
                \includegraphics[height=0.15\vsize]{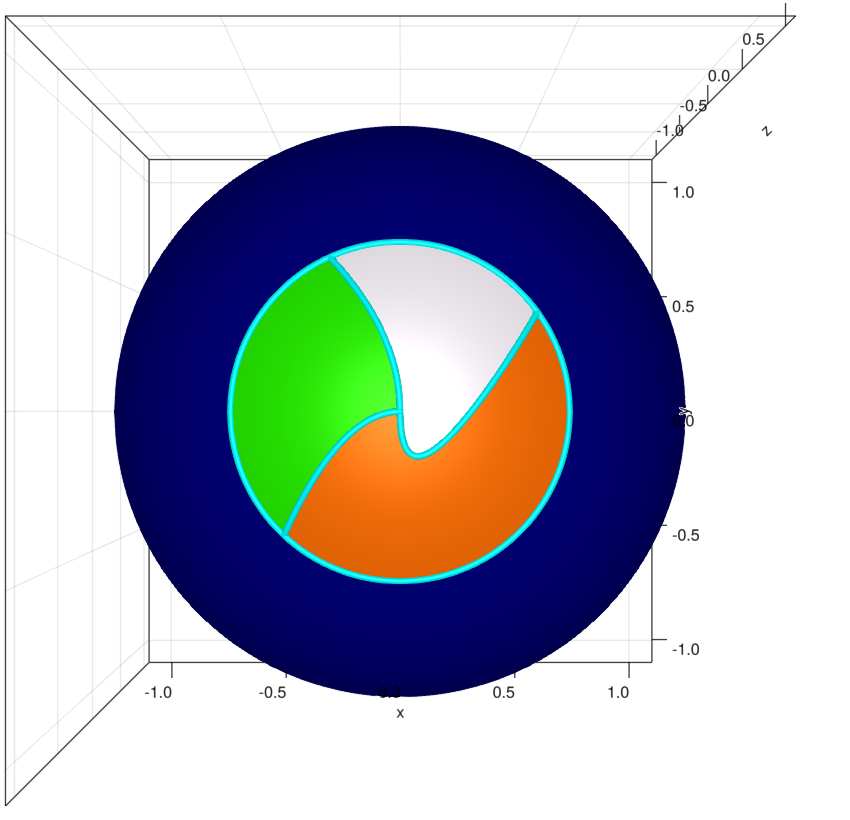}
                \subcaption{}

            \end{minipage} &
            \begin{minipage}[t]{0.25\columnwidth}
                \centering
                \includegraphics[height=0.15\vsize]{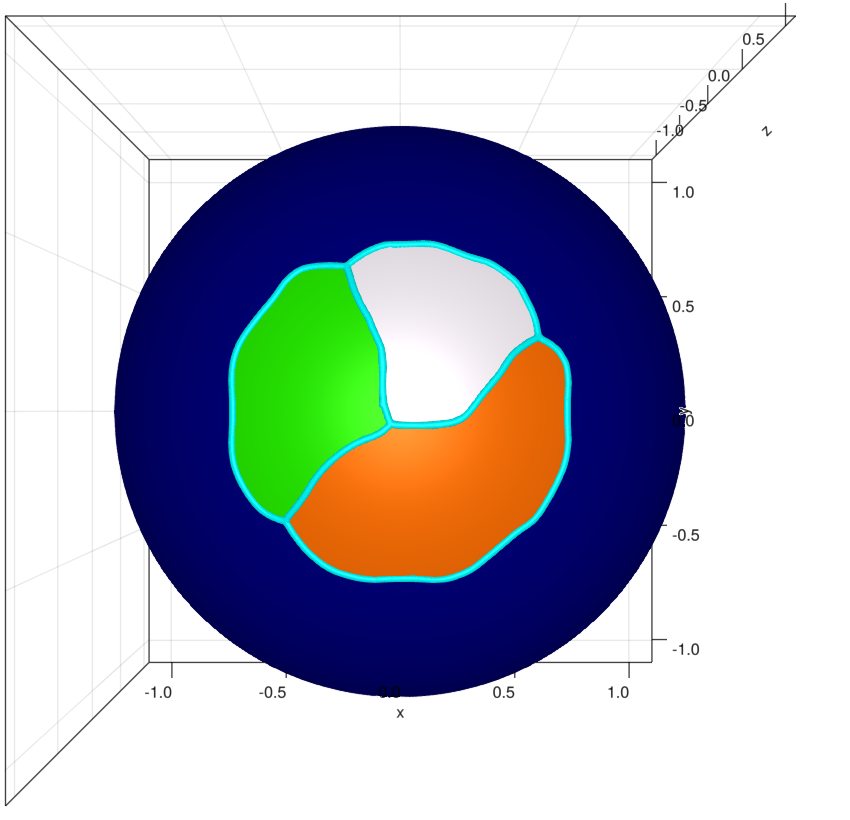}
                \subcaption{}
                \label{}
            \end{minipage}
            \\
            \begin{minipage}[t]{0.25\columnwidth}
                \centering
                \includegraphics[height=0.15\vsize]{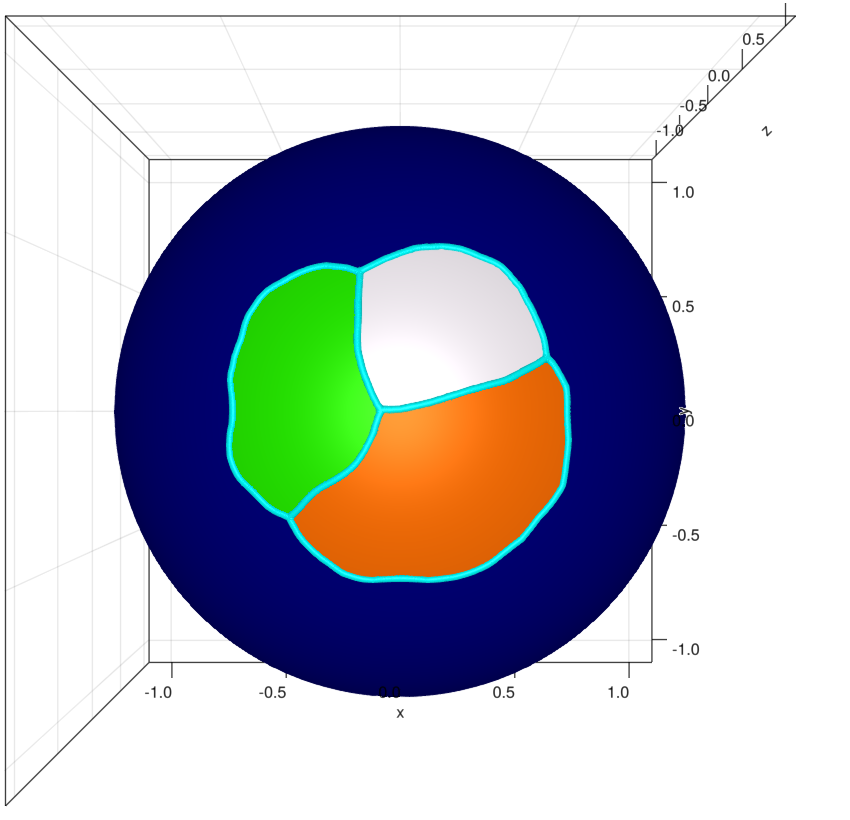}
                \subcaption{}

            \end{minipage} &
            \begin{minipage}[t]{0.25\columnwidth}
                \centering
                \includegraphics[height=0.15\vsize]{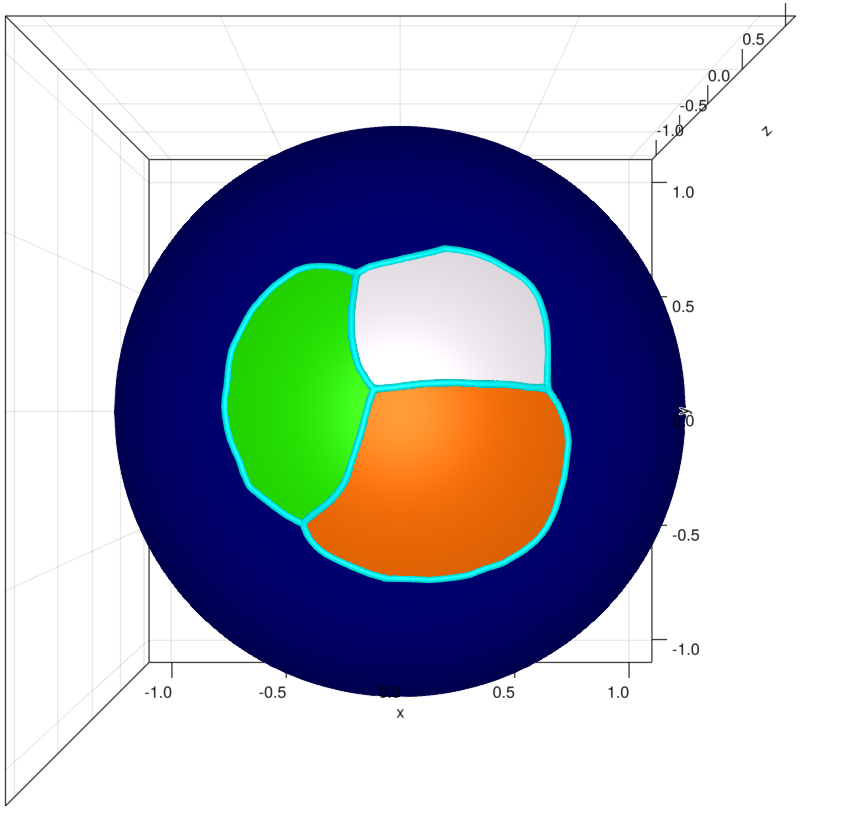}
                \subcaption{}

            \end{minipage}
            \\
            \begin{minipage}[t]{0.25\columnwidth}
                \centering
                \includegraphics[height=0.15\vsize]{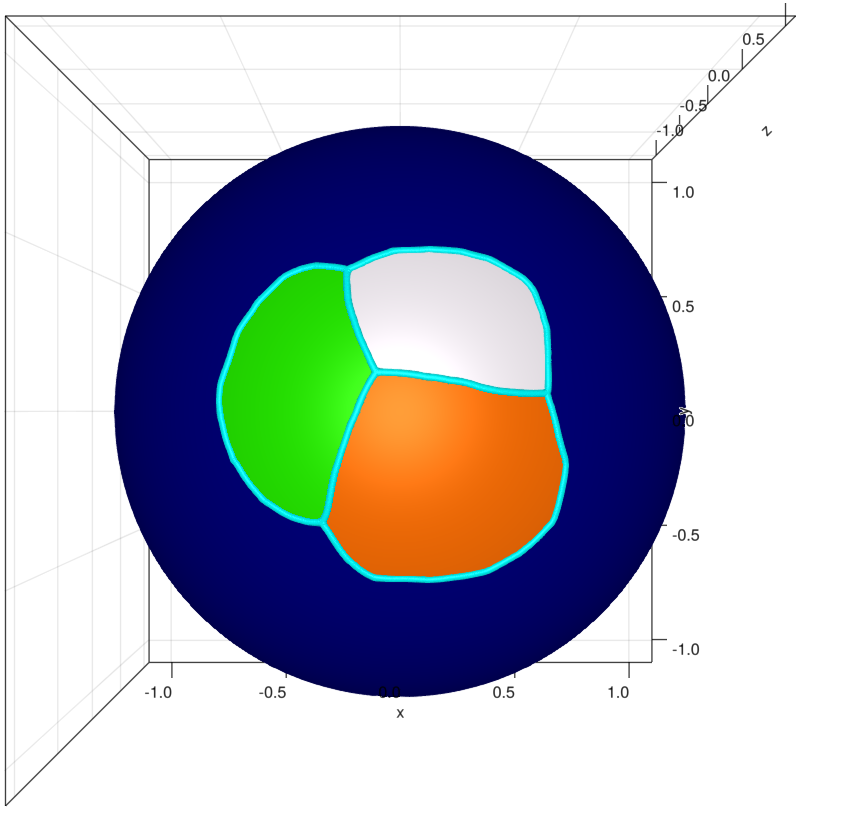}
                \subcaption{}

            \end{minipage} &
            \begin{minipage}[t]{0.25\columnwidth}
                \centering
                \includegraphics[height=0.15\vsize]{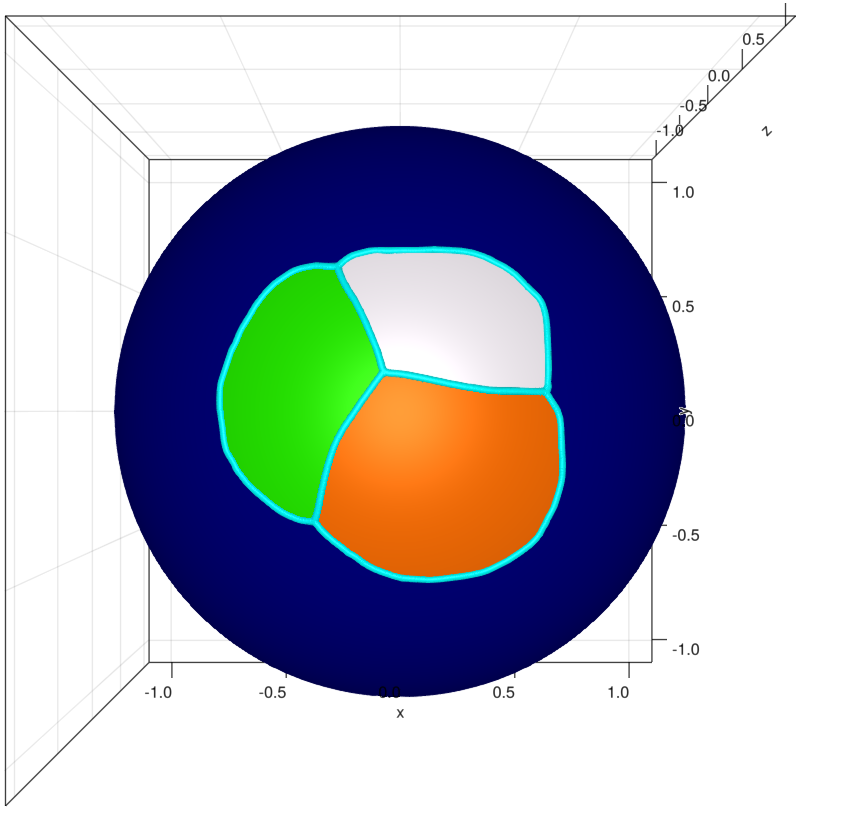}
                \subcaption{}

            \end{minipage}
            \\
            \begin{minipage}[t]{0.25\columnwidth}
                \centering
                \includegraphics[height=0.15\vsize]{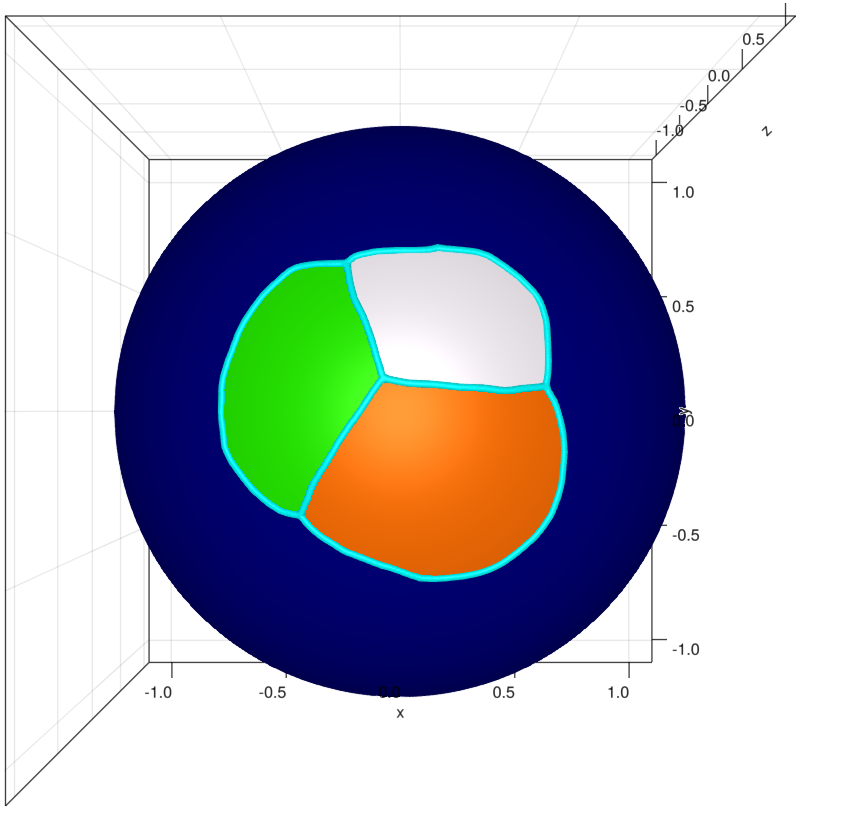}
                \subcaption{}

            \end{minipage} &
            \begin{minipage}[t]{0.25\columnwidth}
                \centering
                \includegraphics[height=0.15\vsize]{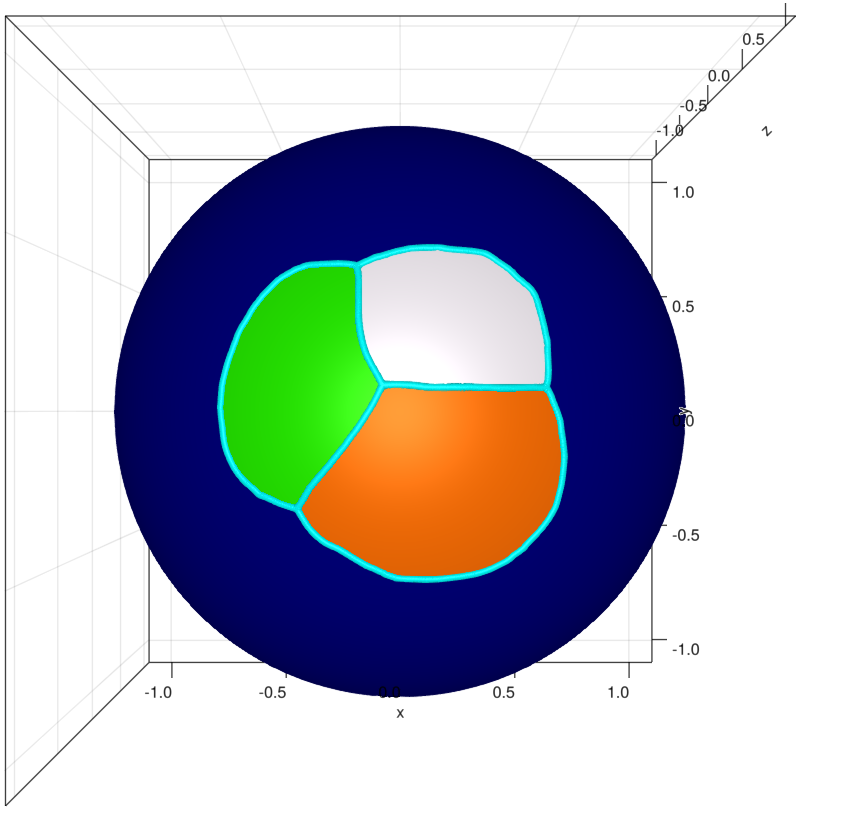}
                \subcaption{}

            \end{minipage}
            \\
            \begin{minipage}[t]{0.25\columnwidth}
                \centering
                \includegraphics[height=0.15\vsize]{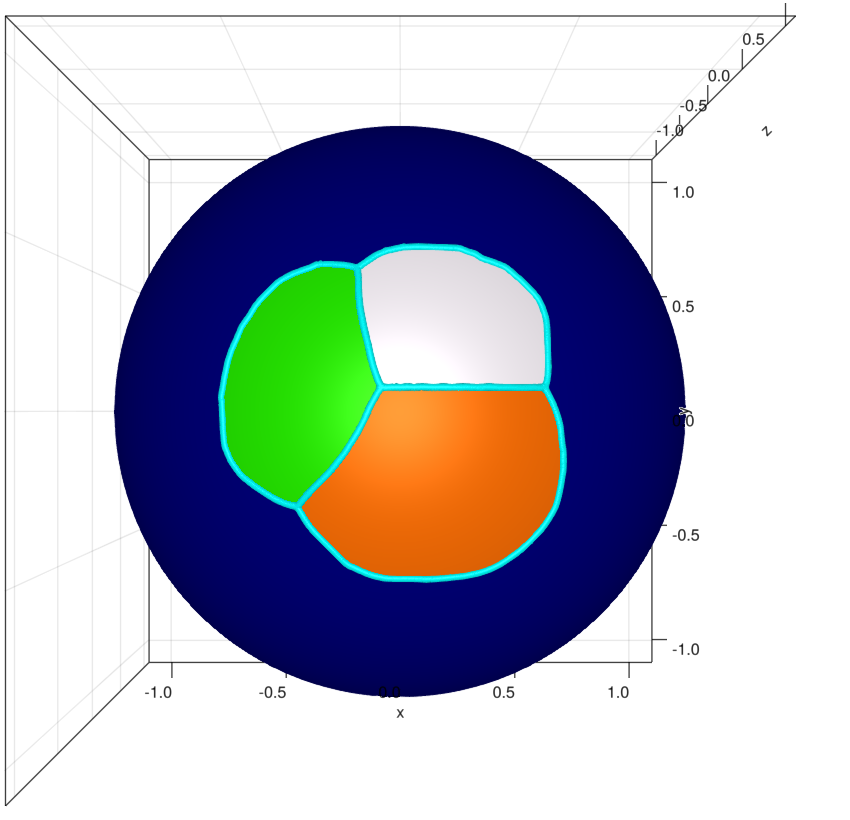}
                \subcaption{}

            \end{minipage} &
            \begin{minipage}[t]{0.25\columnwidth}
                \centering
                \includegraphics[height=0.15\vsize]{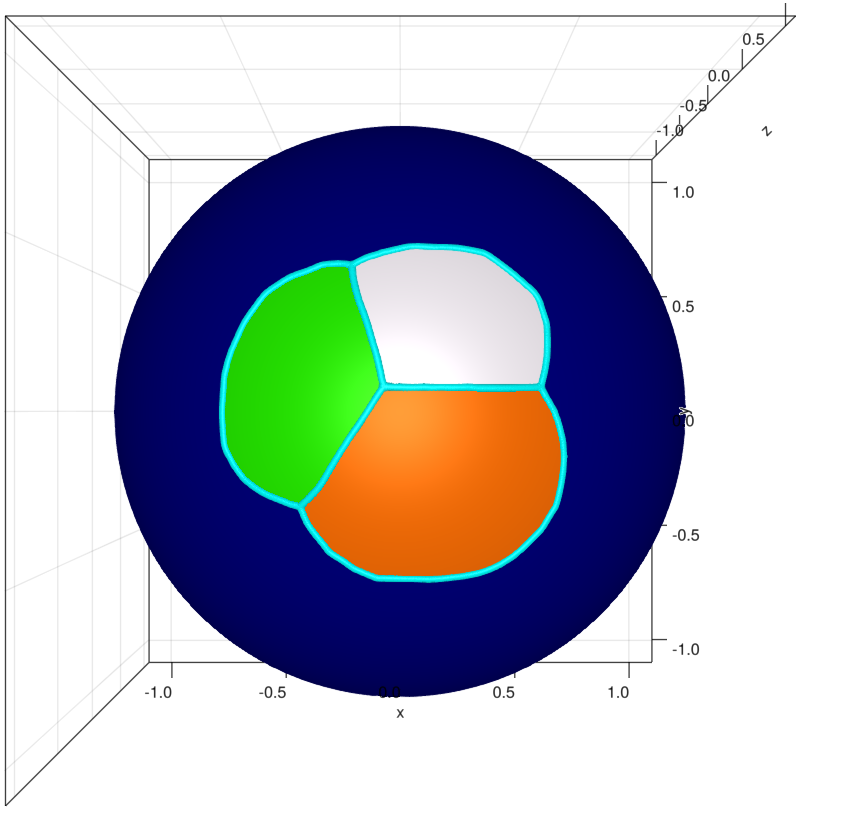}
                \subcaption{}

            \end{minipage}
        \end{tabular}
    }
    \caption{Result of four-phase hyperbolic MCF on the unit sphere (with area preservation). Arranged in alphabetical order. The time intervals from the initial time to the time the interface reaches a nearly stationary state are equally spaced.}
    \label{res:sp4_hmcf_cons}
\end{figure}

Here, we introduce the parameters used in the numerical calculation of hyperbolic mean curvature flow on surfaces. 
Numerical calculations were carried out in the two and four-phase setting on the unit sphere, both with and without area preservation. 
The discussion of the results is presented in section \ref{結果に対する考察}.
\begin{description}
    \item[Result of two-phase hyperbolic MCF on the unit sphere (without area preservation)]~\\
    We used the two-phase initial condition shown in Figure \ref{ini_sp_phase2} and the algorithm described in section \ref{algo-shmboc} for the numerical calculation. 
    Parameters were set as follows:
    $$
        \alpha=0.1, \quad\Delta x=0.05, \quad h= \Delta x,\quad\tau=5h,\quad \rho=0
    $$
    The numerical result is shown in Figure \ref{res:sp2_hmcf}.
    \item[Result of two-phase hyperbolic MCF on the unit sphere (with area preservation)]~\\

    We used the two-phase initial condition shown in Figure \ref{ini_sp_phase2} and the algorithm described in section \ref{algo-shmboc} for the numerical calculation. 
    Parameters were set as follows:
    $$
        \alpha=0.1,\quad \Delta x=0.05, \quad h= \Delta x,\quad\tau=5h,\quad \rho=10^3
    $$
    The numerical result is shown in Figure \ref{res:sp2_hmcf_cons}.
    \item[Result of four-phase hyperbolic MCF on the unit sphere (without area preservation)]~\\
    We used the four-phase initial condition shown in Figure \ref{ini_sp_phase4} and the algorithm described in section \ref{algo-shmboc} for the numerical calculation. 
    Parameters were set as follows:
    $$
        \alpha=1,\quad\Delta x=0.01, \quad h=7.84\Delta x/1000,\quad\tau=200h,\quad \rho=0
    $$
    The numerical result is shown in Figure \ref{res:sp4_hmcf}.
    \item[Result of four-phase hyperbolic MCF on the unit sphere (with area preservation)]~\\

    We used the four-phase initial condition shown in Figure \ref{ini_sp_phase4} and the algorithm described in section \ref{algo-shmboc} for the numerical calculation. 
    Parameters were set as follows:
    $$
        \alpha=1,\quad\Delta x=0.01, \quad h= 1.4\Delta x/100,\quad\tau=100h,\quad \rho=10^{6}
    $$
    The numerical result is shown in Figure \ref{res:sp4_hmcf_cons}.
\end{description}
\noindent

\subsection{Numerical error analysis of the area-preserving condition in the two-phase setting}\label{2相領域での面積保存条件の数値誤差解析}

We investigated how well the area enclosed by the interface in the two-phase setting on the unit sphere (shown in Figure \ref{ini_sp_phase2}), is preserved under the area-constrained mean curvature flow, and hyperbolic mean curvature flow. Here, we numerically solved the mean curvature flow and the hyperbolic mean curvature flow on the unit sphere with area preservation 
using the algorithm ``Surface MBO for multiphase regions with the area preservation condition" introduced in section \ref{algo-smboc}, 
and the algorithm ``Surface HMBO for multiphase regions with the area preservation condition" introduced in section \ref{algo-shmboc}, respectively. 
Note that setting $K=2$ in the above algorithms yields an approximation method for two-phase regions.

\begin{figure}[H]
    \begin{center}
        \fbox{\includegraphics[bb=0 0 354 249
                ,scale=0.85]{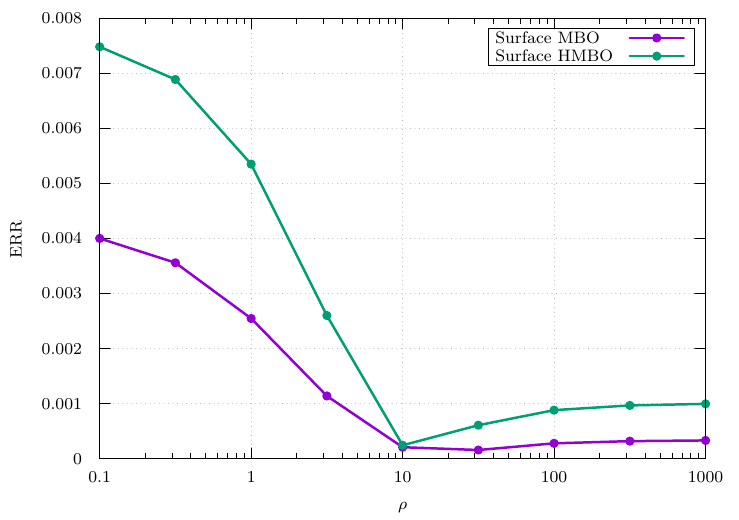}}
        \caption{Numerical errors of area preservation (two-phase), 
        The horizontal axis represents $\rho$, the vertical axis represents ERR defined in Equation \eqref{area_err}. The horizontal axis is shown on a logarithmic scale. 
        ``Surface MBO" corresponds to the two-phase MCF on the unit sphere, 
        and ``Surface HMBO" corresponds to the two-phase hyperbolic MCF on the unit sphere. 
        We confirm that, as $\rho$ increases, ERR tends to decrease. In both cases, we note that the error in the enclosed area is less than 0.001 for $\rho\geq10$.}\label{fig:mbohmbo_area}
    \end{center}
\end{figure}

Numerical errors were investigated as follows. 
The approximation of $A^1$ obtained 
in Step 3 of the algorithm ``Surface MBO for multiphase regions with the area preservation" and in Step 4 of ``Surface HMBO for multiphase regions with the area preservation" is denoted by $V_0$. 
The surface SDVF $z_{2}^{\epsilon,\lambda}$ obtained by executing the above algorithms for one step is used to calculate $V^1_{\boldsymbol{w}}$ using Eq. 
\eqref{cons_param_mbo}. Here, $z_{2}^{\epsilon,\lambda}$ is calculated in Steps 4.f and 5.f of the above algorithms, respectively. 
The solution obtained by executing the algorithm for one step corresponds to the time $\tau$. 
The approximate value of $V^1_{\boldsymbol{w}}$ is denoted by $V_\tau$ and the value of the error $\text{ERR}$ is defined as follows:
\begin{align}\label{area_err}
\text{ERR}=|V_0-V_\tau|.
\end{align}
We investigated the response of ERR to changes in $\rho$ (the penalty parameter for the area preservation) for the minimizing movements Eq. \eqref{mboF} and Eq. \eqref{F}. 
The parameters were as follows:

\begin{description}
    \setlength{\leftskip}{1.0cm}
    \item[Two-phase MCF on the unit sphere (with area preservation)]~\\
    $$
        \alpha=0.05,\quad\Delta x=0.05, \quad h= \Delta x^2/6,\quad\tau=100h,\quad 10^{-1}\le\rho\le10^3.
    $$
    \item[Two-phase hyperbolic MCF on the unit sphere (with area preservation)]~\\
    $$
        \alpha=0.1,\quad \Delta x=0.05, \quad h= \Delta x,\quad\tau=5h,\quad  10^{-1}\le\rho\le10^3.
    $$
\end{description}
The parameters used in the calculation of ``two-phase MCF on the unit sphere (with area preservation)'' were the same as those used in section \ref{results_MBO}, except for the values of $\rho$. 
The result obtained by evolving the system with $\rho=10^3$ is shown in Figure \ref{res:sp2_mcf_cons}. 
The parameters used in the calculation of ``two-phase hyperbolic MCF on the unit sphere (with area preservation)'' were the same as those used in section \ref{results_HMBO}, except for the value of $\rho$. 
The result obtained by evolving the system with $\rho=10^3$ is shown in Figure \ref{res:sp2_hmcf_cons}. 
Table \ref{tab:mbohmbo_area} shows the specific values of $\rho$ used for the numerical error analysis. 
The results are presented in Table \ref{tab:mbohmbo_area} and Figure \ref{fig:mbohmbo_area}, where the x-axis uses a logarithmic scale. 
Discussions of the results is presented in section \ref{結果に対する考察}.

\begin{table}[H]
    \caption{Numerical errors of area preservation (two-phase)}
    \label{tab:mbohmbo_area}
    \centering
    \begin{tabular}{lcc}
        \hline

        $\rho$      & ERR (MCF) & ERR (HMCF) \\

        \hline \hline
        $10^{-1}$   & 4.000e-03       & 7.482e-03             \\
        $10^{-0.5}$ & 3.557e-03       & 6.888e-03             \\
        $10^{0}$    & 2.545e-03       & 5.350e-03             \\
        $10^{0.5}$  & 1.138e-03       & 2.599e-03             \\
        $10^{1} $   & 2.056e-04       & 2.409e-04             \\
        $10^{1.5}$  & 1.553e-04       & 6.060e-04             \\
        $10^{2}$    & 2.769e-04       & 8.801e-04             \\
        $10^{2.5}$  & 3.166e-04       & 9.653e-04             \\
        $10^{3}$    & 3.299e-04       & 9.934e-04             \\
        \hline
    \end{tabular}
\end{table}
\subsection{Discussion}\label{結果に対する考察}
In this section, we will explain the results of the numerical calculations conducted in sections \ref{results_MBO} to \ref{2相領域での面積保存条件の数値誤差解析}. Our results are summarized in the following order: surface MCF, surface hyperbolic MCF, and numerical error analysis of the area preservation conditions in the two-phase setting. Following this, we discuss our observations regarding the motion of interfaces in the simulations, properties of the functionals used in our approximation methods, and issues related to area preservation.

Numerical results of surface MCF involving interfaces in the two-phase setting (Figure \ref{res:sp2_mcf}) without area preservation show that curves on the unit sphere disappear over time. Similar to the flat setting, their length decreases, and the interface becomes nearly circular before contracting to a single point. 
When the area preservation condition (Figure \ref{res:sp2_mcf_cons}) is prescribed, the curvature of the interfaces tend to decrease over time, and the interface converges to a circular shape. 
The area remains approximately constant, and the curve approaches a stationary state at the terminal time. 
In the four-phase setting (Figure \ref{res:sp4_mcf}) without area preservation, the interface smooths itself while the length of the network decreases with time. Each interface moves to maintain the junctions where they intersect before contracting to a single point. 
When the area preservation condition is prescribed (see Figure \ref{res:sp4_mcf_cons} and Figure \ref{res:sp4_mcf_cons_2}), the curvature of the interfaces tend to decrease over time, while the junctions are maintained. However, even in the steady state, we note some slight irregularities near the junction (see Figure \ref{res:sp4_mcf_cons}). Figure \ref{res:sp4_mcf_cons_2} confirms that the curvature of the curve is still relatively large. The areas of each phase changed slightly compared to the initial state. Both observations can be attributed to various factors including: ambient mesh spacing, interpolation methods, and optimization stopping criteria.

Next, we explain the numerical results for the hyperbolic mean curvature flow on surfaces. In the two-phase setting without area-preservation (Figure \ref{res:sp2_hmcf}), the length of the curve decreases while oscillating and  approaching a circular shape, before contracting to a point. Under the area-preservation condition (Figure \ref{res:sp2_hmcf_cons}), the shape of interface converges to a circle over time and the area remained approximately constant. After becoming approximately circular, the curve continued to move atop the surface of the sphere.

In the four-phase setting, the numerical results of the hyperbolic mean curvature flow (Figure \ref{res:sp4_hmcf}) showed that, without the area-preservation condition, the interfaces oscillated with time while the total length of the network decreased. 
The interfaces evolved while maintaining junctions, and eventually the network contracted to a single point. 
When the area-preserving condition was applied, the interfaces oscillated while while preserving their areas and maintaining junctions over time. 
The area of each phase remained almost constant throughout the evolution, and eventually reached a nearly stationary state.

In the two-phase setting, surface MCF and hyperbolic MCF with area preservation both showed a tendency for the error ERR (equation \eqref{area_err}) to decrease as $\rho$ increases.
In Figure \ref{fig:mbohmbo_area}, ERR decreased as $\rho$ increased over the range $10^{-1}\le\rho\le10$, and after $\rho>10$, ERR increased and became almost constant and less than 0.001.
Overall, we observe that if we approximate the area-preserving MCF and hyperbolic mean curvature flow using the proposed numerical method, the value of ERR will decrease as $\rho$ increases.
However, it is expected that increasing $\rho$ beyond a certain value will not reduce ERR below a certain threshold.

In the numerical results (Figure \ref{res:sp4_mcf_cons} and Figure \ref{res:sp4_mcf_cons_2}) for the interfaces moving according to surface MCF with area preservation in the multiphase setting, a slight irregularity was observed at the point of stationary state. 
This observation indicates that there are points on the interface with relatively large curvature, which is contrary to the expected result. 
The reason for this stagnation, similar the original MBO, can be attributed to the fact that even though there are points with large curvature along the interface, the curvature may still be too small to resolve for a given threshold length.
Relatedly, setting an interface with sufficiently large curvatures along its initial curve tends to eliminate points with large curvature at the steady state. 
Alternatively, using a smaller spatial discretization tends to alleviate such constrains on the interface's motion. This, of course, leads to an increase in the computational time required by the method. 
In fact, all the methods developed in this study require a relatively fine spatial grid, and refining it further causes a significant increase in the required computation time. 
For example, changing the spatial grid width from 0.05 to 0.01 for the method ``Surface MBO for multiphase regions with the area preservation'' described in section \ref{algo-smboc} increases the require computational time by a factor of 60. 
Consequently, improving the methods used in this study (especially their computation time) is an important future task.

Regarding the surface hyperbolic MCF, the oscillation of the interface tends to decrease with time (Figure \ref{res:sp2_hmcf_cons} and Figure \ref{res:sp4_hmcf_cons}). From the point of view of conservation of energy, this observation is unexpected. One of the reasons for this is thought to be the use the MM used in the algorithm created.
In MM, it is known that the energy of the obtained numerical solution decreases compared to the exact solution as time increases \cite{newbook}.
Since the algorithm created reconstructs the interface based on the numerical solution obtained from the MM, it is understandable that the kinetic energy ofthe numerical solution of the interface decreases with time.

In this research, we have used minimizing movements to impose area constraints on surface-constrained multiphase interfacial motions.  Since minimizing movements require one to minimize a functional, we have also considered the influence of the functional used in this process and on the corresponding numerical results. Again, the functional used in our method for dealing with the area-constrained curvature flows in the multiphase setting (Equations \eqref{mboF} and \eqref{F}) is expressed as follows:
\begin{align}\label{F_}
    \mathcal{F}_m(\boldsymbol{w})=\int_{\Omega_\lambda}\left(F(\boldsymbol{w},\boldsymbol{w}_{m-1},\boldsymbol{w}_{m-2})+\alpha\frac{|\nabla\boldsymbol{w}|^2}{2}\right)d\boldsymbol{x}
    +\rho\sum^K_{i=1}|A^i-V^i_{\boldsymbol{w}}|^2
\end{align}
where $V^i_{\boldsymbol{w}}$ is the $V^i_{\boldsymbol{w}}$ included in Equations \eqref{mboF} and \eqref{F}, and $F$ is expressed as follows.
\begin{align}\label{F_F}
    \displaystyle
    F(\boldsymbol{w},\boldsymbol{w}_{m-1},\boldsymbol{w}_{m-2})=
    \begin{dcases}
        \displaystyle
        \frac{|\boldsymbol{w}-\boldsymbol{w}_{m-1}|^2}{2h},\quad                         & \text{Surface MBO}  \\
        \displaystyle
        \frac{|\boldsymbol{w}-2\boldsymbol{w}_{m-1}+\boldsymbol{w}_{m-2}|^2}{2h^2},\quad & \text{Surface HMBO}
    \end{dcases}
\end{align}
For simplicity, let us define
\begin{align}
    J(\boldsymbol{w}) & =\alpha\frac{|\nabla\boldsymbol{w}|^2}{2}\label{JJ}     \\
    P(\boldsymbol{w}) & =\rho\sum^K_{i=1}|A^i-V^i_{\boldsymbol{w}}|^2\label{PP}
\end{align}
Then, the functional in equation \eqref{F_} can be expressed as follows:
\begin{align}\label{F_F_F}
    \mathcal{F}_m(\boldsymbol{w})=\int_{\Omega_\lambda}\left(F(\boldsymbol{w},\boldsymbol{w}_{m-1},\boldsymbol{w}_{m-2})+J(\boldsymbol{w})\right)d\boldsymbol{x}
    +P(\boldsymbol{w})
\end{align}
The functional $P$ determines the penalty of to the area constraints. If we set $P(\boldsymbol{w})=0$, then equation \eqref{F_F_F} corresponds to the functional used in the absence of area preservation.
Emphasis of each area constraint is controlled through the value of $\rho$. 
However, if $\rho$ is taken too large so that $F$ and $J$ are significantly smaller than $P$, then the minimizing scheme will tend to focus only on the penalty term. That is, during the process of minimizing the functional, the significance of $F$ and $J$ are diminished when compared to that of $P$.
As a result, the approximate solution may deviate from the expected result. One may observe a jagged interface, even after several minimizers and at the stationary state. 
On the other hand, if $P$ is too small, the area constrain of each phase will not be satisfied at an acceptable level. 
Therefore, in order to approximate the motion of each interface following the mean curvature flow or the hyperbolic mean curvature flow while satisfying the area constraints, it may be necessary to adjust the ratio of the magnitudes of $F$, $J$, and $P$.  
Such an approach is would avoid large differences in the magnitudes of $F$, $J$, and $P$. 
However, it is not clear what ratio the magnitudes of $F$, $J$, and $P$ should satisfy at present. 
We would like to return to this and related topics in a future study.

In imparting the area constraints a top the surface $S$, we numerically solved the constrained partial differential equations in the tubular region $\Omega_\lambda$.
The functionals used in the surface-type MM (Eq. \ref{mboF}, Eq. \ref{F}) were designed to conserve volume in $\Omega_\lambda$, where the width of $\Omega_\lambda$ is a constant $\lambda$ (Eq. \ref{getband_cont}). 

\begin{figure}[H]
    \begin{center}
        \fbox{\includegraphics[bb=0cm 0cm 11cm 6cm,scale=0.85]{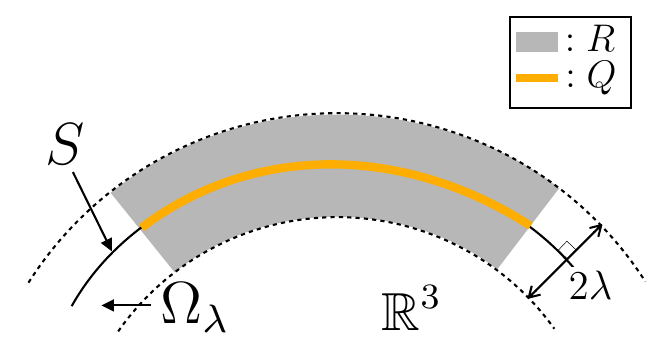}}
        \caption{A cross-section of the surface $S$. Here, $S$ denotes the surface of interest, $\Omega_\lambda$ is the region around $S$ (see definition in equation \eqref{getband_cont}), $Q$ is the region enclosed by the interface on the surface $S$, and $R$ is the region obtained by extending $Q$ in the normal direction of $S$ to $\Omega_\lambda$.}\label{fig:bandingRQ}
    \end{center}
\end{figure}

Consider the case of an interface in the two phase setting. 
Let $Q$ be the region enclosed by the interface on the surface $S$, let $A$ be the area of $Q$, let $R$ be the region obtained by extending $Q$ in the normal direction of the surface $S$ to $\Omega_\lambda$, and let $V$ be the volume of $R$. 
Figure \ref{fig:bandingRQ} shows a schematic diagram of the relationships between $S$, $\Omega_\lambda$, $R$, and $Q$.
In this case, assuming that $\lambda$ is sufficiently small, we note that we can approximate $V$ as
\begin{align}\label{2alambda}
    V\approx2A\lambda
\end{align} 
In section \ref{2相領域での面積保存条件の数値誤差解析}, we investigated the numerical error of the area preservation for two-phase regions. 
We observed that the mean curvature flow and the hyperbolic mean curvature flow conserve area at higher precisions as $\rho$ is increased. 
Since $V$ is approximated by Eq. \eqref{2alambda} and $\lambda$ is a constant, it is expected that for two-phase regions on a surface, increasing $\rho$ will better conserve the area $A$ surrounded by the interface on the surface. 
\\
{\bf{Remark:}}
The value of the width $\lambda$ of $\Omega_\lambda$ used in Section \ref{2相領域での面積保存条件の数値誤差解析} is given by
\begin{align*}
    \lambda & =\sqrt{17}\Delta x \\
            & \approx 0.2
\end{align*}
This is obtained by substituting $p=3$ and $\Delta x=0.05$ into equation \eqref{lambda}.

The computational results in sections \ref{results_MBO} and \ref{results_HMBO} showed that area preservation can be approximately satisfied on surfaces even in the multiphase case.
However, we have not yet performed a numerical error analysis to describe the relationship between the parameters of the computational algorithm, and the area preservation condition for cases other than the two-phase setting. We would like to treat this in a separate study.



\section{Summary}\label{まとめ}
This study developed approximations methods for surface-constrained mean curvature flow and hyperbolic mean curvature flow of interfaces.
This was achieved by first creating approximation methods for a surface partial differential equations by combining the closest point method with minimizing movements. We then extended the methods to implement the conventional MBO and HMBO algorithms on surfaces.
In addition, we constructed the surface-signed distance vector field to distinguish multiphase geometries on surfaces.

Numerical error analyses of our methods were performed for the surface heat and wave equations, and convergence with respect to the spatial discretization was investigated. It was found that the numerical solution of the partial differential equation on the surface obtained by our  approximation methods converges to the exact solution.

By using the surface version of the signed distance vector field, we extended the MBO and HMBO algorithms to the surface-constrained setting. These were used to perform numerical calculations of mean curvature flow and hyperbolic mean curvature flow for two and four phase interfacial motions.

The numerical error of the prescribed area in the two-phase setting for mean curvature flow and hyperbolic mean curvature flow on surfaces was evaluated. Our results confirm that increasing the value of the penalty parameter $\rho$ leads to higher precision in the area preservation.

Improvements to our approximation methods could be made by adjusting the energy functionals used in the MM method. Namely, it is known that, by using appropriate functionals, energy conservation can be realized \cite{newbook}. Therefore, creating approximation methods that conserve energy and performing their numerical error analysis for equations such as the surface-wave equation is an important future task. This is expected to clarify questions about the energy dissipation of the interface in the HMBO algorithm. In addition, we would like to design generalized surface-type threshold dynamics which impart damping terms on target interfacial motion.

\newpage
\section{Appendix}

\subsection{Surface HMBO and Initial Velocity for Multiphase Regions}\label{多相領域に対するSurface HMBOアルゴリズムと初期速度}

In the surface HMBO for multiphase regions (see Section \ref{algo-mHMBO} and Section \ref{algo-shmboc}) one needs to determine $\boldsymbol{P}_{1}$ from the initial shape $\boldsymbol{P}_0$ and initial velocities $\boldsymbol{v}_0^i$ of each phase. Here we describe the method.
\noindent
Let $\Gamma_{ij}$ be the interface between phase $i$ and phase $j$. We define the following:
$$
    \Gamma=\bigcup_{i,j}\Gamma_{ij}, \quad \mathcal{F}_{ij}=\Gamma-\Gamma_{ij}, \quad \mathcal{J}_{ij}=\Gamma_{ij}\cap \Sigma_{\Gamma_{ij}}\\ 
$$
$$
    \Sigma_{\Gamma_{ij}}=\{C|\{\Gamma_{ij}\cap C\}\neq\varnothing , C \in \mathcal{F}_{ij}\}
$$
$\Gamma$ represents the union of all the interfaces, $\mathcal{F}_{ij}$ represents the interfaces other than $\Gamma{ij}$, and  $\mathcal{J}_{ij}$ represents the endpoints of $\Gamma{ij}$. Also, $\Sigma_{\Gamma_{ij}}$ represents the set of interfaces connected to $\Gamma_{ij}$.
In the multiphase setting, the surface HMCF is represented by the following nonlinear partial differential equation. For each interface $\Gamma_{ij}$, we have:
\begin{align}\label{shmcf_eq}
    \displaystyle
    \begin{cases}
        \displaystyle{\frac{d^2 \Gamma_{ij}}{dt^2}=-\kappa_{ij}\boldsymbol{n}_{ij}},\quad & t>0                                                                        \\[8pt]
        \Gamma_{ij}(t=0)=\Gamma_{ij}^0                                                                                     \\[8pt]
        \displaystyle{\frac{d\Gamma_{ij}}{dt}(t=0)}=V_{ij}^0\boldsymbol{n}_{ij}^0                                                                                      \\[8pt]
        \Gamma_{ij}(P)=\sigma(P), \quad                                                   & t\geq 0,\quad \sigma\in\Sigma_{\Gamma_{ij}}, \quad P \in \mathcal{J}_{ij}.
    \end{cases}
\end{align}
Here, $\kappa_{ij}$ denotes the mean curvature of $\Gamma_{ij}$, $\boldsymbol{n}_{ij}$ represents the outward unit normal vector of $\Gamma{ij}$, $\Gamma_{ij}^0$ represents the interface between phases $i$ and $j$ at the initial time, $V_{ij}^0$ represents the initial velocity of $\Gamma_{ij}^0$, and $\boldsymbol{n}_{ij}^0$ represents the outward unit normal vector of $\Gamma_{ij}^0$.
The fourth equation in Equation \eqref{shmcf_eq} represents the continuity condition that the interface $\Gamma_{ij}$ should satisfy.
Without the continuity condition, each interface may move independently over time, which could lead to the loss of junctions.

In the Surface HMBO for multiphase regions, after determining the initial shape $\boldsymbol{P}_0$, the interface set ${\Gamma_{ij}^0}$ is defined for each interface. Then, for each interface, the set $\{\Gamma_{ij}^1\}$ is defined as follows:
$$\Gamma_{ij}^1=\Gamma_{ij}^0+hV_{ij}^0\boldsymbol{n}_{ij}^0$$
Following this, the regions $\boldsymbol{P}_1$ for each phase are determined from the set $\{\Gamma_{ij}^1\}$.
\subsection{Numerical error analysis of the surface MBO for two-phase regions}\label{2相領域に対するSurface MBOアルゴリズムの数値誤差解析}
\begin{figure}[H]
    \begin{center}
        \fbox{\includegraphics[bb=0cm 0cm 8cm 8cm,scale=0.8]{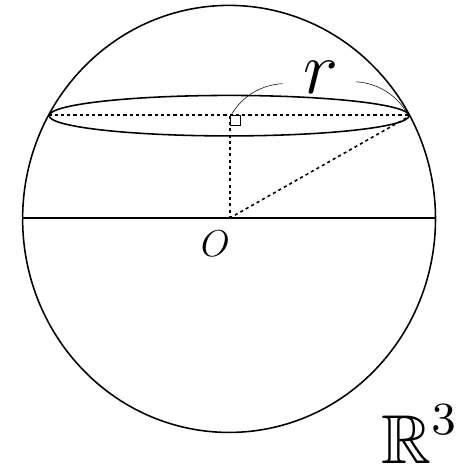}}
        \caption{A circle with radius $r$ on the unit sphere}\label{circleonsphere}
    \end{center}
\end{figure}

In this section, we present the results of an numerical error analysis for the mean curvature flow on surfaces using the algorithm introduced in section \ref{algo-2MBO}.
The analysis focuses on the motion of a circle a top the unit sphere.
As shown in Figure \ref{circleonsphere}, let $r$ denote the radius of the circle on the unit sphere.
A circular interface moving by the mean curvature flow on the surface of the unit sphere satisfies the following differential equation: \cite{mc_surf}:
\begin{align}\label{mcf_surf_de}
    \begin{cases}
        \displaystyle{\frac{dr}{dt}=\frac{r^2-1}{r}},\quad t>0 \\[8pt]
        \displaystyle{ r(0)=r_0}
    \end{cases}
\end{align}
where $r$ denotes the radius of the circle at time $t$ and $r_0>0$ denotes the radius of the circle at the initial time.
The exact solution of Eq. \eqref{mcf_surf_de} can be obtained as follows:
\begin{align}
    r(t)=\sqrt{1-(1-r_0^2)\exp(2t)}
\end{align}
We investigated the numerical error's dependence on the grid spacing $\Delta x$ (the discretization of the ambient space surrounding the unit sphere). Calculations were done using the algorithm introduced in section \ref{algo-2MBO} and we employ minimizing movements to solve Eq. \eqref{algo-2MBO-heat}. The computation implementation of the MM employs the technique introduced in section \ref{曲面上偏微分方程式の数値計算}. We set $r_0=2/3$, the time step was $h=\Delta x^2/5$, and the threshold length was $\tau=0.03$. 
Let $\bar r(t)$ denote the average radius of the data points that approximate the interface at time $t$.

\noindent
{\bf{Remark:}}
The average radius of the data points approximating the interface is defined as follows. Assume that at time $t$, the interface is computed and represented by $M$ points. For $i=1,2,\cdots,M$, the coordinates of the $i$-th point are denoted as $(x_i,y_i,z_i)$. Then, the average radius $\bar{r}(t)$ is calculated as follows: 
\begin{align}
    \displaystyle{\bar r(t)=\frac{\sum_{i=1}^M \sqrt{x_i^2+y_i^2}}{M}\label{平均半径}}
\end{align}

\begin{figure}[H]
    \begin{center}
        \fbox{\includegraphics[bb=0 0 354 248,scale=0.8]{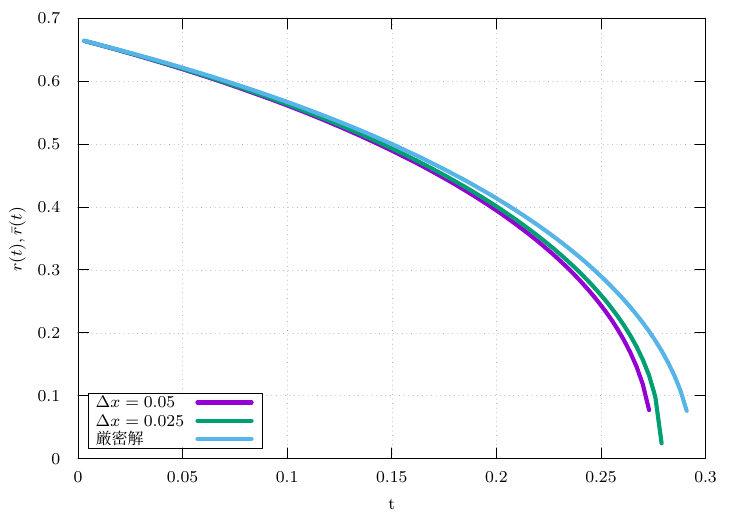}}
        \caption{Relationship between $\Delta x$ and the radius of the circle, with time $t$ on the horizontal axis and the radius on the vertical axis.}\label{res:err_mbos}
    \end{center}
\end{figure}

The results obtained for $\Delta x = 0.05$ and $0.025$ are shown in Figure \ref{res:err_mbos}. The figure shows the exact solution $r(t)$ and the average radius $\bar r(t)$ obtained our method. The points where the curves are interrupted in Figure \ref{res:err_mbos} correspond to the times that the interface could no longer be detected.

The exact solution at $t=0.24$ is approximately $r(0.24)\thickapprox 0.31965745$, while the average radius $\bar r(t)$ obtained from the numerical solution is 0.28291438 for $\Delta x=0.05$ and 0.2955379 for $\Delta x=0.025$.

Although the numerical errors are relatively small at the beginning of the computations, due to the coarsening of the numerical grid for small interfaces, the numerical errors tend to increase as time increases. We also note that the average radius of the numerical solution tends to be smaller than that of the exact solution. Regarding convergence, we indeed observe the tendency of numerical errors to decrease as $\Delta x$ becomes decreases.

When the time step size $h$ is proportional to $\Delta x^2$ and $\tau$ is fixed, reducing $\Delta x$ is expected to improve the accuracy of the  approximation.

\subsection{Numerical error analysis of surface HMBO for two-phase regions}\label{2相領域に対するSurface HMBOアルゴリズムの数値誤差解析}
In this section, show the results of a numerical error analysis for the hyperbolic mean curvature flow on a surface using the algorithm introduced in section \ref{algo-2HMBO}.
Our analysis focuses on the motion of a circle on the surface of the unit sphere.
As shown in Figure \ref{circleonsphere}, let $r(t)$ denote the radius of a circle constrained to the unit sphere.

Since the hyperbolic mean curvature flow represents the motion in which the normal acceleration of the interface is proportional to the mean curvature, the test problem corresponding to Equation \eqref{mcf_surf_de} as follows:
\begin{align}\label{hmcf_surf_de}
    \begin{cases}
        \displaystyle{\frac{d^2r}{dt^2}=\frac{r^2-1}{r}},\quad t>0 \\[8pt]
        \displaystyle{\left.\frac{dr}{dt}\right|_{t=0}=V_0   }     \\[8pt]
        r(0)=r_0
    \end{cases}
\end{align}
Here, $V_0$ is the initial speed of the interface, and $r_0$ is the radius of the circle at the initial time.

We assume that the exact solution represented by the numerical solution obtained by numerically solving Eq. \eqref{hmcf_surf_de}.
We compare this numerical solution with that obtained by our own computational algorithm.
Our computations use DifferentialEquations.jl \cite{des} to numerically solve Eq.n \eqref{hmcf_surf_de}.
Similar to before, we investigated the error's dependence on the spatial discretization $\Delta x$ used in the surrounding space.
We use the algorithm introduced in section \ref{algo-2HMBO} for our numerical calculations.
Minimizing movements were used to solve Eq. \eqref{algo-2HMBO-wave} in the algorithm of section \ref{曲面上偏微分方程式の数値計算}.
Our calculations set $r_0=2/3$, $V_0=0$ and the time step $h$ was assigned to $h=\Delta x/10$. The threshold time $\tau$ was set to $\tau=0.01$.
As before, we define $\bar r(t)$ as the average radius of the data points that approximate the interface at time $t$ using Eq. \eqref{平均半径}.

\begin{figure}[H]
    \begin{center}
        \fbox{\includegraphics[bb=0 0 354 248,scale=0.8]{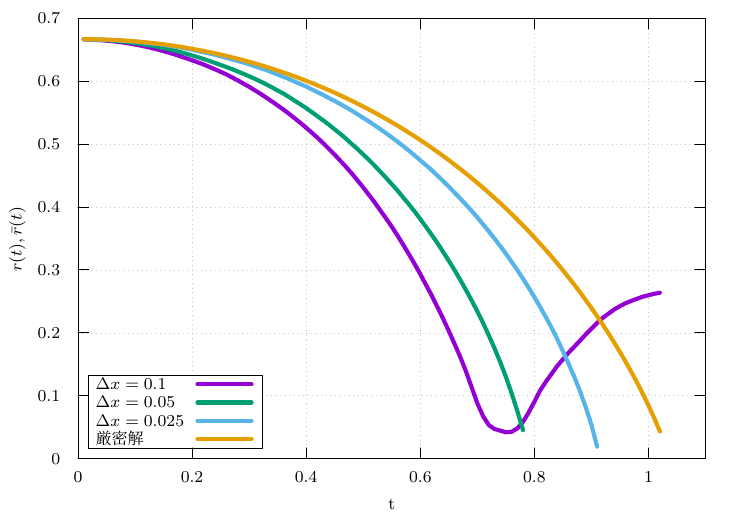}}
        \caption{Relationship between $\Delta x$ and radius, with time $t$ on the horizontal axis and radius on the vertical axis.}\label{res:err_hmbos}
    \end{center}
\end{figure}

We present the numerical results for $\Delta x=0.1, 0.05, 0.025$ in Figure \ref{res:err_hmbos}. The figure shows the numerical solution $r(t)$ obtained by numerically solving equation \eqref{hmcf_surf_de} and the average radius $\bar r(t)$ obtained using our own method. The points where the curves are interrupted in the lower right of Figure \ref{res:err_hmbos} correspond to the time that the interface has dissappeared.

The value of $r(t)$ at $t=0.6$ is approximately $r(0.6)\thickapprox 0.50635371$, while the average radius $\bar r(t)$ obtained from the numerical solution is 0.293178268 for $\Delta x=0.1$, 0.380734809 for $\Delta x=0.05$, and 0.474013006 for $\Delta x=0.025$.

Except for $\Delta x=0.1$, the average radius continues to decrease over time, and the interface can no longer be detected. However, for $\Delta x=0.1$, the average radius starts to increase at some point. This is because the interface that initially shrinks inward eventually becomes a point and starts to expand outward. After the interface has contracted, the subsequent expansion is an interesting feather of the hyperbolic mean curvature flow. A detailed analysis of this phenomenon is a future research topic.

In all cases, the numerical error is relatively small at the beginning of the calculations, but tends to increases over time. Of course, decreasing $\Delta x$ tends to decrease the numerical error.  Setting the time step $h$ proportional to $\Delta x$ and fixing $\tau$ should lead to improved accuracy for decreasing $\Delta x$.

\subsection{Implementation methods}\label{実装方法について}
Here we explain a few important details about the implementation of the numerical algorithms introduced in section \ref{algo-MBO} and section \ref{algo-HMBO}. The numerical algorithms described above involve the calculation of functional values and integrals, which require discretizations a computer.

Similar to section \ref{algo-MBO} and section \ref{algo-HMBO}, we consider a smooth surface $S$ in a three-dimensional Euclidean space and discretize a sufficiently large domain $\Omega_\lambda$ that covers $S$. As before, we refer to the discretized space as $\Omega_\lambda^D$. As a result of the discretization, $\Omega_\lambda^D$ is divided into $N_x$ points in the $x$-direction, $N_y$ points in the $y$-direction, and $N_z$ points in the $z$-direction, as in section \ref{get_dxdydz}. The interval between divisions is assumed to be equal in all three directions and is denoted by $\Delta x$.
\begin{description}
    \item[Approximation of functional values]\mbox{}\\
    In our method, it is necessary to compute the value of the functionals included in the first term of Eq. \eqref{mboF} and Eq.  \eqref{F}. We will explain the approximation using Eq. \eqref{mboF} as an example. The first term of Eq. \eqref{mboF} is expressed by
    \begin{align}\label{mboF_2}
        \mathcal{F}_m(\boldsymbol{w})=\int_{\Omega_\lambda}\left(\frac{|\boldsymbol{w}-\boldsymbol{w}_{m-1}|^2}{2h}+\alpha\frac{|\nabla\boldsymbol{w}|^2}{2}\right)d\boldsymbol{x}.
    \end{align}
    In equation \eqref{mboF_2}, function $\boldsymbol{w}: \Omega_\lambda\rightarrow\mathbb{R}^{K-1}(K:\text{number of phases})$ is a vector-valued and so the functional value can be computed as follows:
    \begin{align}\label{mboF_3}
        \mathcal{F}_m(\boldsymbol{w})=\sum_{i=1}^{K-1}\int_{\Omega_\lambda}\left(\frac{|w^i-w_{m-1}^i|^2}{2h}+\alpha\frac{|\nabla w^i|^2}{2}\right)d\boldsymbol{x}
    \end{align}
    Here, $w^i$ is the $i$th component of $\boldsymbol{w}$. The integrals included in equation \eqref{mboF_3} are approximated using the same technique as in equation \eqref{functional}. For equation \eqref{F}, we used the same method as in the calculation of equation \eqref{functional2}.

    \item[Approximation of $V^i_{\boldsymbol{w}}$]\mbox{}\\
 The value of $V^i_{\boldsymbol{w}}$ appearing in Eq. \eqref{mboF} and Eq. \eqref{F} represents the volume of the extension of phase $i$ within $\Omega_\lambda$. It can be approximated using $H^\epsilon$ and $\phi^i_{\boldsymbol{w}}$ defined in Eq. \eqref{cons_param_mbo} as follows:
    \begin{align*}
        V^i_{\boldsymbol{w}} & =\int_{\Omega_\lambda}H^\epsilon(\phi^i_{\boldsymbol{w}})d\boldsymbol{x},                                                       \\
                             & \approx \Delta x^3 \sum_{\boldsymbol{x}_{i,j,k}\in \Omega_\lambda^D}H^\epsilon(\phi^i_{\boldsymbol{w}}(\boldsymbol{x}_{i,j,k}))
    \end{align*}
    \item[Calculation of $\boldsymbol{w}^S$] \mbox{}\\
        The methods developed here include a step where one must extract the values of $\boldsymbol{w}_{M}$ at the points of $S$:
        $$
            \boldsymbol{w}^S(\boldsymbol{x})=\boldsymbol{w}_M(\boldsymbol{x}),\quad\boldsymbol{x}\in S\cap\Omega_\lambda
        $$
        In the numerical calculations, $\Omega_\lambda$ is discretized. Therefore, an interpolation on $\Omega_\lambda^D$ (the discretized grid points of $\Omega_\lambda$) is required in order to obtain the values on the surface $S$.
        In this study, we have used third-order polynomial interpolations.

    \item[Creation of geodesic distance functions] \mbox{}\\
        When constructing the surface SDVF, signed distance functions on the surface are required. These calculation are not very simple, and we use a method based on the Fast Marching Method \cite{aokishibata} to construct the signed distance vector field on the surface.
\end{description}

\newpage

\bibliography{mainbib}
\bibliographystyle{plain}

\end{document}